\documentclass[a4paper,12pt]{amsart}
\usepackage{amssymb,amscd,amsmath}
\usepackage[dvips]{graphicx}

\theoremstyle{plain}
\newtheorem{theo}{Theorem}[section]
\newtheorem{prop}[theo]{Proposition}
\newtheorem{lemm}[theo]{Lemma}
\newtheorem{coro}[theo]{Corollary}

\theoremstyle{definition}
\newtheorem*{defi}{Definition}
\newtheorem{exam}[theo]{Example}

\theoremstyle{remark}
\newtheorem*{rema}{Remark}

\renewcommand{\u}{u_i^I}
\renewcommand{\Lambda}{C}

\newcommand{\e}{e}
\newcommand{\I}{\mathcal I}
\newcommand{\bI}{\bar \I}
\newcommand{\J}{\mathcal J}
\newcommand{\bJ}{\bar \J}
\newcommand{\uj}{u_j^I}
\newcommand{\field}[1]{\mathbb{#1}}
\newcommand{\C}{\field{C}}

\newcommand{\Q}{\field{Q}}
\newcommand{\R}{\field{R}}
\newcommand{\Z}{\field{Z}}
\newcommand{\PO}{\mathcal{P}}
\newcommand{\DHFP}{\DHF_{\PO}}
\newcommand{\WNP}{\WN_{\PO}}
\newcommand{\augmented}{augmented simplicial set}
\newcommand{\PsiM}{{\rho_M}}
\newcommand{\PsiMi}{{\rho_{M_i}}}
\newcommand{\FL}{\mathcal{F}_L}
\newcommand{\SigmaM}{\Sigma(M)}
\newcommand{\SigmaMsubi}{\Sigma(M)_i}
\newcommand{\SigmasubMi}{\Sigma(M_i)}
\newcommand{\SigmaMsubK}{\Sigma(M)_K}
\newcommand{\SigmasubMK}{\Sigma(M_K)}
\newcommand{\SigmaMsubJ}{\Sigma(M)_J}
\newcommand{\SigmasubMJ}{\Sigma(M_J)}
\newcommand{\LambdaM}{\Lambda(M)}
\newcommand{\wM}{w(M)}
\newcommand{\wMsubJ}{w(M)_J}
\newcommand{\wS}{w(\Stwoalpha)}
\newcommand{\DeltaM}{\Delta(M)}
\newcommand{\DeltaMsubi}{\Delta(M)_i}
\newcommand{\DeltasubMi}{\Delta(M_i)}
\newcommand{\DeltasubMK}{\Delta(M_K)}
\newcommand{\DeltaMsubK}{\Delta(M)_K}

\newcommand{\DeltaMsubJ}{\Delta(M)_J}
\newcommand{\SigmaS}{\Sigma(\Stwoalpha)}
\newcommand{\Stwoalpha}{S^2_\alpha}
\newcommand{\DeltaS}{\Delta(\Stwoalpha)}
\newcommand{\m}{\nu}
\newcommand{\mP}{\nu P}
\newcommand{\SM}{S(M)}
\newcommand{\SMsubi}{S(M)_i}
\newcommand{\SsubMi}{S(M_i)}
\newcommand{\SMsubI}{S(M)_I}

\newcommand{\Si}{S_i}
\newcommand{\tilSi}{\tilde{S}_i}
\newcommand{\tilT}{\tilde{T}}
\newcommand{\tilN}{\tilde{N}}
\newcommand{\tilv}{\tilde{v}}
\newcommand{\tilx}{\tilde{x}}
\newcommand{\Pn}{\field{P}^n}

\renewcommand{\a}{\mathbf a}

\DeclareMathOperator{\Int}{Int} \DeclareMathOperator{\Hom}{Hom}
\DeclareMathOperator{\rank}{rank} \DeclareMathOperator{\sign}{sign}
\DeclareMathOperator{\DHF}{DH} \DeclareMathOperator{\WN}{WN}
\DeclareMathOperator{\HP}{HP} \DeclareMathOperator{\vol}{vol}
\DeclareMathOperator{\Sign}{Sign}
\DeclareMathOperator{\Spin^c}{Spin^c}
\DeclareMathOperator{\Ker}{Ker}
\DeclareMathOperator{\ch}{ch}
\DeclareMathOperator{\Cone}{Cone}

\title{Theory of multi-fans}
\author{Akio Hattori}
\address{Graduate School of Mathematical Science, University of Tokyo,
  Tokyo, Japan}
\email{hattori@ms.u-tokyo.ac.jp}
\author{Mikiya Masuda}
\address{Department of Mathematics, Osaka City University, Osaka, Japan} 
\email{masuda@sci.osaka-cu.ac.jp}

\begin{document}
\bigskip
\bigskip


\maketitle

\section{Introduction} \label{sec:intro}
The purpose of the present paper is to develop a theory of multi-fans
which is an outgrowth of our study initiated in the work \cite{Masuda}
on the topology of torus manifolds (the precise definition will be
given later).
A multi-fan is a combinatorial object generalizing the notion of a
fan in algebraic geometry. Our theory is combinatorial by nature 
but it is built so as to keep a close connection with the topology
of torus manifolds.
\par
It is known that there is a one-to-one
correspondence between toric varieties and fans. A toric variety
is a normal complex algebraic variety of dimension $n$ with a 
$(\C^*)^n$-action 
having a dense orbit. The dense orbit is unique and isomorphic to $(\C^*)^n$, 
and other orbits have smaller dimensions. The fan associated with the toric 
variety is a collection of cones in $\R^n$ with apex at the origin. 
To each orbit there corresponds a cone of dimension equal to the codimension of 
the orbit. Thus the origin is the cone corresponding to the dense orbit, 
one-dimensional cones correspond to maximal singular orbits and so on. 
The important point is the fact that the original toric variety can be
reconstructed from the associated fan, and algebro-geometric properties
of the toric variety can be described in terms of combinatorial data of the 
associated fan.
\par
If one restricts the action of $(\C^*)^n$ to the usual torus 
$T=(S^1)^n$, one can still find the fan, because the orbit types of
the action of the total group $(\C^*)^n$ can be detected by the isotropy 
types of the action of the subgroup $T$. 
Take a circle subgroup
$S$ of $T$ which appears as an isotropy subgroup of the action. Then
each connected component of the closure of the set of those points whose
isotropy subgroup equals $S$ is a $T$-invariant submanifold of real
codimension $2$, and contains a unique $(\C^*)^n$ orbit of
complex codimension $1$. We shall call such a submanifold a
characteristic submanifold. If $M_1,\ldots ,M_k$ are characteristic
submanifolds such that $M_1\cap \cdots \cap M_k$ is non-empty, then
the submanifold $ M_1\cap \cdots \cap M_k$ contains a unique
$(\C^*)^n$-orbit of
complex codimension $k$. This suggests the following
definition of torus manifolds and associated multi-fans.
\par
Let $M$ be an oriented closed manifold of dimension $2n$ with an
effective action of an $n$ dimensional torus $T$ with non-empty fixed
point set $M^T$. A closed, connected, codimension two submanifold of $M$
will be called characteristic if it is a connected component of
the fixed point set of a certain circle subgroup $S$ of $T$, and if it
contains at least one $T$-fixed point. The manifold $M$ together with
a prefered orientation of each characteristic submanifold will be
called a torus manifold. The multi-fan associated with the torus manifold
$M$ involves cones in the Lie
algebra $L(T)$ of $T$, with apex at the origin. 
If $M_i$ is a characteristic submanifold and
$S_i$ is the circle subgroup of $T$ which pointwise fixes $M_i$,
then $S_i$
together with the orientation of $M_i$ determines an element $v_i$
of $\Hom(S^1,T)$, and hence a one dimensional cone in the vector space
$\Hom(S^1,T)\otimes \R$ canonically identified with $L(T)$. If
$M_{i_1},\ldots ,M_{i_k}$ are characteristic submanifolds such that
their intersection contains at least one $T$-fixed point, and if
$v_{i_1},\ldots ,v_{i_k}$ are the corresponding elements in
$\Hom(S^1,T)$, then the $k$-dimensional cone spanned
by $v_{i_1},\ldots ,v_{i_k}$ lies in the multi-fan associated with
$M$. It should be
noted that the intersection of characteristic submanifolds may not
be connected in contrast with the case of toric manifolds where
the intersection is always connected. 
For example, the intersection of
a family of $n$
characteristic submanifolds is a finite set consisting of 
$T$-fixed points. These data are also incorporated in the definition 
of the associated multi-fan in Section 2.
\par
One of the differences between a fan and a multi-fan is that, while
cones in a fan intersect only at their faces and their union covers
the space $L(T)$ just once without overlap for complete toric varieties,
it happens that the union of cones in a multi-fan covers $L(T)$ with
overlap
for torus manifolds. Also the same multi-fan
corresponds to different torus manifolds. Nevertheless
it turns out that
important topological invariants of a torus manifold can be described in 
terms of the associated multi-fan. In fact it is furthermore possible
to develop an abstract theory of multi-fans and to define various
``topological'' invariants of a multi-fan in such a way that, when the
multi-fan
is associated with a torus manifold, they coincide with the ordinary
topological invariants of the manifold. For example, the
``multiplicity of overlap'', which we call
the degree of the multi-fan, equals the Todd genus for a unitary torus 
manifold (unitary toric manifold in the terminology in \cite{Masuda};
the precise definition will be given in Section 9). 
\par
Another feature of the theory of toric varieties is the correspondence
between ample line bundles over a complete toric variety and convex 
polytopes. From a topological point of view this can be explained in the 
following way.
Let $(M,\omega)$ be a compact symplectic manifold with a Hamiltonian 
$T$-action, and let $\Psi :M\to L(T)^*$ be an associated moment map.
Then it is well-known (\cite{Atiyah}, \cite{GS})
that the image $P$ of $\Psi$ is a convex polytope. Moreover, if the de Rham 
cohomology class of $\omega$ is an integral class, then the polytope $P$
is a lattice polytope up to translations in $L(T)^*$ identified with $\R^n$.
Delzant \cite{Delzant} showed that the original symplectic manifold 
$(M,\omega)$ is equivariantly symplectomorphic to a complete non-singular
toric variety and the form $\omega$ is transformed into the first Chern form of
an ample line bundle $L$ over the toric variety. It is known that the number
of lattice points in $P$ is equal to the Riemann-Roch number
\[ \int_M e^{c_1(L)}{\mathcal{T}}(M) \] 
where $ \mathcal{T}$ is the Todd class of $M$, see e.g. \cite{Fulton}. This
sort of phenomenon was generalized to ``presymplectic'' toric manifolds
by Karshon and Tolman \cite{KT}, then to $\Spin^c$ toric manifolds by
Grossberg and Karshon \cite{GK} and also to unitary toric manifolds by the
second-named author \cite{Masuda} in the form which relates the equivariant
index
of the line bundle $L$ regarded as an element of $K_{T}(M)$ to the
Duistermaat-Heckman measure of the moment map associated with $L$. In these
extended cases the form $\omega$ may be degenerate or the line bundle
may not be ample, and consequently the image of the moment map may not
be convex any longer. This leads us to consider more general figures
which we call
multi-polytopes. A multi-polytope is a pair of a multi-fan and an
arrangement of affine hyperplanes in $L(T)^*$.
A similar notion was introduced by Karshon and Tolman
\cite{KT} and also by Khovanskii and Pukhlikov
\cite{KP} for
ordinary fans under the name twisted polytope and virtual polytope
respectively. We shall develop a combinatorial theory of
multi-polytopes as well; we define the
Duistermaat-Heckman measure and the equivariant index
in a purely combinatorial
fashion for multi-fans and multi-polytopes, and generalize the above
results in the combinatorial context. Also we shall introduce a
combinatorial counterpart of
a moment map which can be used to interpret the combinatorial
Duistermaat-Heckman measure.
\par
In carrying out the above program, the use of equivariant homology and
cohomology plays an important role. First note that the group
$\Hom(S^1,T)$ can be canonically identified with the equivariant
integral homology group $H_2(BT)$, and hence the vector space
$L(T)$ with $H_2(BT,\R)$. In this way we regard vectors $v_i$ in
a multi-fan as lying in $H_2(BT,\R)$.
On the other hand a characteristic submanifold $M_i$ with a fixed
orientation determines a
cohomology class $\xi_i$ in $H_{T}^2(M)$, the equivariant Poincar\'e
dual of $M_i$. These cohomology classes are fundamental for describing
the first Chern class of an equivariant line bundle over $M$. This fact
enables us to associate a multi-polytope and a generalization of
the Duistermaat-Heckman
measure with an equivariant line bundle. To a $T$-line bundle $L$
whose equivariant first Chern class has the form
$ c_1^T(L)=\sum c_i\xi_i $,
we associate an arrangement of affine hyperplanes $F_i$
in $H^2(BT;\R)=L(T)^*$ defined by
\[ F_i=\{u\in H^2(BT;\R) \mid \langle u,v_i\rangle =c_i \}. \]
This arrangement defines the multi-polytope associated with the line bundle
$L$. Moreover it is possible to define
the equivariant cohomology of a complete simplicial multi-fan and extend
the results to such abstract multi-fans and multi-polytopes.
\par
If $v_{i_1},\ldots ,v_{i_n}$ are primitive vectors generating an
$n$-dimensional cone in the multi-fan associated with a torus
manifold, then they form a basis of $\Hom(S^1,T)$. However, in the
definition of abstract multi-fans, this condition is not 
postulated. From this point of view, it is natural to deal with torus 
orbifolds besides torus manifolds. This can be achieved without much 
change technically.
More importantly every complete
simplicial multi-fan (the precise definition will be given later)
can be realized as a multi-fan associated with a torus orbifold
in dimensions greater than 2. In dimensions 1 and 2,  realizable 
multi-fans are characterized.

Concerning the realization problem we are not sure at this moment
whether every non-singular complete simplicial multi-fan is realized
as the multi-fan associated with a torus manifold. In any case it 
should be noticed that a multi-fan may correspond to more than one 
torus manifolds unlike the case of toric varieties. 

We now explain the contents of each section. In Section \ref{sec:multi-fan}
we give a
definition of a multi-fan and introduce certain related notions. The
completeness of multi-fans is most important. It is a generalization
of the notion of completeness of fans. But the definition takes a somewhat
sophisticated form. Section \ref{sec:T_y_genus} is devoted to the
$T_y$-genus of  a 
complete multi-fan. It is defined in such a way that, when the multi-fan
is associated with a unitary torus manifold $M$, it coincides
with the $T_y$-genus of $M$. In Lemma \ref{lemm:h=e} we exhibit an equality
which is an analogue of the relation between $h$-vectors and $f$-vectors
in combinatorics (see e.g. \cite{Stanley}), and which, we hope, sheds more
insight on that relation. 
\par
In Sections \ref{sec:multi-polytope} and \ref{sec:DH} the notion of
a multi-polytope and the associated Duistermaat-Heckman function are
defined. As explained above, a multi-polytope is a pair
$\PO =(\Delta, \mathcal{F})$ of an $n$-dimensional complete multi-fan
$\Delta$ and an arrangement of hyperplanes $\mathcal{F}=\{F_i\}$
in $H^2(BT;\R)$
with the same index set
as the set of $1$-dimensional cones in $\Delta$.
It is called simple
if the multi-fan $\Delta$ is simplicial. The Duistermaat-Heckman function
$\DHFP$ associated with a simple multi-polytope $\PO$ is a locally
constant integer-valued function with bounded support
defined on the complement of the hyperplanes
$\{F_i\}$. The wall crossing formula 
(Lemma \ref{lemm:transition_formula_for_DH})
which describes the difference of the values of the function on adjacent
components is important for later use.
In Section \ref{sec:WN} another locally constant function on the
complement of
the hyperplanes $\{F_i\}$ in a multi-polytope $\PO$,
called the winding number, is introduced. It
satisfies a wall crossing formula entirely similar to the
Duistermaat-Heckman function. When
the multi-fan $\Delta$ is associated with a torus manifold
or a torus orbifold $M$ and if there
is an equivariant complex line bundle $L$ over $M$, then there is a simple
multi-polytope $\PO$ naturally associated with $L$, and the winding number
$\WNP$
is closely related to the moment map of $L$. In fact it can be regarded
as the density function of the Duistermaat-Heckman measure associated
with the
moment map. Theorem \ref{theo:DH=WN}, the main theorem
in Section \ref{sec:WN}, states that the Duistermaat-Heckman function
and the winding number coincide for any simple multi-polytope.        
\par
Section \ref{sec:Ehrhart_polynomial} is devoted to a generalization
of the Ehrhart polynomial
to multi-polytopes. If $P$ is a convex lattice polytope and if
$\nu P$ denotes the multiplied polytope by a positive integer
$\nu$, then the number of lattice
points $\sharp(\nu P)$ contained in $\nu P$ is developed as a polynomial
in $\nu$ . It is called the Ehrhart polynomial of $P$.
The generalization to multi-polytopes is straightforward and properties
similar to that of the ordinary Ehrhart polynomial hold
(Theorem \ref{theo:Ehrhart_polynomial_for_PO}). If $\PO$ is a simple
lattice multi-polytope, 
then the associated Ehrhart polynomial $\sharp(\nu\PO)$ is
defined by
\[ \sharp(\nu\PO)=\sum_{u\in H^2(BT;\Z)}\DHF_{{\nu\PO}_+}(u), \]
where ${\PO}_+$ denotes a multi-polytope obtained from $\PO$ by a small
enlargement. Lemma \ref{lemm:key_lemma} is crucial
for the proof of Theorem \ref{theo:Ehrhart_polynomial_for_PO} and for
the later development of the theory. Its corollary,
Corollary \ref{coro:key_corollary}, gives a localization
formula for the Laurent polynomial
$\sum_{u\in H^2(BT;\Z)}\DHF_{{\PO}_+}(u)t^u$ regarded as a character
of $T$. It can be considered as a combinatorial
generalization of Theorem \ref{theo:multiplicity}. It reduces to
$\sharp\PO$ when evaluated at the identity. Using this
fact, in Section \ref{sec:cohomological_formula},
a cohomological formula expressing $\sharp\PO$ in terms of the
``Todd class'' of the multi-fan and the first ``Chern class'' of the
multi-polytope is given in Theorem \ref{theo:integral_formula}. 
The formula can be thought of as a
generalization of the formula expressing the number of lattice points
in a convex lattice polytope by the Riemann-Roch number of the
corresponding ample line bundle. The argument is completely combinatorial.
We define the equivariant cohomology $H_T^*(\Delta)$ of a multi-fan $\Delta$
which is a module over $H^*(BT)$, the index map (Gysin homomorphism)
$\pi_! :H_T^*(\Delta) \to H^{*-2n}(BT)$,
the cohomology $H^*(\Delta)$ of $\Delta$ and finally the evaluation on
the ``fundamental class''.
As a corollary a generalization of Khovanskii-Pukhlikov formula (\cite
{KP}) for simple lattice multi-polytopes is given in Theorem
\ref{theo:KP_formula}.
\par
In Section \ref{sec:torus_manifold} it is shown
how to associate a multi-fan with a torus manifold. It is also shown that
the associated multi-fan is complete. Then, in
Section \ref{sec:T_y_genus_of_a_torus_manifold},
the $T_y$-genus of a general torus manifold is defined and is proved to
coincide with the $T_y$-genus of the associated multi-fan
in Theorem \ref{theo:T_y_genus of a torus manifold}. As a corollary
a formula for the signature of a torus manifold is given. In the same spirit
the definition of the equivariant index of a line bundle over a general torus
manifold is given in Section \ref{sec:equivariant_index}
using a localization formula which holds
in the case of unitary torus manifolds. The main
theorem of this section, Theorem \ref{theo:multiplicity},
gives a formula describing
that equivariant index using the winding number. It generalizes
the results of \cite{KT}, \cite{GK} and \cite{Masuda} as indicated before.
Results of Section \ref{sec:DH} and \ref{sec:WN} are crucially used
here. 
\par
In Section \ref{sec:orbifolds} necessary changes to deal with torus
orbifolds are explained briefly.
One of the remarkable points is that the torus action and the
orbifold structure are closely related to each other for a torus orbifold as
is explained in Lemma \ref{lemm:isotropy}. In
the last section the realization problem is dealt with. Main results
of the section are Theorems~\ref{theo:2-dim torus orbifold}, 
\ref{theo:4-dim torus orbifold} and \ref{theo:higher dim torus orbifold}.

\bigskip
\section{Multi-fans} \label{sec:multi-fan}

In \cite{Masuda}, we introduced the notion of 
a unitary toric manifold, which contains a compact non-singular toric variety 
as an example, and associated with it a combinatorial object 
called a multi-fan, which is a more general notion than a complete 
non-singular fan.  
In this section, we define a multi-fan in a combinatorial way and 
in full generality.  The reader will find that our notion of a multi-fan 
is a complete generalization of a fan.  We also define the completeness 
and non-singularity of a multi-fan, which generalize the corresponding 
notion of a fan.  
To do this, we begin with reviewing the definition of a fan.  

Let $N$ be a lattice of rank $n$, which is isomorphic to $\Z^n$.  We denote 
the real vector space $N\otimes \R$ by $N_{\R}$.  A subset $\sigma$ of $N_\R$ 
is called a \emph{strongly convex rational polyhedral cone} (with apex at the 
origin) if there exits a finite number of vectors $v_1,\dots,v_m$ in $N$ 
such that 
\[ \sigma=\{ r_1v_1+\dots +r_mv_m\mid r_i\in\R \text{ and } r_i\ge 0 
\text{ for all $i$}\}\]
and $\sigma\cap (-\sigma)=\{ 0\}$.  Here \lq\lq rational" means that it is 
generated by vectors in the lattice $N$, and \lq\lq strong" convexity that 
it contains no line through the origin.  We will often call a strongly convex 
rational polyhedral cone in $N_\R$ simply a \emph{cone} in $N$.  The dimension 
$\dim\sigma$ of a cone $\sigma$ is the dimension of the linear space spanned 
by vectors in $\sigma$.  A subset $\tau$ of $\sigma$ is called a \emph{face} 
of $\sigma$ if there is a linear function $\ell\colon 
N_\R\to \R$ such that $\ell$ takes 
nonnegative values on $\sigma$ and $\tau=\ell^{-1}(0)\cap \sigma$.  A cone is 
regarded as a face of itself, while others are called \emph{proper} faces. 

\begin{defi} A fan $\Delta$ in $N$ is a set of a finite number of strongly 
convex rational polyhedral cones in $N_\R$ such that 
\begin{enumerate}
\item Each face of a cone in $\Delta$ is also a cone in $\Delta$;
\item The intersection of two cones in $\Delta$ is a face of each. 
\end{enumerate}
\end{defi}

\begin{defi} A fan $\Delta$ is said to be \emph{complete} if the union of 
cones in $\Delta$ covers the entire space $N_\R$. 
\end{defi}

A cone is called \emph{simplicial} if it is generated 
by linearly independent vectors.  If the generating vectors can be taken as 
a part of a basis of $N$, then the cone is called \emph{non-singular}. 

\begin{defi} A fan $\Delta$ is said to be \emph{simplicial} (resp. 
\emph{non-singular}) if every cone in $\Delta$ is simplicial (resp. 
non-singular).  
\end{defi}

The basic theory of toric varieties tells us that a fan is complete 
(resp. simplicial or non-singular) if and only if the corresponding 
toric variety is compact (resp. an orbifold or non-singular).  

For each $\tau\in\Delta$, we define $N^\tau$ to be the quotient lattice of $N$ 
by the sublattice generated (as a group) by $\tau\cap N$; so the rank of 
$N^\tau$ is $n-\dim\tau$.  We consider cones 
in $\Delta$ that contain $\tau$ as a face, and project them on $(N^\tau)_\R$.  
These projected cones form a fan in $N^\tau$, which we denote by 
$\Delta_\tau$ and call the 
\emph{projected fan} with respect to $\tau$.  
The dimensions of the projected cones decrease by 
$\dim\tau$. The completeness, simpliciality and non-singularity of $\Delta$ 
are inherited to $\Delta_\tau$ for any $\tau$.  

We now generalize these notions of a fan.  Let $N$ be as before.  Denote 
by $\Cone(N)$ the set of all cones in $N$.  
An ordinary fan is a subset of $\Cone(N)$.  The set $\Cone(N)$ has 
a (strict) partial ordering $\prec$ defined by: $\tau\prec\sigma$ 
if and only if $\tau$ is a proper face of $\sigma$.  
The cone $\{0\}$ consisting of the origin is the unique minimum element 
in $\Cone(N)$. 
On the other hand, let $\Sigma$ be a partially ordered finite set with 
a unique minimum element.  We denote the (strict) partial ordering by 
$<$ and the minimum element by $*$.  
An example of $\Sigma$ used later is an abstract simplicial set 
with an empty set added as a member, which we call an 
\emph{\augmented}.  In this case the partial ordering 
is defined by the inclusion relation and the empty set is the unique minimum 
element which may be considered as a $(-1)$-simplex.  
Suppose that there is a map 
\[ \Lambda\colon \Sigma \to \Cone(N)\]
such that 
\begin{enumerate}
\item $\Lambda(*)=\{ 0\}$;
\item If $I<J$ for $I,J\in\Sigma$, then $\Lambda(I)\prec \Lambda(J)$;
\item For any $J\in\Sigma$ the map $C$ restricted on 
$\{I\in \Sigma\mid I\leq J\}$ is an isomorphism of ordered sets onto
$\{\kappa\in \Cone(N)\mid \kappa \preceq \Lambda(J)\}$.  
\end{enumerate}
For an integer $m$ such that $0\le m\le n$, we set 
\[ \Sigma^{(m)}:=\{ I\in \Sigma\mid \dim\Lambda(I)=m\}.\]
One can easily check that $\Sigma^{(m)}$ does not depend on $\Lambda$. 
When $\Sigma$ is an \augmented, $I\in\Sigma$ belongs to $\Sigma^{(m)}$ if and 
only if the cardinality $|I|$ of $I$ is $m$, namely $I$ is an $(m-1)$-simplex. 
Therefore, even if $\Sigma$ is not an \augmented, we use the notation 
$|I|$ for $m$ when $I\in \Sigma^{(m)}$. 

The image $\Lambda(\Sigma)$ is a finite set of cones in $N$.  We may think of 
a pair $(\Sigma,\Lambda)$ as a set of cones in $N$ labeled by 
the ordered set $\Sigma$.  
Cones in an ordinary fan intersect only at their faces, but cones in  
$\Lambda(\Sigma)$ may overlap, even the same cone may appear repeatedly 
with different labels.  
The pair $(\Sigma,\Lambda)$ is almost what we call a multi-fan, but we 
incorporate a pair of weight functions on cones in $\Lambda(\Sigma)$ of the 
highest dimension $n=\rank N$.  More precisely, we consider two functions 
\[ w^\pm\colon \Sigma^{(n)}
\to\Z_{\ge 0}.\]
We assume that $w^+(I)>0$ or
$w^-(I)>0$ for every $I \in \Sigma^{(n)}$.
These two functions have its origin from geometry. In fact if 
$M$ is a torus manifold of dimension $2n$ and if $M_{i_1},\ldots ,M_{i_n}$ 
are characteristic submanifolds such that
their intersection contains at least one $T$-fixed point, then the 
intersection $M_I=\bigcap_\nu M_{i_\nu}$ consists of a finite number of
$T$-fixed points. At each fixed point $p\in M_I$ the tangent space $\tau_p$
has two orientations; one is endowed by the orientation of $M$ and the
other comes from the intersection of the oriented submanifolds $M_{i_\nu}$.
Denoting the ratio of the above two orientations by $\epsilon_p$ we define
the number $w^+(I)$ to be the number of points $p\in M_I$ with $\epsilon_p=+1$
and similarly for $w^-(I)$. More detailed explanation will be given in
Section 9.    

\begin{defi} We call a triple $\Delta:=(\Sigma,\Lambda,w^\pm)$ a 
\emph{multi-fan} in $N$.  We define the dimension of $\Delta$ to be 
the rank of $N$ (or the dimension of $N_\R$). 
\end{defi}

Since an ordinary fan $\Delta$ in $N$ is a subset of $\Cone(N)$, one can view 
it as a multi-fan by taking $\Sigma=\Delta$, 
$\Lambda=\text{the inclusion map}$, 
$w^+=1$, and $w^-=0$.  In a similar way as in the case of ordinary fans, 
we say that a multi-fan 
$\Delta=(\Sigma,\Lambda,w^\pm)$ is \emph{simplicial} 
(resp. \emph{non-singular}) 
if every cone in $\Lambda(\Sigma)$ is simplicial (resp. non-singular).  
The following lemma is easy.  

\begin{lemm} \label{lemm:simplicial}
A multi-fan $\Delta=(\Sigma,\Lambda,w^\pm)$ is simplicial if and only if 
$\Sigma$ is isomorphic to an \augmented\ as partially 
ordered sets. 
\end{lemm}

The definition of completeness of a multi-fan $\Delta$ is rather complicated.
A naive definition of the completeness would be that the union 
of cones in $\Lambda(\Sigma)$ covers the entire space $N_\R$.  
But it turns out that 
this is not a right definition if we look at multi-fans associated 
with unitary torus manifolds, see Section~\ref{sec:torus_manifold}.  
Although the two weighted functions $w^\pm$ are incorporated in the 
definition of a multi-fan, only the difference 
\[ w:=w^+-w^-\]
matters in this paper except Section~\ref{sect:realization}.  
We shall introduce the following intermediate notion of pre-completeness 
at first.  A vector $v\in N_\R$ will be called generic if $v$ does
not lie on any linear subspace spanned by a cone in 
$\Lambda(\Sigma)$ of dimesnsion less than $n$. For a generic
vector $v$ we set $d_v=\sum_{v\in \Lambda(I)}w(I)$, where
the sum is understood to be zero if there is no such $I$.

\begin{defi}
We call a multi-fan $\Delta=(\Sigma,\Lambda,w^\pm)$ of dimension $n$ 
{\it pre-complete} if $\Sigma^{(n)}\not=\emptyset$ and 
the integer $d_v$ is independent of
the choice of generic vectors $v$. We call this integer
the {\it degree} of $\Delta$ and denote it by $\deg(\Delta)$.
\end{defi}

\begin{rema} For an ordinary fan, pre-completeness is the same as completeness.
\end{rema}

To define the completeness for a multi-fan $\Delta$, 
we need to define a projected multi-fan with respect to an element in 
$\Sigma$.  We do it as follows. For each $K\in \Sigma$, we set 
\[ \Sigma_K:=\{ J\in\Sigma\mid K\leq J\}.\]
It inherits the partial ordering from $\Sigma$, and $K$ is the unique 
minimum element in $\Sigma_K$.
A map 
\[ \Lambda_K\colon \Sigma_K\to \Cone(N^{\Lambda(K)})\]
sending $J\in\Sigma_K$ to the cone $\Lambda(J)$ projected on 
$(N^{\Lambda(K)})_\R$ satisfies the three properties above 
required for $\Lambda$. 
We define two functions 
\[ {w_K}^\pm\colon  \Sigma_K^{(n-|K|)}\subset \Sigma^{(n)} \to \Z_{\ge 0}\]
to be the restrictions of $w^\pm$ to $\Sigma_K^{(n-|K|)}$.  The triple 
$\Delta_K:=(\Sigma_K,\Lambda_K,{w_K}^\pm)$ is a multi-fan in 
$N^{\Lambda(K)}$, and this is the desired 
\emph{projected multi-fan} with respect to $K\in \Sigma$. 
When $\Delta$ is an ordinary fan, 
this definition agrees with the previous one.  

\begin{defi} A pre-complete multi-fan $\Delta=(\Sigma,\Lambda,w^\pm)$ is 
said to be 
\emph{complete} if the projected multi-fan $\Delta_K$ is pre-complete for any 
$K\in \Sigma$. 
\end{defi}

\begin{rema} A multi-fan $\Delta$ is complete 
if and only if the projected multi-fan $\Delta_J$ is pre-complete for any 
$J\in\Sigma^{(n-1)}$.  The argument is as follows.  The pre-completeness of 
$\Delta_J$ for $J\in\Sigma^{(n-1)}$ implies that 
$d_v=\sum_{v\in\Lambda(I)}w(I)$ remains unchanged when $v$ gets across the 
codimension one cone $\Lambda(J)$, which means the pre-completeness of 
$\Delta$.  Since $\Sigma_K^{(n-|K|-1)}$ is contained in $\Sigma^{(n-1)}$ 
for any $K\in\Sigma$, the pre-completeness of $\Delta_J$ for any $J\in 
\Sigma^{(n-1)}$ also implies the pre-completeness of $\Delta_K$ for any 
$K\in\Sigma$. 
\end{rema}


\begin{exam} \label{exam:degree_two}
Here is an example of a complete non-singular multi-fan of 
degree two. Let $v_1,\dots,v_5$ be integral vectors shown in 
Figure~\ref{1}, where the dots denote lattice points. 


\begin{figure}[h]
\begin{center}
\setlength{\unitlength}{1mm}
\begin{picture}(50,40)(0,0)
\put(10,10){\includegraphics[scale=.3]{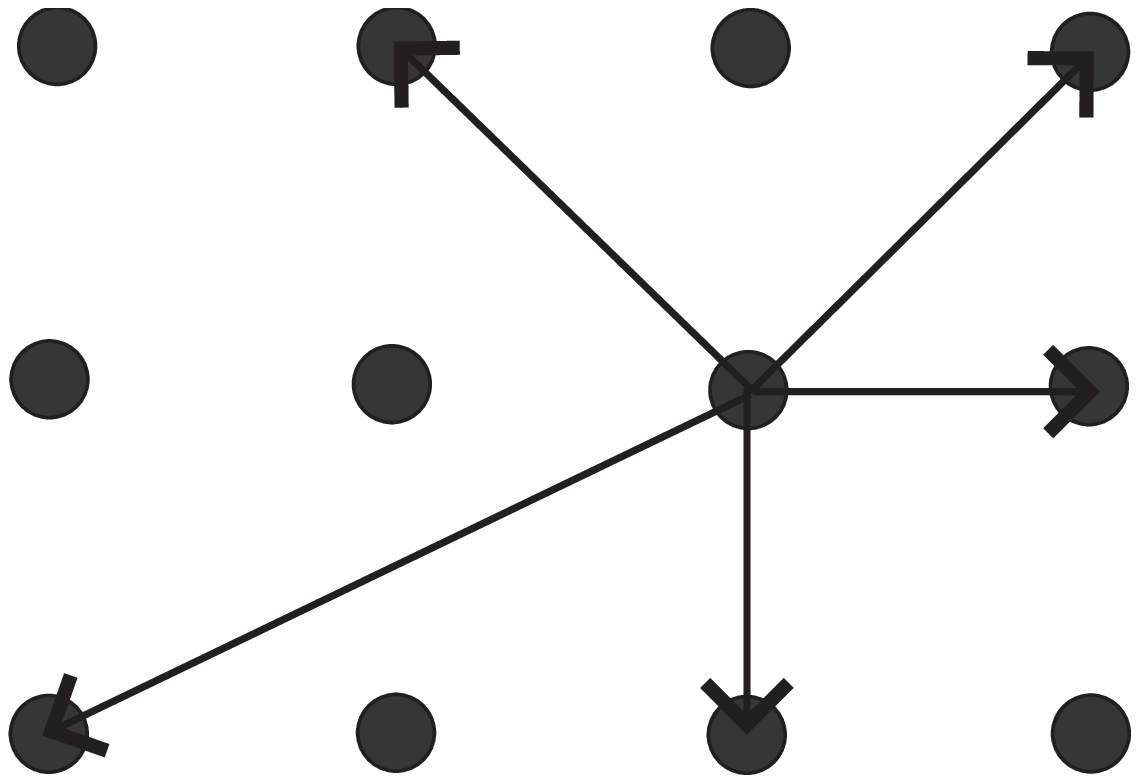}}
\put(47,21){$v_1$}
\put(47,36){$v_4$}
\put(20,36){$v_2$}
\put(7,6){$v_5$}
\put(30,6){$v_3$}
\end{picture}
\end{center}
\caption{}\label{1}
\end{figure}

\bigskip

\noindent
The vectors are rotating around the origin twice in counterclockwise. 
We take 
\[ \Sigma=\{\phi, \{1\},\dots,\{5\},\{1,2\},\{2,3\},\{3,4\},
\{4,5\},\{5,1\}\},\]
define $\Lambda\colon \Sigma\to \Cone(N)$ by 
\begin{align*}
\Lambda(\{i\})&=\text{ the cone spanned by $v_i$,}\\
\Lambda(\{i,i+1\})&=\text{ the cone spanned by $v_i$ and $v_{i+1}$,}
\end{align*}
where $i=1,\dots,5$ and $6$ is understood to be $1$, and take 
$w^\pm$ such that $w=1$ on every two dimensional cone.  Then 
$\Delta=(\Sigma,\Lambda,w^\pm)$ is a complete non-singular two-dimensional 
multi-fan with $\deg(\Delta)=2$.  
\end{exam}

\begin{exam} Here is an example of a complete multi-fan \lq\lq with folds". 
Let $v_1,\dots,v_5$ be vectors shown in Figure~\ref{1_1}.  


\bigskip

\begin{figure}[h]
\setlength{\unitlength}{1mm}
\begin{picture}(40,40)(0,-5)
\put(0,0){\includegraphics[scale=.3]{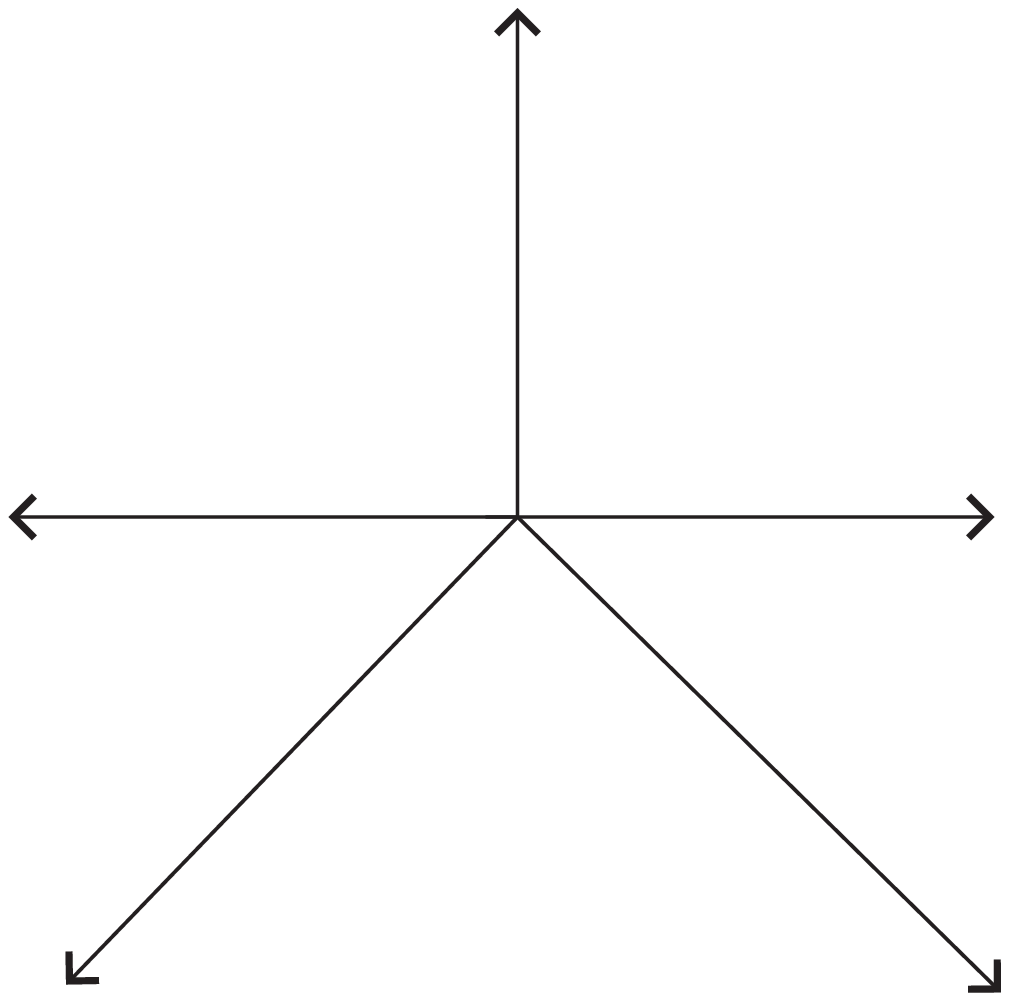}}
\put(15,33){$v_2$}
\put(33,15){$v_1$}
\put(-5,15){$v_4$}
\put(-3,-2){$v_3$}
\put(33,-2){$v_5$}
\end{picture}
\caption{}\label{1_1}
\end{figure}

\bigskip

\noindent
We take the same $\Sigma$ and $\Lambda$ as in 
Example~\ref{exam:degree_two} and take $w^\pm$ 
such that 
\[
w(\{3,4\})=-1 \quad\text{and}\quad w(\{i,i+1\})=1\ \text{ for $i\not=3$.} 
\]
Then $\Delta=(\Sigma,\Lambda,w^\pm)$ is a complete two-dimensional 
multi-fan with $\deg(\Delta)=1$. 

A similar example can be constructed for a number of vectors 
$v_1,\dots,v_d$ $(d\ge 3)$ by defining 
\begin{align*}
w(\{i,i+1\})&=1\quad\text{if $v_i$ and $v_{i+1}$ are rotating in 
counterclockwise,}\\
w(\{i,i+1\})&=-1\quad\text{if $v_i$ and $v_{i+1}$ are rotating in 
clockwise,}
\end{align*}
where $d+1$ is understood to be $1$.  The degree 
$\deg(\Delta)$ is the rotation 
number of the vectors $v_1,\dots,v_d$ around the origin in counterclockwise 
and may not be one.  
\end{exam}

\begin{exam} Here is an example of a multi-fan which is pre-complete but 
not complete. Let $v_1,\dots,v_5$ be vectors shown in Figure~\ref{2}. 

\bigskip

\begin{figure}[h]
\setlength{\unitlength}{1mm}
\begin{picture}(55,45)(0,3)
\put(10,10){\includegraphics[scale=.3]{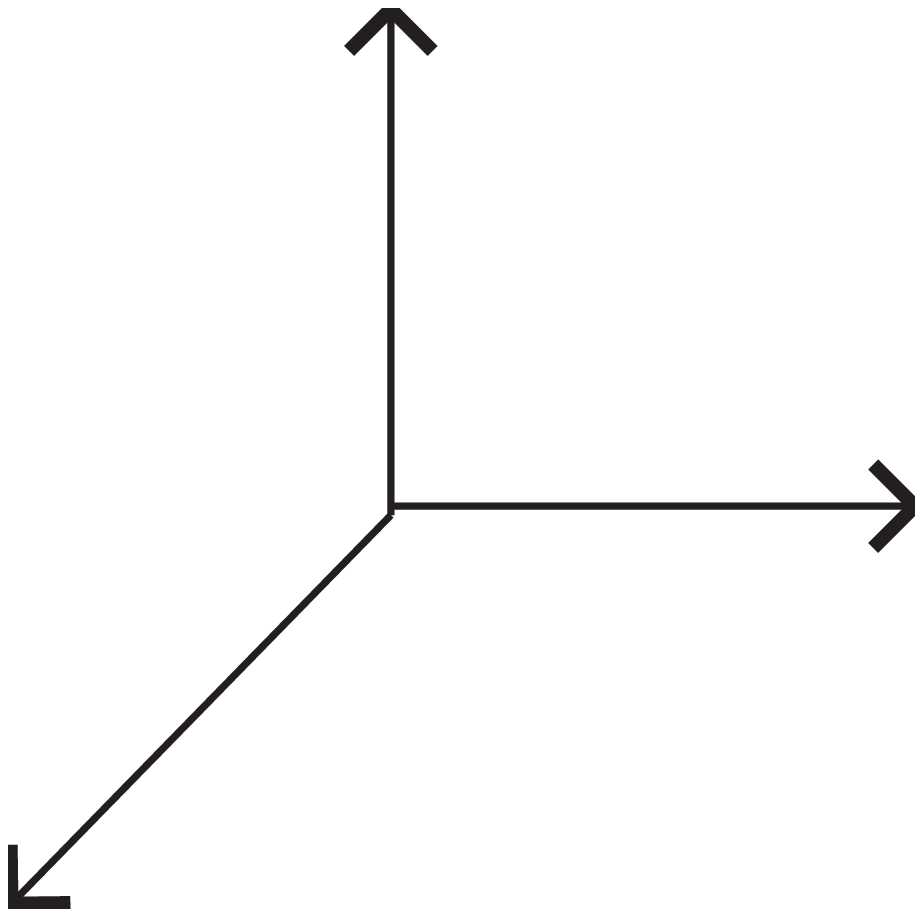}}
\put(5,7){$v_3$}
\put(40,21){$v_1=v_4$}
\put(18,41){$v_2=v_5$}
\end{picture}
\caption{}\label{2}
\end{figure}

\bigskip

\noindent
We take 
\[ \Sigma=\{\phi,\{1\},\dots,\{5\},\{1,2\},\{2,3\},\{3,1\},\{4,5\}\},\]
define $\Lambda\colon \Sigma\to \Cone(N)$ as in Example~\ref{exam:degree_two},  
and take $w^\pm$ such that 
\[ w(\{1,2\})=2,\ w(\{2,3\})=1,\ w(\{3,1\})=1,\ w(\{4,5\})=-1.\]
Then 
$\Delta=(\Sigma,\Lambda,w^\pm)$ is a two-dimensional 
multi-fan which is pre-complete (in fact, $\deg(\Delta)=1$) 
but not complete because the projected multi-fan $\Delta_{\{i\}}$ for 
$i\not=3$ is not pre-complete.  
\end{exam}

So far, we treated \emph{rational} cones that are generated by vectors in 
the lattice $N$.  But, most of the notions introduced above make sense even 
if we allow cones generated by vectors in $N_\R$ which may not be in $N$. 
In fact, the notion of non-singularity requires the lattice $N$, but others do 
not. Therefore, one can define a multi-fan and its completeness and 
simpliciality in this extended category as well.  The reader will find that 
the arguments developed in Sections~\ref{sec:T_y_genus} through 
\ref{sec:WN} work in this extended category. 

\bigskip

\section{$T_y$-genus of a multi-fan} \label{sec:T_y_genus}

A unitary torus manifold $M$ determines 
a complete non-singular multi-fan. 
(This will be discussed and extended to 
torus manifolds in Section~\ref{sec:torus_manifold}.) 
On the other hand, the $T_y$-genus (also called $\chi_y$-genus) for unitary manifolds 
introduced by Hirzebruch in his famous book \cite{Hirzebruch} is defined 
for $M$. Its characteristic power series is given by
$\dfrac{x(1+ye^{-(1+y)x})}{1-e^{-(1+y)x}}$. 
It is a polynomial in one variable $y$ of degree (at most) 
$\frac{1}{2}\dim M$.  
The Kosniowski formula about the $T_y$-genus for unitary 
$S^1$-manifolds (see \cite{HT}, 
\cite{Kawakubo}) and the results in \cite{Masuda} imply that 
the $T_y$-genus of $M$ should be described in terms of the multi-fan 
associated with $M$. 
In this section (and in Section~\ref{sec:T_y_genus_of_a_torus_manifold})
we give the explicit description. In fact, our argument is rather 
more general.  We think of the $T_y$-genus of $M$ as a polynomial invariant 
of the associated multi-fan which is complete and non-singular.  
It turns out that the polynomial invariant can be 
defined not only for the 
multi-fans associated with unitary torus manifolds but 
also for all complete simplicial multi-fans.

Since the lattice $N$ is unnecessary from now until the end of 
Section~\ref{sec:WN}, 
we shall denote the vector space, in which cones sit, by $V$ instead of 
$N_\R$. 
Let $\Delta=(\Sigma,\Lambda,w^\pm)$ be a complete simplicial multi-fan 
defined on $V$.  By Lemma~\ref{lemm:simplicial} we may assume that 
$\Sigma$ is an \augmented, say, consisting of subsets of $\{1, \dots,d\}$ 
and $\Sigma^{(1)}=\{\{1\},\dots,\{d\}\}$ where $d$ is the number of elements 
in $\Sigma^{(1)}$.  For each $i=1,\dots,d$, let $v_i$ denote a nonzero 
vector in the one-dimensional cone $\Lambda(\{i\})$.  Choose a generic 
vector $v\in V$.  Let $I\in \Sigma^{(n)}$.  Since $v_i$'s ($i\in I$) 
are linearly 
independent, $v$ has a unique expression $\sum_{i\in I}a_iv_i$ with 
real numbers $a_i$'s.  The coefficients $a_i$'s are all nonzero because 
$v$ is generic.  We set 
\[ \mu(I):=\sharp\{ i\in I\mid a_i>0\}.
\]
This depends on $v$ although $v$ is not recorded in the notation $\mu(I)$. 

\begin{defi} For an integer $q$ with $0\le q\le n$, we define 
\[ h_q(\Delta):=\sum_{\mu(I)=q}w(I)\quad\text{and}\quad 
\e_q(\Delta):=\sum_{K\in\Sigma^{(q)}}\deg(\Delta_K).
\]
Note that $h_n(\Delta)=\deg(\Delta)=\e_0(\Delta)$, and $\e_q(\Delta)$'s 
are independent of $v$.  If $\Delta$ is a complete simplicial 
multi-fan such that $\deg(\Delta)=1$ and $w(I)=1$ for all $I\in \Sigma^{(n)}$ 
(e.g. this is the case if $\Delta$ is a complete simplicial ordinary fan), 
then $\deg(\Delta_K)$ equals $1$ for all $K\in\Sigma$ and hence $\e_q(\Delta)$ 
agrees with the number of cones of dimension $q$ in the multi-fan.  
\end{defi}

The following lemma reminds us of the relation between the $h$-vectors 
and the $f$-vectors for simplicial sets studied in combinatorics 
(see \cite{Stanley}).  

\begin{lemm} \label{lemm:h=e}
$\displaystyle{\sum_{q=0}^n h_q(\Delta)(s+1)^q=\sum_{m=0}^n 
\e_{n-m}(\Delta)s^m}$ where $s$ is an indeterminate.
\end{lemm}

\begin{proof} The lemma is equivalent to the following equality:
\begin{equation*} 
\sum_{q=m}^n h_q(\Delta)\binom{q}{m}=\e_{n-m}(\Delta).\tag{3.1}
\end{equation*}
It follows from the definition of $h_q(\Delta)$ that 
\begin{equation*}
\text{l.h.s. of (3.1)}=\sum_{q=m}^n\binom{q}{m}\sum_{\mu(I)=q}
w(I).\tag{3.2}
\end{equation*}
On the other hand, we shall rewrite $\e_{n-m}(\Delta)$.  
It follows from the definition of $\deg(\Delta_K)$ that 
\[ \deg(\Delta_K)=\sum_{J\in\Sigma_K^{(n-|K|)} s.t.\ 
v_K\in\Lambda_K(J)}w_K(J)
\]
where $v_K$ denotes the projection image of $v$ on the quotient vector 
space of $V$ by the subspace $V_K$ spanned by the cone $\Lambda(K)$. 
Note that $v_K$ lies in $\Lambda_K(J)$ if and only if $v$ lies in 
$\Lambda(J\cup K)$ modulo $V_K$, and that $w_K(J)=w(J\cup K)$ by definition. 
Therefore, writing $J\cup K$ as $I$, the equality above turns into 
\[ \deg(\Delta_K)=\sum_{I}w(I),
\]
where $I$ runs over elements in $\Sigma^{(n)}$ such that $K\subset I$ 
and $v\in \Lambda(I)$ modulo $V_K$. 
Putting this in the defining equation of $\e_{n-m}(\Delta)$, we have 
\begin{equation*}
 \e_{n-m}(\Delta)=\sum_{K,I}w(I), \tag{3.3}
\end{equation*}
where the sum is taken over elements $K\in \Sigma^{(n-m)}$ and 
$I\in\Sigma^{(n)}$ such that $K\subset I$ and 
$v\in \Lambda(I)$ modulo $V_K$.  Fix $I\in\Sigma^{(n)}$ with $\mu(I)=q$, 
and observe how many times $I$ appears in the above sum. 
It is equal to the number of $K\in\Sigma^{(n-m)}$ 
such that $K\subset I$ and $v\in \Lambda(I)$ modulo $V_K$.  But the 
number of such $K$ is $\binom{q}{m}$.  To see this, we note that 
$\mu(I)=q$ means that $\sharp\{i\in I\mid a_i>0\}=q$ by definition,
where $v=\sum_{i\in I}a_iv_i$, and that 
the condition that $v\in \Lambda(I)$ modulo $V_K$ is equivalent to 
saying that $K$ contains the complement of the set $\{i\in I\mid a_i>0\}$ 
in $I$.  Therefore, any such $K$ is obtained as the complement of 
a subset of $\{i\in I\mid a_i>0\}$ with cardinality $m$, so that 
the number of such $K$ is $\binom{q}{m}$. 
This together with (3.2) and (3.3) proves the equality (3.1). 
\end{proof}

\begin{coro} \begin{enumerate}
\item $h_q(\Delta)$'s are independent of the choice of the generic 
vector $v$.
\item $h_q(\Delta)=h_{n-q}(\Delta)$ for any $q$. 
\end{enumerate}
\end{coro}

\begin{proof} (1) This immediately follows from Lemma~\ref{lemm:h=e} 
because $\e_q(\Delta)$'s are independent of $v$. 

(2) If we take $-v$ instead of $v$, then $\mu(I)$ turns into 
$n-\mu(I)$, so that $h_q(\Delta)$ turns into $h_{n-q}(\Delta)$. 
Since $h_q(\Delta)$'s are independent of $v$ as shown in (1) above, 
this proves $h_q(\Delta)=h_{n-q}(\Delta)$. 
\end{proof}

When $\Delta$ is associated with a unitary torus manifold $M$, 
the $T_y$-genus of $M$ turns out to be given by 
$\sum_{q=0}^n h_q(\Delta)(-y)^q$.  (This will be discussed in 
Section~\ref{sec:T_y_genus_of_a_torus_manifold} later.)  
Motivated by this observation, 

\begin{defi} For a complete simplicial multi-fan $\Delta$, we define 
\[ T_y[\Delta]:=\sum_{q=0}^n h_q(\Delta)(-y)^q \]
and call it the \emph{$T_y$-genus} of $\Delta$.  Note that 
$T_0[\Delta]=h_0(\Delta)=h_n(\Delta)=\deg(\Delta)$. 
\end{defi}

Lemma~\ref{lemm:h=e} can be restated as 

\begin{coro} \label{coro:h=e}
Let $\Delta$ be a complete simplicial multi-fan. Then 
\[ T_y[\Delta]=\sum_{m=0}^n\e_{n-m}(\Delta)(-1-y)^m.\]
\end{coro}


\bigskip

\section{Multi-polytopes} \label{sec:multi-polytope}

\begin{figure}[h]
\setlength{\unitlength}{1mm}
\begin{picture}(60,55)(-5,-5)
\put(0,0){\includegraphics[scale=.3]{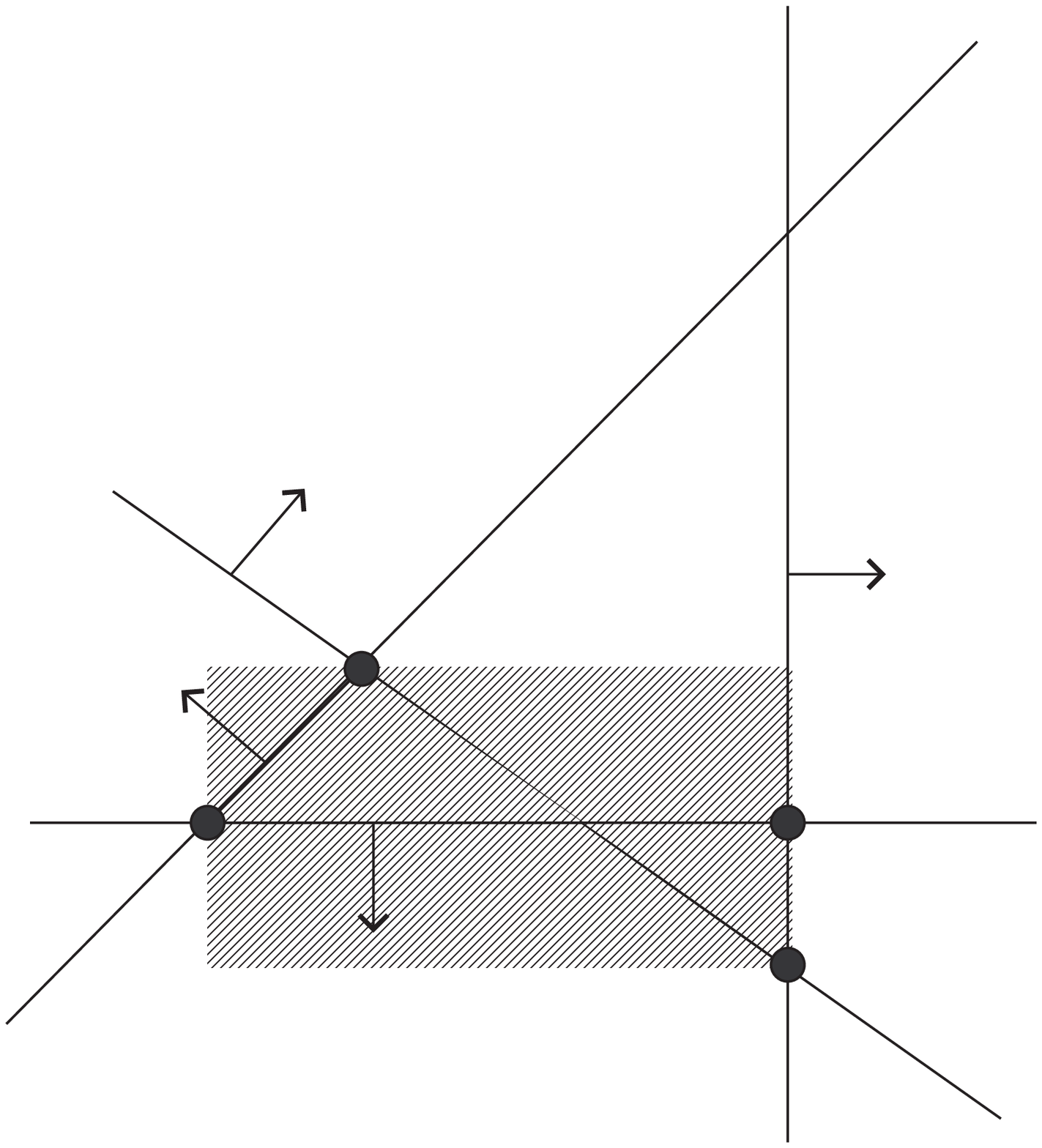}}
\put(3,22){$v_2$}
\put(15,30){$v_1$}
\put(42,25){$v_4$}
\put(15,5){$v_3$}
\put(1,30){$\ell_1$}
\put(47,13){$\ell_3$}
\put(-3,2){$\ell_2$}
\put(30,47){$\ell_4$}
\put(22,-5){$Q$}
\end{picture}
\hspace{1cm}
\setlength{\unitlength}{1mm}
\begin{picture}(60,55)(-5,-5)
\put(0,0){\includegraphics[scale=.3]{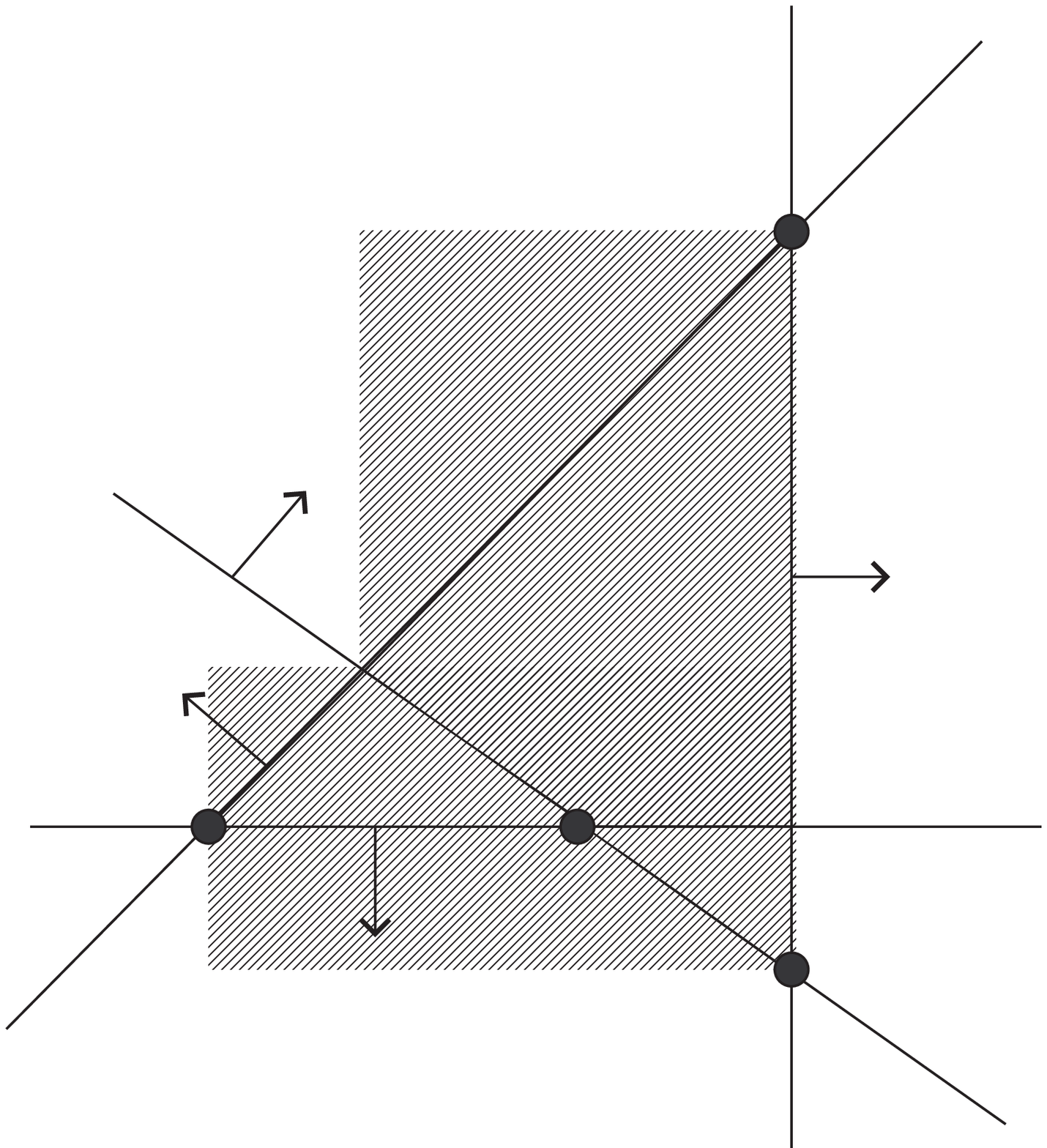}}
\put(3,22){$v_2$}
\put(15,30){$v_1$}
\put(42,25){$v_4$}
\put(15,5){$v_3$}
\put(1,30){$\ell_1$}
\put(47,13){$\ell_3$}
\put(-3,2){$\ell_2$}
\put(30,47){$\ell_4$}
\put(22,-5){$Q'$}
\end{picture}
\caption{}\label{4.2}
\end{figure}

A convex polytope $P$ in $V^*=\Hom(V,\R)$ 
is the convex hull of a finite set of points in $V^*$.  It is  
the intersection of a finite number of half spaces in $V^*$ separated 
by affine hyperplanes, 
so there are a finite number of nonzero vectors $v_1,\dots,
v_d$ in $V$ and real numbers $c_1,\dots,c_d$ such that 
\[ P=\{ u\in V^* \mid \langle u,v_i\rangle \le c_i\text{ for all $i$}
\}, \]
where $\langle\ ,\ \rangle$ denotes the natural pairing between $V^*$ and 
$V$.  (Warning: In this paper, we take $v_i$ to be \lq\lq outward normal" 
to the corresponding face of $P$ contrary to the usual convention in 
algebraic geometry, cf. e.g. \cite{Oda}.) 
The convex polytope $P$ can be recovered from the data 
$\{(v_i,c_i)\mid i=1,\dots,d\}$.  But, a more general figure like 
$Q$ shaded in Figure~\ref{4.2} cannot be determined by the data 
$\{(v_i,c_i)\mid i=1,\dots,d\}$.  
We need to prescribe the vertices of 
$Q$, in other words, which pairs of lines $\ell_i$'s are presumed to 
intersect.  
For instance, if four points $\ell_1\cap \ell_2$, 
$\ell_2\cap \ell_3$, $\ell_3\cap \ell_4$ and $\ell_4\cap\ell_1$ are 
presumed to be vertices (and the others such as $\ell_2\cap\ell_4$ 
are not), then we can find the figure $Q$ in Figure~\ref{4.2}.  
But, if different four points $\ell_1\cap \ell_4$, 
$\ell_4\cap \ell_2$, $\ell_2\cap \ell_3$ and $\ell_3\cap\ell_1$ are 
presumed to be vertices, then we obtain a figure $Q'$ shaded in 
Figure~\ref{4.2}.  

\begin{figure}[h]
\setlength{\unitlength}{1mm}
\begin{picture}(55,60)(0,-10)
\put(5,2.5){\includegraphics[scale=.5]{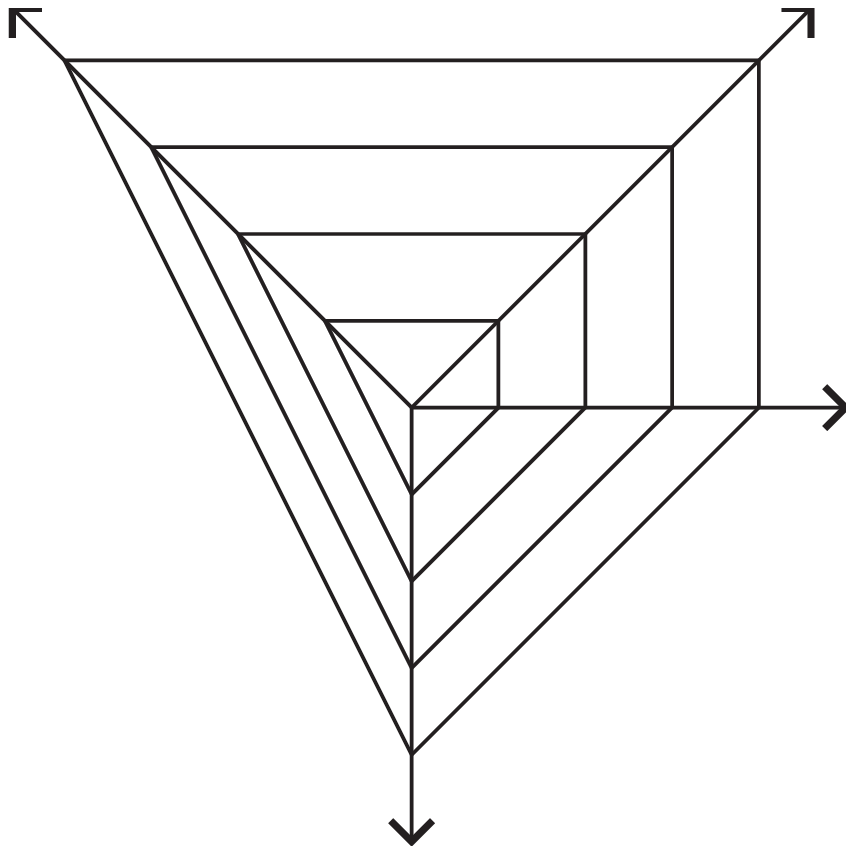}}
\put(5,47){$v_2$}
\put(45,47){$v_1$}
\put(50,25){$v_4$}
\put(26,-1){$v_3$}
\put(22,-7){(1)}
\end{picture}
\hspace{1cm}
\setlength{\unitlength}{1mm}
\begin{picture}(55,60)(0,-10)
\put(5,2.5){\includegraphics[scale=.5]{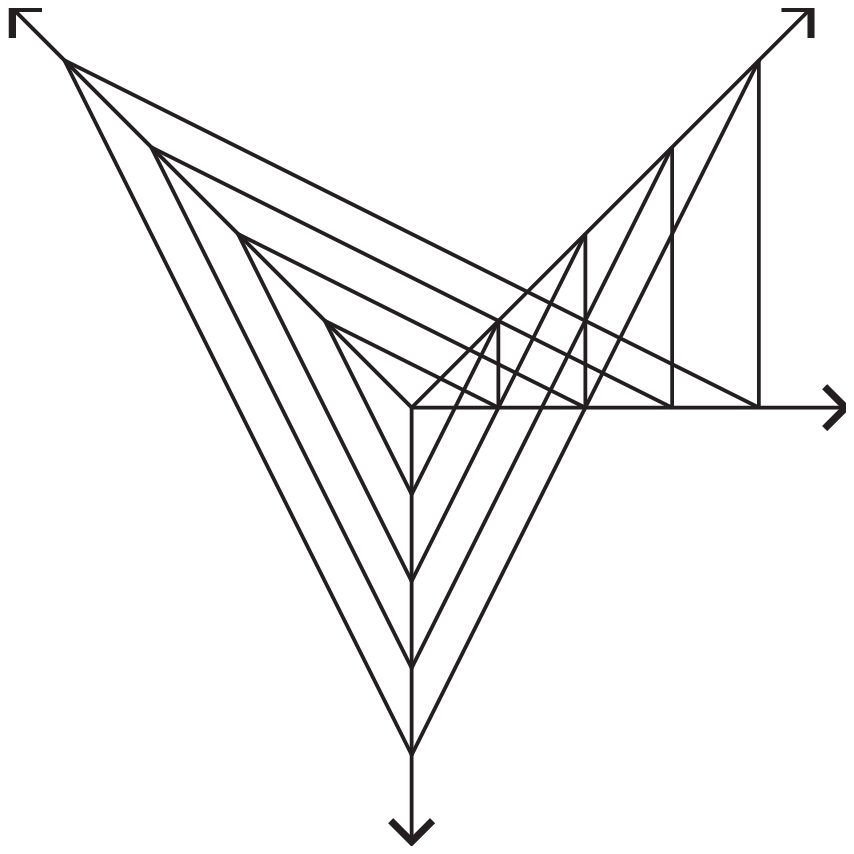}}
\put(5,47){$v_2$}
\put(45,47){$v_1$}
\put(50,25){$v_4$}
\put(26,-1){$v_3$}
\put(22,-7){(2)}
\end{picture}
\caption{}\label{4.3}
\end{figure}

The data of whether two lines $\ell_i$ and $\ell_j$ are 
presumed to intersect is equivalent to the data of whether the 
corresponding vectors $v_i$ and $v_j$ span a cone.  In the former 
(resp. latter) 
example above, resulting cones are four two-dimensional ones shown in 
Figure~\ref{4.3} (1) (resp. (2)). 
Needless to say, $\ell_i$ is \lq perpendicular' to the 
half line spanned by $v_i$.  

A polytope gives rise to a multi-fan in this way.  One notes that a 
convex polytope gives rise to a complete fan.  
Taking this observation into account, we reverse a gear.  We start with a 
complete multi-fan $\Delta=(\Sigma,\Lambda,w^\pm)$. 
Let $\HP(V^*)$ be the set of all affine hyperplanes in $V^*$. 

\begin{defi} Let $\Delta=(\Sigma,\Lambda,w^\pm)$ be a complete multi-fan and 
let $\mathcal F\colon \Sigma^{(1)}\to \HP(V^*)$ be a map such that 
the affine hyperplane 
$\mathcal F(I)$ is \lq perpendicular' to the half line $\Lambda(I)$ 
for each $I\in\Sigma^{(1)}$, i.e., an element in $\Lambda(I)$ 
takes a constant on $\mathcal F(I)$.  We call a pair
$(\Delta,\mathcal F)$ a \emph{multi-polytope} and denote it by $\PO$. 
The dimension of a multi-polytope $\PO$ is defined to be the dimension 
of the multi-fan $\Delta$. 
We say that a multi-polytope $\PO$ 
is \emph{simple} if $\Delta$ is simplicial. 
\end{defi}

\begin{rema} The completeness assumption for $\Delta$ is not needed for
the definition of multi-polytopes. We incorporated it because most of our results
depend on that assumption. Similar notions were introduced by Karshon-Tolman \cite{KT} and 
Khovanskii-Pukhlikov \cite{KP} when $\Delta$ is an ordinary fan. 
They use the terminology \emph{twisted polytope} and \emph{virtual 
polytope} respectively.  The notion of multi-polytopes is a 
direct generalization of that of twisted polytopes. 
The relation between virtual polytopes and multi-polytopes 
is clarified by \cite{Nishimura}.
\end{rema}

\begin{exam} \label{exam:polytope}
A convex polytope determines a complete fan together with an 
arrangement of affine hyperplanes containing the facets 
of the polytope (as explained above), so it uniquely 
determines a multi-polytope.  
\end{exam}

\begin{exam} \label{exam:star}
Associated with the multi-fan in Example~\ref{exam:degree_two}, 
one obtains the arrangement of lines drawn in Figure~\ref{star} with a 
suitable choice of the map $\mathcal F$. 
The pentagon shown up in Figure~\ref{star} produces 
the same arrangement of lines and can be viewed as a multi-polytope 
as explained in Example~\ref{exam:polytope} above, 
but these two multi-polytopes 
are different because the underlying multi-fans are different; one is 
a multi-fan of degree two while the other is an ordinary fan. 
The reader will find 
a star-shaped figure in the former multi-polytope.  

\bigskip
\bigskip

\begin{figure}[h]
\setlength{\unitlength}{1mm}
\begin{picture}(55,70)(0,0)
\put(0,0){\includegraphics[scale=.3]{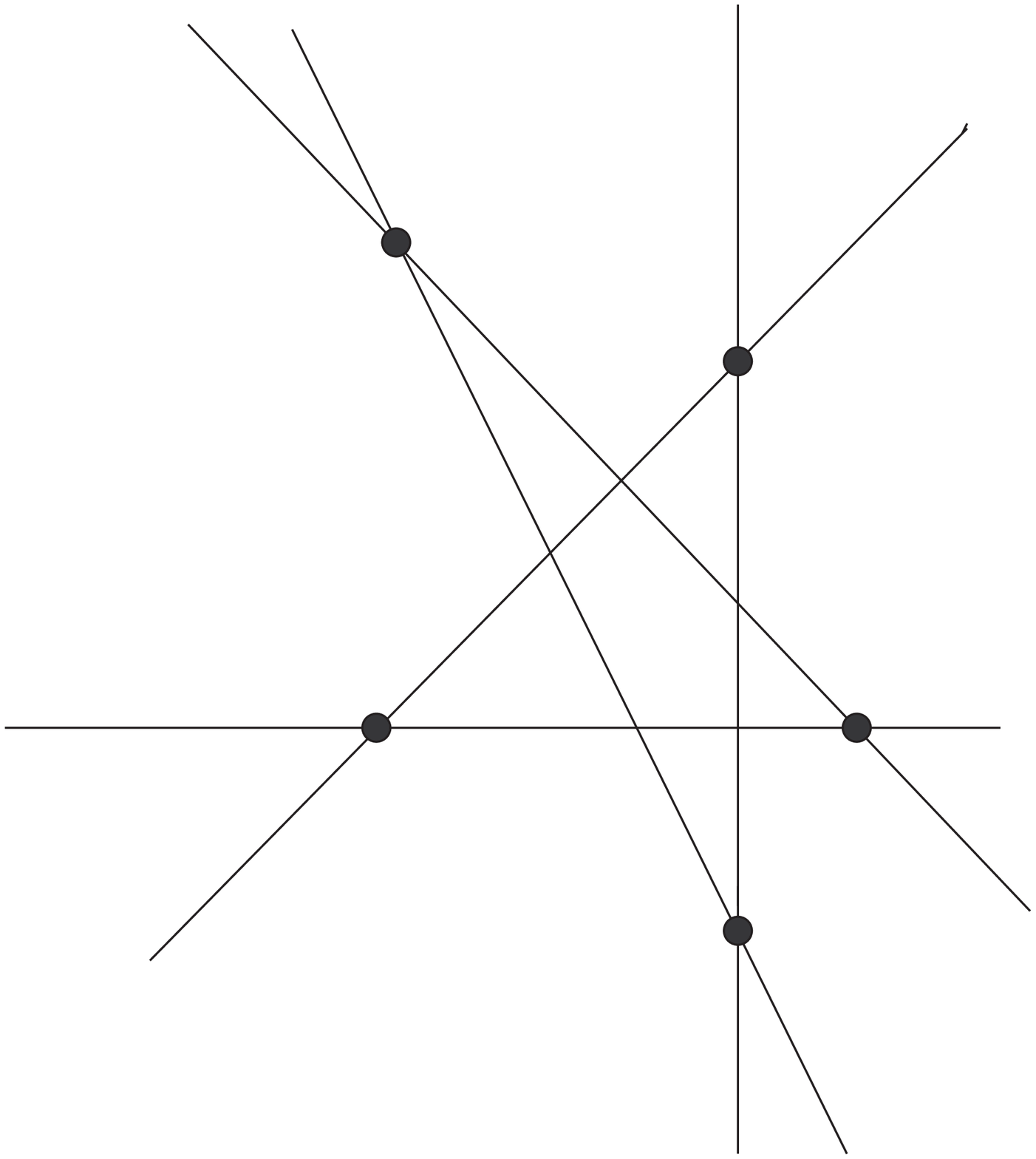}}
\put(37,65){$\mathcal F(\{1\})$}
\put(2,5){$\mathcal F(\{2\})$}
\put(-5,26){$\mathcal F(\{3\})$}
\put(50,8){$\mathcal F(\{4\})$}
\put(45,-5){$\mathcal F(\{5\})$}
\end{picture}
\caption{}\label{star}
\end{figure}

\end{exam}

\bigskip

\section{Duistermaat-Heckman functions} \label{sec:DH}

A multi-polytope $\PO=(\Delta,\mathcal F)$ defines an arrangement of 
affine hyperplanes in $V^*$.  In this section, 
we associate with $\PO$ a function on $V^*$ minus the affine hyperplanes 
when $\PO$ is simple.  
This function is locally constant and Guillemin-Lerman-Sternberg 
formula (\cite{GLS1} \cite{GLS2}) tells us that it agrees with the 
density function of a Duistermaat-Heckman measure when 
$\PO$ arises from a moment map.  

Hereafter our multi-polytope $\PO$ 
is assumed to be simple, so that the multi-fan 
$\Delta=(\Sigma,\Lambda,w^\pm)$ is complete and simplicial unless 
otherwise stated.  
As before, we may assume that $\Sigma$ consists of subsets of 
$\{1,\dots,d\}$ and $\Sigma^{(1)}=\{\{1\},\dots,\{ d\}\}$, and 
denote by $v_i$ a nonzero vector in the one-dimensional cone 
$\Lambda(\{i\})$.  To simplify notation, we denote 
$\mathcal F(\{i\})$ by $F_i$ and set 
\[ F_I:=\cap_{i\in I}F_i\quad\text{for $I\in\Sigma$}.\]
$F_I$ is an affine space of dimension $n-|I|$. 
In particular, if $|I|=n$ (i.e., $I\in\Sigma^{(n)}$), then $F_I$ 
is a point, denoted by $u_I$.   

Suppose $I\in\Sigma^{(n)}$.  Then the set $\{ v_i\mid i\in I\}$ forms a 
basis of $V$. 
Denote its dual basis of $V^*$ by 
$\{\u\mid i\in I\}$, i.e., $\langle \u,v_j\rangle=
\delta_{ij}$ where $\delta_{ij}$ denotes the Kronecker delta. 
Take a generic vector $v\in V$. Then $\langle \u,v\rangle
\not=0$ for all $I\in\Sigma^{(n)}$ and $i\in I$. Set 
\[
(-1)^I:=(-1)^{\sharp\{ i\in I \mid \langle \u,v\rangle>0\}}
\quad\text{and}\quad 
(\u)^+:=\begin{cases} \u\quad&\text{if $\langle \u,v
\rangle > 0$}\\
-\u\quad&\text{if $\langle \u,v\rangle < 0$.}\end{cases}
\]
We denote by $\Lambda^*(I)^+$ the cone 
in $V^*$ spanned by $(\u)^+$'s $(i\in I)$ with apex at $u_I$, and 
by $\phi_I$ its characteristic function. 

\begin{defi} We define a function $\DHFP$ on 
$V^*\backslash\cup_{i=1}^d F_i$ by 
\[ \DHFP:=\sum_{I\in\Sigma^{(n)}}(-1)^Iw(I)\phi_I
\]
and call it the \emph{Duistermaat-Heckman function} associated with 
$\PO$.  
\end{defi}

\begin{rema}
Apparently, the function $\DHFP$ is defined on 
the whole space $V^*$ and 
depends on the choice of the generic vector $v\in V$, but 
we will see in Lemma~\ref{lemm:independence} below that 
it is independent of $v$ on $V^*\backslash \cup F_i$. 
This is the reason why we restricted the domain of the function to 
$V^*\backslash \cup F_i$. 
\end{rema}

For the moment, we shall see the independence of $v$ when 
$\dim \PO=1$. 

\begin{exam} \label{exam:one_dimension}
Suppose $\dim \PO=1$.  We identify $V$ with $\R$, so that $V^*$ is also 
identified with $\R$.  
Let $E$ be the subset of $\{1,\dots,d\}$ such that $i\in E$ 
if and only if $\Lambda(\{i\})$ is the half line consisting of nonnegative 
real numbers.  Then 
the completeness of $\Delta$ means that 
\begin{equation}
\sum_{i\in E}w(\{i\})=\sum_{i\notin E}w(\{i\})=\deg(\Delta).
\tag{5.1}
\end{equation}
Take a nonzero vector $v$. 
Since $V^*$ is identified with $\R$, each affine hyperplane $F_i$ is 
nothing but a real number. 
Suppose that $v$ is toward the positive direction.  Then 
\begin{equation}
(-1)^{\{i\}}=\begin{cases} -1\quad&\text{if $i\in E$}\\
1 \quad&\text{if $i\notin E$}\end{cases}\tag{5.2}
\end{equation}
and the support of the characteristic function $\phi_{\{i\}}$ is 
the half line given by 
\[ \{ u\in \R\mid  F_i\le u \}.\]
Therefore 
\begin{equation}
\DHFP(u)=
\sum_{i\in E\ s.t.\ F_i<u}(-w(\{i\}))+
\sum_{i\notin E\ s.t.\ F_i<u}w(\{i\})\tag{5.3}
\end{equation}
for $u\in \R\backslash\cup F_i$. If $u$ is sufficiently small, then the sum 
above is empty; so it is zero.  
If $u$ is sufficiently large, 
then the sum is also zero by (5.1).  Hence the support of the function 
$\DHFP$ is bounded.  

Now, suppose that $v$ is toward the negative direction.  Then 
$(-1)^{\{i\}}$ above is multiplied by $-1$ 
and the inequality $\le$ above turns into $\ge$.  Therefore 
\begin{equation}
\DHFP(u)=\sum_{i\in E\ s.t.\ u<F_i}w(\{i\})+
\sum_{i\notin E\ s.t.\ u<F_i}(-w(\{i\})).\tag{5.4}
\end{equation}
It follows that 
\[ \text{r.h.s. of (5.3)}-\text{r.h.s. of (5.4)}=-\sum_{i\in E}w(\{i\})
+\sum_{i\notin E}w(\{i\}),\]
which is zero by (5.1).  This shows that the function $\DHFP$ 
is independent of $v$ when $\dim \PO=1$. 
\end{exam}

\begin{exam} 
For the star-shaped multi-polytope in Example~\ref{exam:star}, 
$\DHFP$ takes $2$ on the pentagon, $1$ on the five triangles adjacent 
to the pentagon and $0$ on other (unbounded) regions. 
The check is left to the reader. 
\end{exam}

Assume $n=\dim \Delta >1$. 
For each $\{i\}\in\Sigma^{(1)}$, the projected multi-fan 
$\Delta_{\{i\}}=(\Sigma_{\{i\}},\Lambda_{\{i\}},
w_{\{i\}}^\pm)$, 
which we abbreviate as $\Delta_i=(\Sigma_i,\Lambda_i,w_i^\pm)$, 
is defined on the quotient vector space $V/V_i$ of $V$ by 
the one-dimensional subspace $V_i$ spanned by $v_i$.  
Since $\Delta$ is complete and simplicial, so is $\Delta_i$. 
We identify the dual space $(V/V_i)^*$  with 
\[ (V^*)_i:=\{ u\in V^*\mid \langle u,v_i\rangle=0\}\]
in a natural way.  We choose an element $f_i\in F_i$ arbitrarily and 
translate $F_i$ onto $(V^*)_i$ by $-f_i$.  If $\{i,j\}\in 
\Sigma^{(2)}$, then $F_j$ intersects $F_i$ and their intersection 
will be translated into $(V^*)_i$  by $-f_i$. 
This observation leads us to consider the map 
\[ \mathcal F_i\colon \Sigma_i\to \HP((V^*)_i)\]
sending $\{j\}\in \Sigma_i^{(1)}$ to $F_i\cap F_j$ translated by $-f_i$. 
The pair $\PO_i=(\Delta_i,\mathcal F_i)$ is a multi-polytope in 
$(V/V_i)^*\cong (V^*)_i$.  

Let $I\in \Sigma^{(n)}$ such that $i\in I$.  
Since $\langle \uj,v_i\rangle=\delta_{ij}$, $\uj$ for $j\not=i$ is an element 
of $(V^*)_i$, which we also regard as an element of 
$(V/V_i)^*$ through the isomorphism $(V/V_i)^*\cong (V^*)_i$.  
We denote the projection image of the generic element 
$v\in V$ on $V/V_i$ by $\bar v$. Then we have 
$\langle \bar v,\uj\rangle=\langle v,\uj\rangle$ 
for $j\not=i$, where $\uj$ at the left-hand side is viewed as an 
element of $(V/V_i)^*$ while the one at the right-hand side is viewed 
as an element of $(V^*)_i$.  Since 
$\langle \bar v,\uj\rangle=\langle v,\uj\rangle\not=0$ 
for $j\not=i$, we use $\bar v$ to define $\DHF_{\PO_i}$.  

\begin{lemm} \label{lemm:transition_formula_for_DH}
{\rm (Wall crossing formula.)}
Let $F$ be one of $F_i$'s.  Let $u_\alpha$ and $u_\beta$ 
be elements in $V^*\backslash\cup_{i=1}^d F_i$ such that the segment 
from $u_\alpha$ to $u_\beta$ intersects the wall $F$ transversely at 
$\mu$, and does not intersect any other $F_j\not=F$.  Then 
\[ \DHFP(u_\alpha)-\DHFP(u_\beta)
=\sum_{i:F_i=F}\sign\langle u_\beta-u_\alpha,v_i\rangle 
\DHF_{\PO_i}(\mu-f_i).
\]
\end{lemm}

\begin{proof} For simplicity we assume that there is only one $i$ such that 
$F_i=F$.  We may assume that $\langle u_\beta-u_\alpha,v_i\rangle$ is 
positive without loss of generality.  
The situation is as in Figure~\ref{4.4}.

\begin{figure}[h]
\setlength{\unitlength}{1mm}
\begin{picture}(55,50)(0,0)
\put(-48,0){\includegraphics[scale=.3]{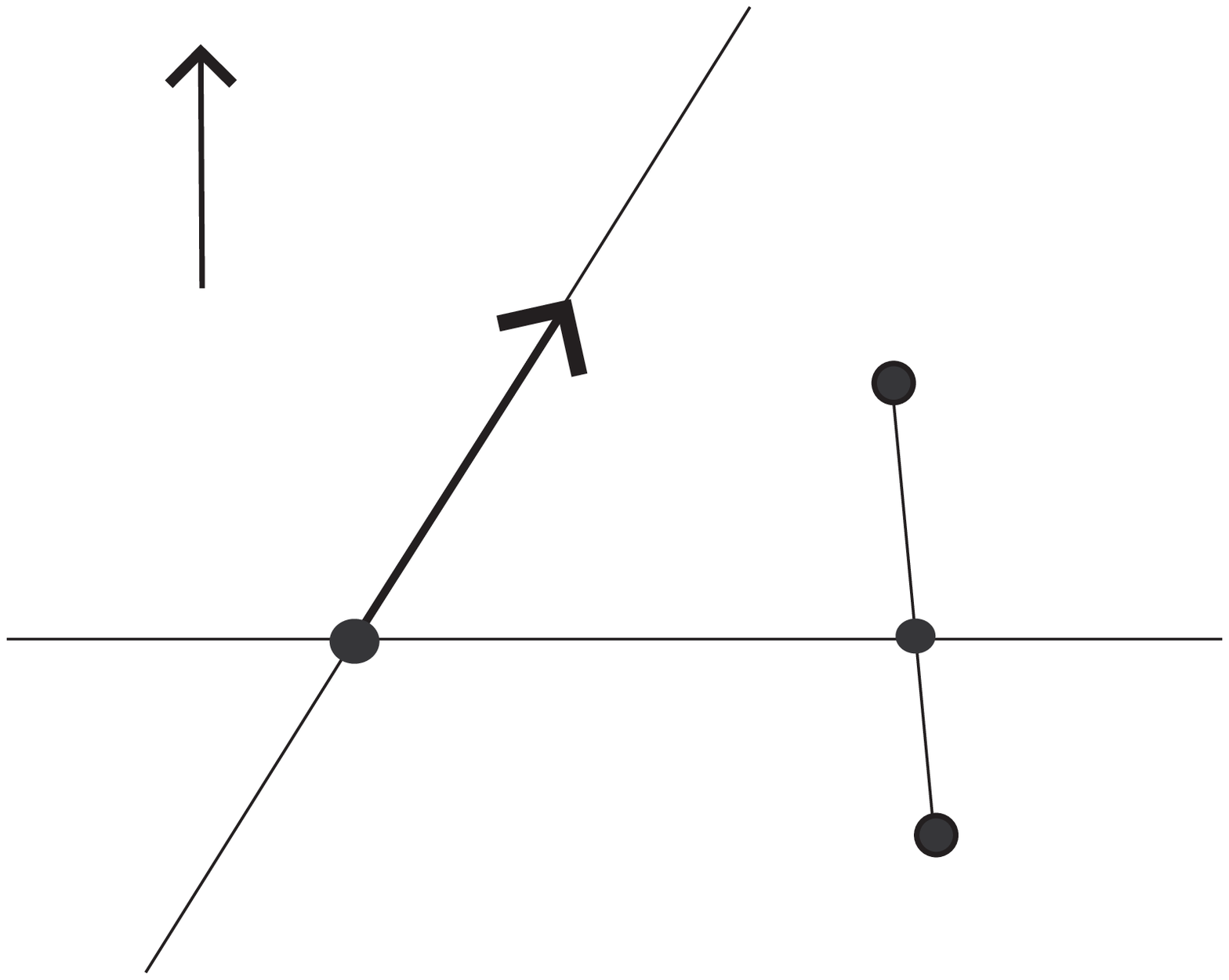}}
\put(12,35){$v_i$}
\put(15,9){$u_I$}
\put(15,23){$u_i^I$}
\put(39,16){$\mu$}
\put(41,6){$u_\alpha$}
\put(38,26){$u_\beta$}
\put(30,43){$F_j$'s}
\put(52,12){$F_i$}
\end{picture}
\caption{}\label{4.4}
\end{figure}

It follows from the definition of $\DHFP$ that the difference between 
$\DHFP(u_\alpha)$ and $\DHFP(u_\beta)$ arises from the cones 
$\Lambda^*(I)^+$'s for $I\in \Sigma^{(n)}$ 
such that $i\in I$ and $\langle u_I,v\rangle <\langle \mu,v\rangle$. 
In fact, one sees that 
\[ \DHFP(u_\alpha)+\sum_I\sign\langle \u,v\rangle (-1)^Iw(I)\phi_I(\mu)
=\DHFP(u_\beta)
\]
where $I$ runs over the elements as above.  
Since $\sign\langle \u,v\rangle(-1)^I=-(-1)^{I\backslash\{i\}}$ 
and $w(I)=w_i(I\backslash\{i\})$, the equality above turns into 
\[ \DHFP(u_\alpha)-\DHFP(u_\beta)=\sum_I(-1)^{I\backslash\{i\}}
w_i(I\backslash\{i\})\phi_I(\mu).\]
Here $\phi_I(\mu)$ may be viewed as the value at $\mu$ of the 
characteristic function of the cones in $F_i$ with apex $u_I$ spanned 
by $(u_j^I)^+$'s $(j\in I, j\not=i)$.  This shows that the right-hand 
side at the equality above agrees with $\DHF_{\PO_i}(\mu-f_i)$, 
proving the lemma. 
\end{proof}

\begin{lemm} \label{lemm:independence}
The support of the function $\DHFP$ is 
bounded, and the function 
is independent of the choice of the generic element $v\in V$. 
\end{lemm}

\begin{proof} Induction on the dimension of simple multi-polytopes $\PO$. 
We have observed the lemma in Example~\ref{exam:one_dimension} when 
$\dim \PO=1$.  Suppose $\dim \PO=n>1$ and suppose that the lemma is true 
for simple multi-polytopes of dimension $n-1$.  
Then the support of $\DHF_{\PO_i}$ is 
bounded by the induction assumption.  This together with 
Lemma~\ref{lemm:transition_formula_for_DH} implies that 
$\DHFP$ takes the same constant on unbounded regions 
in $V^*\backslash\cup F_i$.  On the other hand, it follows from the 
definition of $\DHFP$ that $\DHFP$ 
vanishes on a half space $H_r:=\{ u\in V^*\mid \langle u,v\rangle<r\}$ 
for a sufficiently small real number $r$, because for each $I\in 
\Sigma^{(n)}$ the cone $\Lambda^*(I)^+$ is contained in 
the complement of $H_r$ if $r$ is sufficiently small.  
Therefore the constant which $\DHFP$ takes on the unbounded regions in 
$V^*\backslash\cup F_i$ is zero, proving the former assertion in the 
lemma.  

As for the latter assertion in the lemma, it follows from the 
induction assumption that the right-hand side of the wall crossing formula in 
Lemma~\ref{lemm:transition_formula_for_DH} is independent of $v$, and we 
have seen above that $\DHFP$ vanishes on unbounded 
regions regardless of the choice of $v$.  
Thus, it follows from Lemma~\ref{lemm:transition_formula_for_DH} that 
$\DHFP$ is independent of $v$ on any regions of $V^*\backslash\cup F_i$. 
\end{proof}

\bigskip

\section{Winding numbers} \label{sec:WN}

We continue to assume that our multi-polytope $\PO=(\Delta,\mathcal F)$ 
is simple and that $\Sigma$ 
is an \augmented\ consisting of subsets of $\{1,\dots,d\}$.  
In this section, we associate another 
locally constant function on $V^*\backslash\cup F_i$ with 
$\PO$ from a topological viewpoint, and show that it agrees with the 
Duistermaat-Heckman function defined in Section~\ref{sec:DH}.  

Choose an orientation on $V$ and fix it.  We define an orientation on 
$I=\{i_1,\dots,i_n\}\in \Sigma^{(n)}$ as follows. 
If an ordered basis $(v_{i_1},\dots,
v_{i_n})$ gives the chosen orientation on $V$, 
then we say that the oriented simplex 
$\langle i_1,\dots,i_n\rangle$ has a positive orientation, and otherwise 
a negative orientation.  We define 
\[
\langle I\rangle :=\begin{cases}
\langle i_1,\dots,i_n\rangle \quad&\text{if $\langle i_1,\dots,i_n\rangle$ 
has a positive orientation,}\\
-\langle i_1,\dots,i_n\rangle \quad&\text{if $\langle i_1,\dots,i_n\rangle$ 
has a negative orientation.}
\end{cases}
\]
The completeness of 
$\Delta$ (equivalently, 
the pre-completeness of the projected multi-fan $\Delta_J$ 
for any $J\in \Sigma^{(n-1)}$) implies that 
\[ \sum_{I\in\Sigma^{(n)}}w(I)\langle I\rangle\]
is a cycle in the chain complex of the simplicial set $\Sigma$.  
In fact, the converse also holds, i.e., 
the completeness of $\Delta$ is equivalent to 
$\sum_{I\in\Sigma^{(n)}}w(I)\langle I\rangle$ being a cycle. 
We denote by $[\Delta]$ the homology class that the cycle 
defines in $H_{n-1}(\Sigma)$.  
Actually $[\Delta]$ lies in the reduced homology $\tilde H_{n-1}(\Sigma)$, 
see Example~\ref{exam:DH=WN} discussed later. 

Let $S$ be the realization of the first barycentric subdivision of $\Sigma$. 
For each $i\in \{1,\dots,d\}$, we denote by $S_i$ the union of simplicies 
in $S$ which contain the vertex $\{i\}$, and by $S_I$ the intersection 
$\cap_{i\in I}S_i$ for $I\in\Sigma$. Note that $\partial S_i$ can be 
identified with the realization of 
the first barycentric subdivision of $\Sigma_i$, where $\Sigma_i$ is the 
\augmented\ of the projected multi-fan $\Delta_i=
(\Sigma_i,\Lambda_i,w_i^\pm)$.  

The projected multi-fan $\Delta_i$ is defined 
on $V/V_i$ where $V_i$ is the 
one-dimensional subspace spanned by $v_i$.  We orient $V/V_i$ as follows: 
if an ordered basis $(v_i,v_{j_1},\dots,v_{j_{n-1}})$ defines 
the given orientation on $V$, then we give $V/V_i$ the orientation 
determined by $(v_{j_1},\dots,v_{j_{n-1}})$, and otherwise give the 
opposite orientaiton.  Then $[\Delta_i]$ is defined in 
$\tilde H_{n-2}(\Sigma_i)=\tilde H_{n-2}(\partial S_i)$. 

\begin{lemm} \label{lemm:Delta}
$[\Delta]$ maps to $[\Delta_i]$ through the composition of maps 
\[ \tilde H_{n-1}(\Sigma)=\tilde H_{n-1}(S) \xrightarrow{\iota_*} 
H_{n-1}(S,S\backslash \Int S_i) \xleftarrow[\cong]{\text{excision}} 
H_{n-1}(S_i,\partial S_i) \xrightarrow[\cong]{\partial}
\tilde H_{n-2}(\partial S_i),
\]
where $\iota$ is the inclusion. 
\end{lemm}

\begin{proof} Through $\iota_*$ and the inverse of the 
excision isomorphism, the cycle 
$\sum_{I\in \Sigma^{(n)}} w(I)\langle I\rangle$ maps to $\sum_{i\in I\in 
\Sigma^{(n)}}w(I)\langle I\rangle$. 
We express $\langle I\rangle$ as 
$\epsilon\langle i,j_{1},\dots,j_{n-1}\rangle$ where $\epsilon=+1$ or $-1$ 
and define an oriented $(n-2)$-simplex 
$\langle I\backslash\{i\}\rangle$ in 
$\Sigma_i^{(n-1)}$ by $\epsilon\langle j_1,\dots,j_{n-1}\rangle$.  
It follows that 
\[ \partial (\sum_{i\in I\in\Sigma^{(n)}}w(I)\langle I\rangle)=
\sum_{i\in I\in\Sigma^{(n)}}w(I)\langle I\backslash\{i\}\rangle.\]
Here $w(I)=w_i(I\backslash\{i\})$ by 
the definition of $w_i$, and $i\in I\in \Sigma^{(n)}$ 
if and only of 
$I\backslash\{i\}\in\Sigma_i^{(n-1)}$.  Therefore, the right-hand side 
above reduces to $\sum_{J\in\Sigma_i^{(n-1)}}w_i(J)\langle J\rangle$, that is 
$[\Delta_i]$ in $\tilde H_{n-2}(\partial S_i)$. 
\end{proof}

The following lemma will be used later several times. 

\begin{lemm} \label{lemm:Psi}   
Let $X$ and $Y$ be topological spaces with subspaces $X_i\subset X$ and 
$Y_i\subset Y$ for each $i\in\Sigma^{(1)}$.  For $I\in\Sigma$, we set 
$X_I:=\cap_{i\in I}X_i$ and $Y_I:=\cap_{i\in I}Y_i$.  If 
\begin{enumerate}
\item $X=\cup_{i=1}^d X_i$, 
\item $X_I$'s for $I\in \Sigma^{(n)}$ are disjoint, and
\item $Y_I$ is nonempty and contractible for any non-empty set $I\in\Sigma$,
\end{enumerate}
then there is a continuous 
map $\psi\colon X\to Y$ sending the stratum $X_I$ 
to $Y_I$ for each $I\in\Sigma$, and such a map is unique up to 
homotopy preserving the stratifications. 
\end{lemm}

\begin{proof} \emph{Existence}.  We will construct $\psi$ 
inductively using decending induction on $|I|$.  If $|I|=n$, then 
we map $X_I$ to any point in $Y_I$.  Thus $\psi$ is defined on 
$\cup_{|I|=n}X_I$ with the image in $\cup_{|I|=n}Y_I$.  Let $k$ be a 
nonnegative integer less than $n$ and $|I|=k$.  Suppose that $\psi$ is 
defined on $\cup_{|J|\ge k+1}X_J$ with the image in $\cup_{|J|\ge k+1}Y_J$.
Then 
\[ \psi\colon X_I\cap(\cup_{|J|\ge k+1}X_J)\to Y_I\cap(\cup_{|J|\ge k+1}Y_J)
\subset Y_I\]
extends to a continuous map from $X_I$ to $Y_I$ because $Y_I$ is contractible. 
Thus $\psi$ is defined on $\cup_{|I|\ge k}X_I$ with the image in 
$\cup_{|I|\ge k}Y_I$. This 
completes the induction step, so that we obtain the desired map $\psi$ 
defined on $X$.  

\emph{Uniqueness}.  We construct a homotopy $H:X\times[0,1]\to Y$ of 
given two maps $\psi_0$ and $\psi_1$ in the lemma.  The argument is almost 
same as above.  Since $Y_I$ is contractible, 
$H$ can be defined on $\cup_{|I|=n}X_I
\times [0,1]$ with $\cup_{|I|=n}Y_I$ as the image.  Let $k$ be as above and 
$|I|=k$.   Suppose that $H$ is 
defined on $(\cup_{|J|\ge k+1}X_J)\times[0,1]$ with the image in 
$\cup_{|J|\ge k+1}Y_J$ and that $H$ agrees with $\psi_t$ on 
$(\cup_{|J|\ge k+1}X_J)\times\{t\}$ for $t=0,1$. Then a map 
\[ \begin{split} 
H\cup\psi_0\cup\psi_1\colon &(X_I\cap(\cup_{|J|\ge k+1}X_J))\times[0,1]
\cup X_I\times\{0\}\cup X_I\times \{1\}  \\  
 &\to (Y_I\cap(\cup_{|J|\ge k+1}Y_J))\cup Y_I\cup Y_I=Y_I 
\end{split} \]
extends to a continuous map from $X_I\times[0,1]$ to $Y_I$ because $Y_I$ is 
contractible.  Thus $H$ is defined on 
$(\cup_{|I|\ge k}X_I)\times[0,1]$ with the image in $\cup_{|I|\ge k}Y_I$. This 
completes the induction step, so that we obtain the desired homotopy  $H$ 
defined on $X\times[0,1]$.  
\end{proof}


Lemma~\ref{lemm:Psi} can be applied with $X=S$, $X_i=S_i$, $Y=V^*$ 
and $Y_i=F_i$.  It follows that 
the multi-polytope $\PO$ associates a continuous map 
\[ \Psi\colon S\to \cup_{i=1}^d F_i\subset V^*\]
sending $S_I$ to $F_I$ for each $I\in\Sigma$ by Lemma~\ref{lemm:Psi}, 
and $\Psi$ induces a homomorphism 
\[ \Psi_*\colon \tilde H_{n-1}(S)=\tilde H_{n-1}(\Sigma)\to \tilde H_{n-1}
(V^*\backslash\{u\})\]
for each $u\in V^*\backslash\cup F_i$. 
Such $\Psi$ was first introduced in \cite{Hattori} and plays the role of 
a moment map. 
The orientation on $V$ chosen at the beginning of this section induces 
an orientation on $V^*$ in a natural way.  This determines a fundamental class 
in $H_n(V^*,V^*\backslash\{u\})$ and hence in 
$\tilde H_{n-1}(V^*\backslash\{u\})$ through $\partial \colon 
H_n(V^*,V^*\backslash\{u\})\cong \tilde H_{n-1}(V^*\backslash\{u\})$.  
We denote the fundamental class in $\tilde H_{n-1}(V^*\backslash\{u\})$ by 
$[V^*\backslash\{u\}]$.  

\begin{defi} For each $u\in V^*\backslash\cup F_i$, we define an integer 
$\WNP(u)$ by 
\[ \Psi_*([\Delta])=\WNP(u)[V^*\backslash\{u\}]\]
and call it the \emph{winding number} of the multi-polytope 
$\PO=(\Delta,\mathcal F)$ around $u$. 
\end{defi}

\begin{rema} The function $\WNP$ 
is independent of the choice of an 
orientation on $V$ because if the orientation on $V$ is reversed, then 
$[\Delta]$ and $[V^*\backslash\{u\}]$ are multiplied by $-1$ simultaneously. 
Moreover, it is locally constant and vanishes  
on unbounded regions separated by $F_i$'s, which immediately follows from 
the definition of $\WNP$. 
\end{rema}

We will see in Theorem~\ref{theo:DH=WN} below that $\WNP
=\DHFP$.  For the moment, we shall check this 
coincidence when $\dim\PO=1$. 

\begin{exam} \label{exam:DH=WN}
We use the notation in Example~\ref{exam:one_dimension}.  
We identify $V$ with $\R$, so that $V^*$ is also identified with 
$\R$.  Then $V$ and $V^*$ have standard orientations, and 
since $v_i$ gives the orientation on $V$ if and only if $i\in E$, 
the cycle which 
defines $[\Delta]$ is given by 
\[ \sum_{i\in E}w(\{i\})\langle i\rangle+\sum_{i\notin E}w(\{i\})
(-\langle i\rangle)=
-\sum_{i=1}^d(-1)^{\{i\}}w(\{i\})\langle i\rangle \]
where $(-1)^{\{i\}}$ is the same as in (5.2).  Since $\Delta$ is complete, 
$\sum_{i=1}^d(-1)^{\{i\}}w(\{i\})=0$; so 
$[\Delta]$ actually lies in $\tilde H_0(\Sigma)=\tilde H_0(S)$ and 
one can rewrite the cycle above as 
\[ \sum_{i=1}^d(-1)^{\{i\}}w(\{i\})(\langle j\rangle-\langle i\rangle)\]
for any $j\in\{1,\dots,d\}$.  Since $S_i=\{i\}$ and $\Psi(\{i\})=F_i$, 
$\WNP(u)=0$ unless $u$ is between the minimum 
value and the maximum value of $\{F_1,\dots,F_d\}$.
Suppose $u$ is between them and take $j$ such that $F_j$ is the maximum. 
Then one easily sees that 
\[ \WNP(u)=\sum_{F_i<u}(-1)^{\{i\}}w(\{i\}).\]
This together with (5.3) shows that $\WNP=\DHFP$ when $\dim \PO=1$. \qed
\end{exam}

We will show that $\WN$ satisfies the same wall crossing formula as in 
Lemma~\ref{lemm:transition_formula_for_DH}.  For that, we first state 
a lemma which expresses the winding number as a sum of local winding 
numbers so to speak.  Assume $\dim \PO >1$. 
We orient $F_i$ in such a way that the 
juxtaposition of a normal vector to $F_i$, whose evaluation on $v_i$ is 
positive, and the orientation on $F_i$ agrees with the prescribed 
orientation on $V^*$. By Lemma~\ref{lemm:Psi}, $\Psi$ maps a pair 
$(S_i,\partial S_i)$ into a pair $(F_i,F_i\backslash\{\mu\})$ 
for any $\mu\in F_i\backslash(F_i\cap(\cup_{j\in \Sigma_i^{(1)}}F_j)$. 
If we identify $F_i$ with $(V^*)_i$ through the translation by $-f_i$ as 
before, then the map $\Psi$ restricted to $\partial S_i$ agrees with the 
map (up to homotopy) constructed from the multi-polytope $\PO_i=
(\Delta_i,\mathcal F_i)$.  
It follows that 
\begin{equation} 
\Psi_*([\Delta_i])=\WN_{\PO_i}(\mu-f_i)[F_i\backslash 
\{\mu\}]. \tag{6.1}
\end{equation}

Let $u\in V^*\backslash\cup F_i$.  We choose a generic ray $R$ starting 
from $u$ with direction $\gamma\in V^*$, so that the intersection $F_i\cap R$ 
is one point for each $i$ if it is nonempty.  We denote the point 
$F_i\cap R$ by $R_i$. 

\begin{lemm} \label{lemm:ray}
$\displaystyle{\WNP(u)=
\sum_{i:F_i\cap R\not=\phi}
\sign\langle \gamma,v_i\rangle\WN_{\PO_i}(R_i-f_i)}.$
\end{lemm}

\begin{proof} Consider the following commutative diagram: 
\[
\begin{CD}
\tilde H_{n-1}(S) & \rightarrow & H_{n-1}(S,S\backslash\cup_i \Int S_i)
& \xleftarrow[\cong]{\text{excision}} & \bigoplus_i 
H_{n-1}(S_i,\partial S_i) & \xrightarrow[\cong]{\partial} & \bigoplus_i 
\tilde H_{n-2}(\partial S_i)\\
@V \Psi_*VV  @V \Psi_* VV  @V \Psi_*VV  @V \Psi_* VV \\
\tilde H_{n-1}(V^*\backslash\{u\}) & \xrightarrow[\cong]{} & 
H_{n-1}(V^*\backslash\{u\},V^*
\backslash R) & \leftarrow & \bigoplus_i H_{n-1}(F_i,F_i\backslash \{R_i\})
& \xrightarrow[\cong]{\partial} & 
\bigoplus_i \tilde H_{n-2}(F_i\backslash\{R_i\})
\end{CD}
\]
\noindent
where $i$ runs over the indices of $F_i$'s which intersect $R$. 
The element $[\Delta]\in \tilde H_{n-1}(S)$ maps to $\oplus_i [\Delta_i]\in 
\oplus_i \tilde H_{n-2}(\partial S_i)$ through the upper horizontal sequence 
by Lemma~\ref{lemm:Delta} and down to $\oplus_i 
\WN_{\PO_i}(R_i-f_i)[F_i\backslash\{ R_i\}]$ by (6.1). 

Now we trace the lower horizontal sequence from the right to the left. 
Through the inverse of $\partial$, $[F_i\backslash\{ R_i\}]$ maps to 
the fundamental class $[F_i,F_i\backslash\{R_i\}]$, and further maps to 
$\sign\langle \gamma,v_i\rangle[V^*\backslash\{u\}]\in \tilde H_{n-1}(V^*
\backslash\{u\})$, where the sign arises from the choice of the orientation 
on $F_i$.  These together with the commutativity of the diagram above 
show that 
\[ \Psi_*([\Delta])=\sum_{i:F_i\cap R\not=\phi}\sign\langle \gamma,v_i
\rangle\WN_{\PO_i}(R_i-f_i)[V^*\backslash\{u\}].\]
On the other hand, 
$\Psi_*([\Delta])=\WNP(u)[V^*\backslash\{u\}]$ 
by definition.  The lemma follows by comparing these two equalities. 
\end{proof}

\begin{lemm} \label{lemm:transition_formula_for_WN}
The wall crossing formula as in Lemma~\ref{lemm:transition_formula_for_DH} 
holds for $\WN$ instead of $\DHF$. 
\end{lemm}

\begin{proof} Subtract the equality in 
Lemma~\ref{lemm:ray} for $u=u_\beta$ from that for $u=u_\alpha$.  Since 
one can take $\gamma$ to be $u_\beta-u_\alpha$, the lemma follows. 
\end{proof}

\begin{theo} \label{theo:DH=WN}
$\DHFP=\WNP$ for any simple multi-polytope $\PO$. 
\end{theo}

\begin{proof} The equality is established in Example~\ref{exam:DH=WN} 
when $\dim\PO=1$.  Suppose $\dim\PO=n>1$ and suppose that
the equality holds for simple multi-polytopes of 
dimension $n-1$.  Both $\DHFP$ and 
$\WNP$ are locally constant, satisfy the 
same wall crossing formula (Lemma~\ref{lemm:transition_formula_for_DH}, 
Lemma~\ref{lemm:transition_formula_for_WN}) and  
$\DHF_{\PO_i}=\WN_{\PO_i}$ by 
induction assumption. Therefore, it suffices to see that 
$\DHFP$ and $\WNP$ agree on 
one region. But we know that they vanish on unbounded regions 
(Lemma~\ref{lemm:independence} and the remark after the definition of 
$\WNP$), hence 
they agree on the whole domain. This completes the induction step, proving 
the theorem. 
\end{proof}

\bigskip

\section{Ehrhart polynomials} \label{sec:Ehrhart_polynomial}

Let $P$ be a convex lattice polytope of dimension $n$ in $V^*$, where 
\lq\lq lattice polytope'' means that each vertex of $P$ lies in the lattice 
$N^*=\Hom(N,\Z)$ of $V^*=\Hom(V,\R)$.  For a positive integer $\m$, let 
$\mP:=\{ \m u\mid u\in P\}$.  It is again a convex lattice polytope in $V^*$. 
We denote by $\sharp(\mP)$ (resp. $\sharp(\mP^\circ)$) the number of lattice 
points in $\mP$ (resp. in the interior of $\mP$).  
The lattice $N^*$ determines a volume element on $V^*$ by requiring 
that the volume of the unit cube determined by a basis of $N^*$ is $1$. 
Thus the volume of $P$, denoted by $\vol(P)$, is defined. 
The following theorem is well known. 

\begin{theo}  \label{theo:Ehrhart}  
{\rm (See \cite{Fulton}, \cite{Oda} for example.)} 
Let $P$ be an $n$-dimensional convex lattice polytope. 
\begin{enumerate}
\item $\sharp(\mP)$ and $\sharp(\mP^\circ)$ are polynomials in $\m$ of degree 
$n$. 
\item $\sharp(\mP^\circ)=(-1)^n\sharp(-\mP)$, where $\sharp(-\mP)$ denotes the 
polynomial $\sharp(\mP)$ with $\m$ replaced by $-\m$.  
\item The coefficient of $\m^n$ in $\sharp(\mP)$ is $\vol(P)$ and the 
constant term in $\sharp(\mP)$ is $1$. 
\end{enumerate}
\end{theo}

The fan $\Delta$ associated with $P$ may not be simplicial, but 
if we subdivide $\Delta$, then we can always take a simplicial 
fan that is compatible with $P$.
In this section, we show that 
the theorem above holds for a {\it simple} lattice multi-polytope 
$\PO=(\Delta,\mathcal F)$. For that, we need to define $\sharp(\PO)$ 
and $\sharp(\PO^\circ)$.  This is done as follows. Let $v_i$ $(i=1,\dots,d)$ 
be a primitive integral vector in the half line $\Lambda(\{i\})$.  
In our convention, $v_i$ is chosen \lq\lq outward normal" to the face 
$\mathcal F(\{i\})$ when $\PO$ arises from a convex polytope.  We slightly 
move $\mathcal F(\{i\})$ in the direction of $v_i$ (resp. $-v_i$) for each 
$i$, so that we obtain a map $\mathcal F_+$ (resp. $\mathcal F_-$) $\colon 
\Sigma^{(1)}\to\HP(V^*)$.  We denote the multi-polytopes 
$(\Delta,\mathcal F_+)$ and $(\Delta,\mathcal F_-)$ by $\PO_+$ and $\PO_-$ 
respectively.  Since the affine hyperplanes $\mathcal F_\pm(\{i\})$'s miss 
the lattice $N^*$, the functions $\DHF_{\PO_\pm}$ and $\WN_{\PO_\pm}$ are 
defined on $N^*$.  

\begin{defi} We define 
\begin{align*}
\sharp(\PO)&:=\sum_{u\in N^*}\DHF_{\PO_+}(u)=\sum_{u\in N^*}\WN_{\PO_+}(u),\\
\sharp(\PO^\circ)&:
=\sum_{u\in N^*}\DHF_{\PO_-}(u)=\sum_{u\in N^*}\WN_{\PO_-}(u).
\end{align*}
\end{defi}

When $\PO$ arises from a convex polytope $P$, $\DHF_{\PO_+}=\WN_{\PO_+}$ 
(resp. $\DHF_{\PO_-}=\WN_{\PO_-}$) takes $1$ on $u$ in $P$ 
(resp. in the interior of $P$) and $0$ otherwise.  
Therefore, $\sharp(\PO)$ (resp. $\sharp(\PO^\circ)$) agrees with the number of 
lattice points in $P$ (resp. in the interior of $P$) in this case. 

Denote the volume element on $V^*$ by $dV^*$, and 
define the volume $\vol(\PO)$ of $\PO$ by 
\[ \vol(\PO):=\int_{V^*} \DHF_\PO\,dV^*=\int_{V^*}\WN_\PO\,dV^*.\]
When $\PO$ arises from a (convex) polytope $P$, 
$\vol(\PO)$ agrees with the actual volume of $P$, but otherwise 
it can be zero or negative. 

For a (not necessarily positive) integer $\m$, we denote 
$(\Delta,\m\mathcal F)$ by $\m\PO$, where 
\[(\m\mathcal F)(\{i\}):=\{u\in V^*\mid \langle u,v_i\rangle=\m c_i\}\]
when $\mathcal F(\{ i\})=\{ u\in V^*\mid \langle u,v_i\rangle=c_i\}$ 
for a constant $c_i$. 

\begin{theo} \label{theo:Ehrhart_polynomial_for_PO}
Let $\PO=(\Delta,\mathcal F)$ be a simple lattice 
multi-polytope of dimension $n$.  
\begin{enumerate}
\item $\sharp(\m\PO)$ and $\sharp(\m\PO^\circ)$ are polynomials in $\m$ 
of degree 
(at most) $n$. 
\item $\sharp(\m\PO^\circ)=(-1)^n\sharp(-\m\PO)$ for any integer $\m$. 
\item The coefficient of $\m^n$ in $\sharp(\m\PO)$ is $\vol(\PO)$ and the 
constant term in $\sharp(\m\PO)$ is $\deg(\Delta)$. 
$($See Section~\ref{sec:multi-fan} for $\deg(\Delta)$.$)$ 
\end{enumerate}
\end{theo}

In order to prove this theorem, we need some notations and 
a lemma.  
Basic ideas in the following arguments are in \cite{BV} and \cite{BV2}. 
Let $I\in \Sigma^{(n)}$.  Although the 
integral vectors $\{ v_i\mid i\in I\}$ are not necessarily a basis of the 
lattice $N$, they are linearly independent.  Therefore, the sublattice $N_I$ 
of $N$ generated by $v_i$'s $(i\in I)$ is of the same rank as $N$, hence 
$N/N_I$ is a finite group.  
Needless to say, $N/N_I$ is trivial for any $I\in\Sigma^{(n)}$ if $\Delta$ 
is non-singular. 
For $u\in N_I^*=\Hom(N_I,\Z)\supset N^*$ and 
$g\in N/N_I$, we define 
\begin{equation} 
\chi_I(u,g):=\exp(2\pi \sqrt{-1}\langle u,v_g\rangle)
\tag{7.1}
\end{equation}
where $v_g\in N$ is a representative of $g$.  The right-hand side does not 
depend on the choice of the representative $v_g$, and $\chi_I(u,\ )$ (resp. 
$\chi(\ ,g)$) is a homomorphism from $N/N_I$ (resp. $N_I^*$) to $\C^*$.  
Note that $\chi_I(u,\ )\colon N/N_I\to \C^*$ is trivial if and 
only if $u\in N^*$.  It follows that 
\begin{equation}
\sum_{g\in N/N_I}\chi_I(u,g)=\begin{cases} |N/N_I| \quad 
&\text{if $u\in N^*$,}\\
0 &\text{otherwise.}\end{cases}\tag{7.2}
\end{equation}

\begin{lemm} \label{lemm:key_lemma} 
For each $I\in\Sigma^{(n)}$ let $u_I$ 
be the corresponding vertex of $\PO$ and let 
$\{\u\mid i\in I\}$ be the dual basis of $\{v_i\mid i\in I\}$ 
as in Section~\ref{sec:DH}.  Then, for $v\in N$ such that 
$\langle \u,v\rangle$ is a nonzero integer 
for any $I\in \Sigma^{(n)}$ and $i\in I$, 
we have 
\[ \sum_{I\in\Sigma^{(n)}}\frac{w(I)z^{\langle u_I,v\rangle}}
{|N/N_I|}\sum_{g\in N/N_I}
\frac{1}{\prod_{i\in I}(1-\chi_I(\u,g)z^{-\langle \u,v\rangle})}
=\sum_{u\in N^*}\DHF_{\PO_+}(u)z^{\langle u,v\rangle}\]
as functions of $z\in\C$.  
\end{lemm}

\begin{proof} The Maclaurin expansion of $1/(1-az^{-m})$ 
$(a\in\C^*,\ m\in\Z)$ is given by 
\[ \begin{cases}
-a^{-1}z^m-a^{-2}z^{2m}-\dots \qquad&\text{if $m>0$}\\
1+az^{-m}+a^2z^{-2m}+\dots \qquad &\text{if $m<0$.}
\end{cases}
\]
Taking this into account, we expand the sum 
\[
S_I:=\sum_{g\in N/N_I}\frac{1}{\prod_{i\in I}(1-\chi_I(\u,g)
z^{-\langle \u,v\rangle})}
\]
into Maclaurin series and get 
\begin{align*}
S_I=& \sum_{g\in N/N_I}(-1)^I\prod_{i\in I}\sum_{\{b_i\}}
(\chi_I(\u,g)^{-b_i}z^{b_i\langle \u,v\rangle})\\
=& \sum_{g\in N/N_I}(-1)^I\sum_{\{b_i\}}\chi_I(-\sum_{i\in I}b_i\u,g)
z^{\langle \sum_{i\in I}b_i\u,v\rangle},
\end{align*}
where the summation $\displaystyle{\sum_{\{b_i\}}}$ 
runs over the collection of such 
$\{ b_i\mid i\in I,\ b_i\in \Z\}$ that 
\begin{equation}
\text{$b_i\ge 1$ for $i$ with $\langle \u,v\rangle>0$ and 
$b_i\le 0$ for $i$ with $\langle \u,v\rangle <0$},\tag{7.3}
\end{equation} 
(see Section~\ref{sec:DH} for $(-1)^I$).
Since 
\[
\sum_{g\in N/N_I}\chi_I(-\sum_{i\in I}b_i\u,g)
=\begin{cases}
|N/N_I| \quad&\text{if $\sum_{i\in I} b_i\u\in N^*$,}\\
0 &\text{otherwise,}
\end{cases}
\]
by (7.2), the Maclaurin expansion of the left-hand side of the equality 
in Lemma \ref{lemm:key_lemma} has the form 
\[ \sum_{u\in N^*}
  \left(\sum_{I\in\Sigma^{(n)}}(-1)^Iw(I)\phi_I'(u)\right)z^{\langle u,v\rangle}\]
where 
\[ \phi_I'(u)=\begin{cases} 
1\quad&\text{if $u=u_I+\sum_{i\in I}b_i\u$, $b_i$'s are 
as in (7.3) and $\sum_{i\in I} b_i\u\in N^*$,}\\
0 &\text{otherwise.}
\end{cases}
\]
One easily checks that $\sum_{I\in\Sigma^{(n)}}(-1)^Iw(I)\phi_I'(u)$ agrees 
with $\DHF_{\PO_+}(u)$, proving the lemma. 
\end{proof}

\begin{proof}[Proof of Theorem~\ref{theo:Ehrhart_polynomial_for_PO}.] 
We shall prove (2) first.  It suffices to prove $\sharp(\PO^\circ)=
(-1)^n\sharp(-\PO)$. Since $\sharp(\PO^\circ)=\sum_{u\in N^*}\WN_{\PO_-}(u)$ 
by definition, it suffices to prove that 
\begin{equation}
\WN_{\PO_-}(u)=(-1)^n\WN_{(-\PO)_+}(u) \quad\text{for any $u\in N^*$.}
\tag{7.4}
\end{equation}
Let $\Psi_{\PO_-}$ and $\Psi_{(-\PO)_+}$ be the maps introduced in 
Section~\ref{sec:WN} which are associated with 
multi-polytopes $\PO_-$ and $(-\PO)_+$ respectively.  We note that 
$\Psi_{\PO_-}$ and $-\Psi_{(-\PO)_+}$ considered as maps from $S$ to 
$V^*\backslash\{u\}$ for $u\in N^*$ are homotopic.  Since the 
multiplication by $-1$ on $V^*$ sends the fundamental class 
$[V^*\backslash\{-u\}]$ to $(-1)^n[V^*\backslash\{u\}]$, we obtain (7.4). 

We shall prove (1).  Because of (2), it suffices to prove (1) for 
$\sharp(\m\PO)$.  We apply Lemma~\ref{lemm:key_lemma} to $\m\PO$ in place 
of $\PO$ (so that $u_I$ is replaced by $\m u_I$), and 
approach $z$ to $1$ in the equality.  Since the right-hand side approaches 
$\sharp(\m\PO)$, it suffices to show that the left-hand side 
approaches a polynomial in $\m$ of degree at most $n$. 
When $g\in N/N_I$ is the identity element, $\chi_I(\u,g)=1$. 
Therefore, the term in the summand $\sum_{g\in N/N_I}$  
in the left-hand side has a pole at $z=1$ of 
degree exactly $n$ when $g$ is the identity element, and 
of degree at most $n$ otherwise. Thus the left-hand side of the equality 
in Lemma~\ref{lemm:key_lemma} applied to $\m\PO$ can be written as 
\[ \frac{\sum_{I\in\Sigma^{(n)}} 
z^{\m\langle u_I,v\rangle}h_I(z)}{(1-z)^nf(z)}
\]
where $h_I(z)$ and $f(z)$ are polynomials in $z$ and $f(1)\not=0$. 
Then the repeated use of L'Hospital's Theorem implies that when $z$ 
approaches $1$, the limit of the above rational function is a 
polynomial in $\m$ of degree at most $n$. 

Finally we prove (3). Since 
\[ \sharp(\m\PO)=\sum_{u\in H^2(BT)}\DHF_{(\m\PO)_+}
(u)=\sum_{u\in H^2(BT)/\m}\DHF_{\PO_+}(u), 
\]
it follows from the definition of definite integral that 
\[
\lim_{\m\to\infty}\frac{1}{\m^n}\sharp(\m\PO)=\lim_{n\to\infty}\frac{1}{\m^n}
\sum_{u\in H^2(BT)/\m}\DHF_{\PO_+}(u)=\int_{V^*}\DHF_\PO\,dV^*=\vol(\PO),
\]
proving that the coefficient of $\m^n$ in $\sharp(\m\PO)$ is $\vol(\PO)$. 

We apply Lemma~\ref{lemm:key_lemma} to $0\PO$, that is $\m\PO$ with $\m=0$. 
Then the $u_I$ in the lemma is zero, and $\DHF_{(0\PO)_+}(u)=\WN_{(0\PO)_+}
(u)=0$ unless $u=0$ because the origin is the only vertex of $0\PO$ so that 
the vertices of $(0\PO)_+$ are very close to the origin.  Thus the right-hand 
side of the equality in the lemma applied to $0\PO$ 
is a constant, say $c$, which is nothing 
but the constant term in $\sharp(\m\PO)$.  Now we approach $z$ to $\infty$. 
Then the equality reduces to 
\[ \sum_{v\in\Lambda(I)}w(I)=c\]
because $\langle\u,v\rangle>0$ for all $i\in I$ if and only if $v=\sum_{i\in I}
a_iv_i$ with $a_i>0$ for all $i\in I$, and the latter is equivalent to 
saying that $v$ belongs to the cone $\Lambda(I)$ spanned by $v_i$'s 
$(i\in I)$.  Since $\sum_{v\in C(I)}w(I)=\deg(\Delta)$ 
by definition, the constant term in $\sharp(\m\PO)$, that is $c$, agrees with 
$\deg(\Delta)$. 
\end{proof}

Let $N^*_\Delta$ be the lattice of $N^*_\R$ generated by all 
$\u$'s for $I\in\Sigma^{(n)}$ and $i\in I$.  
If $\Delta$ is non-singular, then $N^*_\Delta=N^*$. 
The group ring $\C[N^*_\Delta]$ is an integral domain, and it has 
a basis $t^u$ $(u\in N^*_\Delta)$ as a complex vector space with 
multiplication determined by the addition in $N^*_\Delta$:
\[ t^u\cdot t^{u'}:=t^{u+u'}.\]
The quotient field of $\C[N^*_\Delta]$ will be denoted by
$\C(N^*_\Delta)$. It contains $\C[N^*_\Delta]$. 
Each $v\in N$ such that $\langle\u,v\rangle$ is an integer for any 
$I\in\Sigma^{(n)}$ and $i\in I$ 
determines a map from $\C[N^*_\Delta]$ to a Laurent polynomial 
ring $\C[z,z^{-1}]$ sending $t^u$ to $z^{\langle u,v\rangle}$.
This map extends to a map from $\C(N^*_\Delta)$ to $\C(z)$, the field of 
rational functions in $z$. 
Since Lemma~\ref{lemm:key_lemma} holds for any such $v$ that 
$\langle\u,v\rangle\not=0$, we obtain 

\begin{coro}\label{coro:key_corollary} Let the notation be the same 
as in Lemma~\ref{lemm:key_lemma}.  Then 
\[ \sum_{I\in\Sigma^{(n)}}\frac{w(I)t^{u_I}}{|N/N_I|}\sum_{g\in N/N_I}
\frac{1}{\prod_{i\in I}(1-\chi_I(\u,g)t^{-\u})}
=\sum_{u\in N^*}\DHF_{\PO_+}(u)t^{u}\in \C[N^*]\]
as elements in $\C(N^*_\Delta)$. 
In particular, if the multi-fan $\Delta$ is non-singular, then 
$N^*_\Delta=N^*$ and 
\[ \sum_{I\in\Sigma^{(n)}}
\frac{w(I)t^{u_I}}{\prod_{i\in I}(1-t^{-\u})}
=\sum_{u\in N^*}\DHF_{\PO_+}(u)t^{u}.\]
\end{coro}

For a later use, we shall rewrite $\chi_I(\u,g)$.  Consider a 
homomorphism $\eta\colon \R^d\to N_\R$ 
mapping $\a=(a_1,\dots,a_d)\in\R^d$ to $\sum_{i=1}^d a_iv_i\in N_\R$. 
For $I\in\Sigma^{(n)}$, we define 
\[
G_I':=\{\a\in\R^d\mid \eta(\a)\in N 
\text{ and }a_j=0 \text{ for } j\notin I\}
\]
and define $G_I$ to be the projection image of $G_I'$ on $\R^d/\Z^d$. 
Since $v_i$'s $(i\in I)$ are linearly independent and belong to $N$, 
$G_I$ is a finite subgroup of $\R^d/\Z^d$ and 
$\eta$ restricted to $G_I'$ induces an isomorphism 
\[
\eta_I\colon G_I\cong N/N_I.
\]
Note that $\eta_I([\a])=[\sum_{i\in I}a_iv_i]$ where $[\ \ ]$ denotes the 
equivalence class.  

On the other hand, for $i=1,\dots,d$, let 
\[
\rho_i\colon \R^d/\Z^d\to \C^*
\]
be a homomorphism defined by $\rho_i([\a])=\exp(2\pi\sqrt{-1}a_i)$.  

\begin{lemm} 
For $[\a]\in G_I\subset \R^d/\Z^d$ and $i\in I$, we have 
$\rho_i([\a])=\chi_I(\u,\eta_I([\a]))$. 
\end{lemm}

\begin{proof} 
Since $\eta_I([\a])=[\sum_{i\in I}a_iv_i]$ and 
$\langle\u,\sum_{i\in I} a_iv_i\rangle=a_i$, it follows from 
the definition (7.1) of $\chi_I$ that 
$\chi_I(\u,\eta_I([\a]))=\exp(2\pi\sqrt{-1}a_i)$, 
which is equal to $\rho_i([\a])$ by definition.  
\end{proof}

Since $G_I$ is isomorphic to $N/N_I$, 
Corollary~\ref{coro:key_corollary} can be restated as follows. 

\begin{coro}\label{coro:rho} 
Let the notation be as above.  Then 
\[ \sum_{I\in\Sigma^{(n)}}\frac{w(I)t^{u_I}}{|G_I|}\sum_{g\in G_I}
\frac{1}{\prod_{i\in I}(1-\rho_i(g)t^{-\u})}
=\sum_{u\in N^*}\DHF_{\PO_+}(u)t^{u}\in \C[N^*]\]
as elements in $\C(N^*_\Delta)$. 
\end{coro}

\bigskip

\section{Cohomological formula for $\sharp(\PO)$} 
\label{sec:cohomological_formula}

Motivated by the geometrical observation which will be explained in subsequent 
sections~\ref{sec:torus_manifold} and ~\ref{sec:equivariant_index}, 
we define the \lq\lq (equivariant) cohomology" of a complete 
simplicial multi-fan and the \lq\lq (equivariant) first Chern class" of a 
multi-polytope.  
We then define an index map \lq\lq in cohomology" and establish a 
\lq\lq cohomological formula" describing $\sharp(\PO)$ for a lattice 
multi-polytope.  
This cohomological formula is a counterpart in combinatorics to the 
Hirzebruch-Riemann-Roch formula applied to a complex $T$-line bundle 
over a torus manifold. As an application of the cohomological formula, 
we show that the Khovanskii-Pukhlikov formula for 
a simple lattice convex polytope (\cite{KK}, \cite{KP}, \cite{CS1}, \cite{CS2},
\cite{Guill}, \cite{BV}, \cite{BV2}) can be generalized 
to a simple lattice multi-polytope. 

Let $T$ be a compact torus of dimension $n=\rank_{\Z}N$ and let 
$BT$ be the classifying space of $T$.  Then $H_2(BT)$ is canonically 
isomorphic to $\Hom(S^1,T)$, the group consisting of homomorphisms 
from $S^1$ to $T$.  In fact, a homomorphism $f\colon S^1\to T$ induces 
a continuous map $Bf\colon BS^1\to BT$ and once we fix a generator 
$\alpha$ of $H_2(BS^1)\cong\Z$, $(Bf)_*\alpha$ defines an element 
of $H_2(BT)$.  The correspondence $: f\to (Bf)_*\alpha$ 
is known to be an isomorphism 
from $\Hom(S^1,T)$ to $H_2(BT)$.  In the following we assume 
$N=H_2(BT)$ and identify it with $\Hom(S^1,T)$.  Then $N^*=H^2(BT)$ 
is identified with $\Hom(T,S^1)$ and 
the group ring $\C[N^*]$ can be identified with the representation 
ring $R(T)$ of $T$. 

Let $\Delta=(\Sigma,\Lambda,w^\pm)$ be a complete simplicial 
multi-fan in $N$.  
Let $v_i\in H_2(BT)$ be a unique primitive vector in $\Lambda(\{i\})$ 
for each $i=1,\dots,d$ as before. 
Motivated by the description of the equivariant cohomology of a 
compact non-singular toric variety (see 
Proposition~\ref{prop:ring_structure} in the next section), 
we define $H^*_T(\Delta)$ to be 
the face ring of the \augmented\ $\Sigma$, i.e., 
\[ H^*_T(\Delta):=\Z[x_1,\dots,x_d]/(x_I\mid I\notin\Sigma),
\] 
where $x_I=\prod_{i\in I}x_i$ and the degree of $x_i$ is two, and 
call $H^*_T(\Delta)$ the \emph{equivariant cohomology} of $\Delta$. 
We also define 
a homomorphism $\pi^*\colon H^2(BT)\to H^2_T(\Delta)$ by 
\begin{equation}
\pi^*(u)=\sum_{i=1}^d\langle u,v_i\rangle x_i,\tag{8.1}
\end{equation}
where $\langle\ ,\ \rangle$ denotes the natural pairing between 
cohomology and homology.  
It extends to an algebra homomorphism $H^*(BT)\to H^*_T(\Delta)$, which 
we also denote by $\pi^*$.  One can think of $H^*_T(\Delta)$ as 
a module (or more generally an algebra) over $H^*(BT)$ through $\pi^*$.  

In the following we will mainly work with $\Q$ coefficients but the 
argument will work with $\Z$ coefficients when the multi-fan $\Delta$ 
is non-singular.  Any homomorphism $f\colon A\to B$ between additive 
groups induces a homomorphism $: A\otimes\Q\to B\otimes\Q$ 
(or $A\otimes\R\to B\otimes \R$), which we also 
denote by $f$.  

\begin{lemm} \label{lemm:generator} Any element in 
$H^*_T(\Delta)\otimes\Q$ can be 
written in the form $\sum_{J\in\Sigma}\pi^*(a_J)x_J$ with 
$a_J\in H^*(BT;\Q)$ (not necessarily uniquely), 
in other words, $H_T^*(\Delta)\otimes\Q$ is generated by 
$x_J$'s $(J\in\Sigma)$ as an $H^*(BT;\Q)$-module.  
\end{lemm}

\begin{proof} Let $\I$ denote a finite set which consists of elements 
in $\{1,\dots,d\}$ taken with multiplicity, i.e., elements in $\{1,\dots,
d\}$ may appear in $\I$ repeatedly.  Set $x_\I:=\prod_{i\in\I}x_i$ and 
denote by $\bI$ the subset of $\{1,\dots,d\}$ consisting of elements 
appearing in $\I$.  It follows from the definition that $H^*_T(\Delta)$ 
is additively generated by $x_\I$'s such that 
$\bI\in\Sigma$, so it suffices to prove the lemma for such 
$x_\I$.  We shall prove it by induction on $[\I]:=|\I|-|\bI|$.  

If $[\I]=0$, then $\I=\bI\in\Sigma$; so $x_\I$ is obviously of the form 
in the lemma in this case.  Suppose $[\I]\ge 1$.  Then there is an 
$i\in\I$ which appears in $\I$ at least twice.  Set $\J:=\I\backslash
\{i\}$. Then $\bJ=\bI\in\Sigma$ and $[\J]=[\I]-1$.  Multiplying the 
both sides at (8.1) by $x_\J$, we obtain 
\[ \pi^*(u)x_\J=\langle u,v_i\rangle x_\I+\sum_{k\not=i}\langle u,
v_k\rangle x_{\J\cup\{k\}}
\]
for any $u\in H^2(BT;\Q)$. We choose $u$ such that $\langle u,v_i
\rangle=1$ and $\langle u,v_j\rangle=0$ for all $j\in \J$ different 
from $i$.  (Such $u$ exists because $\{ v_j\mid j\in \bJ\}$ is a subset 
of a basis of $N_\Q$.)  Then the equality above reduces to 
\[ x_\I=\pi^*(u)x_\J-\sum_{k\not=i,k\notin\J}\langle u,v_k\rangle 
x_{\J\cup\{k\}}.\]
Here $[\J\cup\{k\}]=[\J](=[\I]-1)$ for $k\notin\J$, so 
the right-hand side above are of the form in the lemma by the induction 
assumption, showing that so is $x_\I$.  This completes the induction 
step and proves the lemma. 
\end{proof}

For $I\in\Sigma^{(n)}$, let $\{\u\mid i\in I\}$ be the dual basis of 
$\{v_i\mid i\in I\}$ as before.  We define a ring homomorphism 
$\iota_I^*\colon H^*_T(\Delta)\otimes\Q\to H^*(BT;\Q)$ by 
\[ \iota_I^*(x_i)=\begin{cases} \u \quad&\text{if $i\in I$,}\\
0 \quad&\text{otherwise.}
\end{cases}
\]
This map is well-defined because $x_J$ for $J\notin\Sigma$, which 
is zero in $H^*_T(\Delta)\otimes\Q$, maps to zero through $\iota_I^*$. 

\begin{lemm} \label{lemm:injectivity} 
The composition $\iota_I^*\circ\pi^*$ is the identity map. In particular 
$\iota_I^*$ is an $H^*(BT;\Q)$-module map.   
\end{lemm}

\begin{proof} Both $\pi^*$ and $\iota_I^*$ are ring homomorphisms 
and $H^*(BT)$ is a polynomial ring generated by elements in $H^2(BT)$, 
so it suffices to check the lemma on 
$H^2(BT)$.  Let $u\in H^2(BT)$.  It follows from the definitions 
of $\pi^*$ and $\iota_I^*$ that 
\[(\iota_I^*\circ\pi^*)(u)=\iota_I^*(\sum_{i=1}^d\langle u,v_i\rangle 
x_i)=\sum_{i=1}^d\langle u,v_i\rangle \u,
\]
which agrees with $u$ because $\{\u\mid i\in I\}$ is the dual basis of 
$\{v_i\mid i\in I\}$.  Since $u$ is arbitrary, this proves that 
$\iota_I^*\circ\pi^*$ is the identity on $H^2(BT)$. 
\end{proof}

A multi-polytope $\PO=(\Delta,\mathcal F)$ is associated with 
real numbers $c_i$'s by 
\[ \mathcal F(\{i\})=\{ u\in H^2(BT;\R)\mid \langle u,v_i\rangle=c_i\},
\]
and these numbers determine an element $c_1^T(\PO):=
\sum_{i=1}^d c_ix_i$ of $H^2_T(\Delta)\otimes\R$, 
which we call the \emph{equivariant first Chern class} of $\PO$. 
This gives a bijective correspondence between the set of 
multi-polytopes defined on $\Delta$ and $H^2_T(\Delta)\otimes\R$. 
Note that $\iota_I^*(c_1^T(\PO))$ agrees with the vertex 
$u_I=\cap_{i\in I}\mathcal F(\{i\})$, see Section 5. 
When $\Delta$ is non-singular, $\PO$ is a lattice multi-polytope 
if and only if the $c_i$'s are all integers, but otherwise the 
\lq\lq if'' part does not hold, in other words, an element of 
$H^2_T(\Delta)$ is not necessarily realized as the equivariant 
first Chern class of a lattice multi-polytope.  However, there is 
a nonzero integer $m$ such that $mx$ for any $x\in H_T^2(\Delta)$ 
is realized as the equivariant first Chern class 
of a lattice multi-polytope 
because $|N/N_I|\iota_I^*(x)$'s lie in $H^2(BT)$.  

We set $H^{**}(BT;\Q)=\prod_{q=0}^\infty H^q(BT;\Q)$. It is a formal power
series ring.
\begin{lemm} \label{lemm:chern_character}
For any $J\in\Sigma$, the element 
\[ 
\sum_{I\in \Sigma^{(n)}}\frac{w(I)\iota_I^*(\prod_{j\in J}(e^{mx_j}-1))}
{|G_I|}\sum_{g\in G_I}\frac{1}{\prod_{i\in I}(1-\rho_i(g)e^{-\u})}
\]
in the quotient field of $H^{**}(BT;\Q)$ actually belongs to $H^{**}(BT;\Q)$.
\end{lemm}

\begin{proof} Since $\prod_{j\in J}(e^{mx_j}-1)$ is a linear combination 
of $\prod_{k\in K}e^{mx_k}=e^{m\sum_{k\in K}x_k}$ for $K\in \Sigma$, 
it suffices to show that 
\begin{equation}
\sum_{I\in\Sigma^{(n)}}\frac{w(I)\iota_I^*(e^{m\sum_{k\in K}x_k})}
{|G_I|}\sum_{g\in G_I}\frac{1}{\prod_{i\in I}(1-\rho_i(g)e^{-\u})}
\in H^{**}(BT;\Q).\tag{8.2}
\end{equation}
As remarked above, 
$m\sum_{k\in K}x_k$ is realized as the equivariant first Chern 
class of a lattice multi-polytope, so 
it follows from Corollary~\ref{coro:rho} that 
\[ 
\sum_{I\in \Sigma^{(n)}}\frac{w(I)t^{\iota_I^*(m\sum x_k)}}{|G_I|}
\sum_{g\in G_I}\frac{1}{\prod_{i\in I}(1-\rho_i(g)t^{-\u})}\in 
\C[N^*]=R(T).
\]
The Chern character $ \C[N^*]=R(T)\to H^{**}(BT;\Q)$ mapping $t^u$ 
to $e^u$ extends to a map from $\C[N^*_\Delta]$ and it further extends 
to a map between the quotient fields.  Sending the element above 
by this extended Chern character, we obtain (8.2).  
\end{proof}

Let $S$ be the multiplicative set consisting of nonzero homogeneous 
elements of positive degree in $H^*(BT;\Q)$.  
Since $H^*(BT;\Q)$ is a polynomial ring, it can be thought of as a subring of 
the localized ring $S^{-1}H^*(BT;\Q)$.  We define the index map 
\[ \pi_!\colon H^*_T(\Delta)\otimes\Q\to S^{-1}H^*(BT;\Q)\]
\lq\lq in cohomology" by
\[ \pi_!(x):=\sum_{I\in\Sigma^{(n)}}\frac{w(I)\iota_I^*(x)}
{|G_I|\prod_{i\in I}\u}
\]
(cf. \cite[(3.8)]{AB}). 
This map decreases degrees by $2n$, and is an $H^*(BT;\Q)$-module map 
by Lemma~\ref{lemm:injectivity}.

\begin{lemm} \label{lemm:integrality}
The image of $\pi_!$ lies in $H^*(BT;\Q)$.
\end{lemm}

\begin{proof}
Since $\pi_!$ is an $H^*(BT;\Q)$-module map, it suffices to check the lemma 
for elements $x_J$'s $(J\in\Sigma)$ by Lemma~\ref{lemm:generator}. 
We distinguish two cases. 

\emph{Case 1.} The case where $|J|=n$, i.e., $J\in \Sigma^{(n)}$. 
In this case 
\[ \iota_I^*(x_J)=\begin{cases} \prod_{i\in I}\u \quad&\text{if 
$I=J$,}\\
0 &\text{otherwise.}
\end{cases}\]
Therefore 
\[
\pi_!(x_J)=\sum_{I\in \Sigma^{(n)}}\frac{w(I)\iota_I^*(x_J)}
{|G_I|\prod_{i\in I}\u}=\frac{w(J)}{|G_J|}\in H^0(BT;\Q).
\]

\emph{Case 2.} The case where $|J|<n$.  In this case we will show that 
$\pi_!(x_J)=0$.  
Since $\rho_i(g)=1$ for any $i\in I$ if and only if $g$ is the 
identity, and 
\begin{align*}
\prod_{i\in I}(1-e^{-\u})&=(\prod_{i\in I}\u)(1+\text{higher degree term})\\
\prod_{j\in J}(e^{mx_j}-1)&=m^{|J|}x_J(1+\text{higher degree term}),
\end{align*}
the term of lowest degree in Lemma~\ref{lemm:chern_character} 
(up to a nonzero constant multiple) is 
\[ \sum_{I\in\Sigma^{(n)}}\frac{w(I)\iota^*_I(x_J)}{|G_I|\prod_{i\in I}\u},
\]
that is, $\pi_!(x_J)$, and Lemma~\ref{lemm:chern_character} 
tells us that it is an element of 
$H^*(BT;\Q)$. This means that $\pi_!(x_J)=0$ because the degree of
$\pi_!(x_J)$ is equal to $2|J|-2n<0$.  
\end{proof}

Now, motivated by the description of the cohomology ring 
of a compact non-singular toric variety (see p.106 in \cite{Fulton}), 
we define $H^*(\Delta)$ to be the quotient ring of $H_T^*(\Delta)$ 
by the ideal generated by $\pi^*(H^2(BT))$, in other words, 
\[ H^*(\Delta):=\Z[x_1,\dots,x_d]/\mathfrak A,\]
where $\mathfrak A$ is the ideal generated by all 
\begin{enumerate}
\item $x_I$ for $I\notin\Sigma$,
\item $\sum_{i=1}^d\langle u,v_i\rangle x_i$ for $u\in N$. 
\end{enumerate}
Since $\pi_!$ is an $H^*(BT;\Q)$-module map and $H^*(BT;\Q)/(H^2(BT;\Q))$ 
is isomorphic to $H^0(BT;\Q)=\Q$, $\pi_!$ induces a homomorphism 
\[ \int_\Delta\colon H^*(\Delta)\otimes\Q\to \Q,\]
where only elements of degree $2n$ in $H^*(\Delta)\otimes\Q$ survive through 
the map $\displaystyle{\int_\Delta}$.  

Remember that $G_I$ is a finite subgroup of $\R^d/\Z^d$.  We 
denote by $G_\Delta$ the union of $G_I$ over all $I\in\Sigma^{(n)}$.
Since $\rho_i$ is defined on $\R^d/\Z^d$, $\rho_i(g)$ makes sense 
for $g\in G_\Delta$.  It follows from the definition of $G_I$ and 
$\rho_i$ that if $g\in G_I$, then $\rho_i(g)=1$ for $i\notin I$. 

We define the \emph{equivariant Todd class} $\mathcal T^T(\Delta)$ 
of the complete simplicial multi-fan $\Delta$ by 
\[ \mathcal T^T(\Delta):=\sum_{g\in G_\Delta}
\prod_{i=1}^d\frac{x_i}{1-\rho_i(g)e^{-x_i}}
\in H^{**}_T(\Delta)\otimes\Q,\]
and the \emph{Todd class} 
$\mathcal T(\Delta)$ of $\Delta$ by 
\[ \mathcal T(\Delta):=\sum_{g\in G_\Delta}
\prod_{i=1}^d\frac{\bar x_i}{1-\rho_i(g)e^{-\bar x_i}}
\in H^{**}(\Delta)\otimes\Q,\]
where $\bar x_i$ denotes the image of $x_i\in H^*_T(\Delta)$ in 
$H^*(\Delta)$ (cf. \cite{BV2}). 
We also define the \emph{first Chern class} $c_1(\PO)$ of a 
multi-polytope $\PO$ 
defined on $\Delta$ to be the image of 
$c_1^T(\PO)\in H^2_T(\Delta)\otimes\R$ 
in $H^2(\Delta)\otimes\R$. 

\begin{theo} \label{theo:integral_formula}
If $\PO$ is a simple lattice multi-polytope, then 
$\displaystyle{\int_\Delta e^{c_1(\PO)}\mathcal T(\Delta)=\sharp(\PO).}$
\end{theo}

\begin{proof} We shall compute 
$\pi_!(e^{c_1^T(\PO)}\mathcal T^T(\Delta))$.  
For that, we need to see $\iota^*_I(\mathcal T^T(\Delta))$.  
Let $g\in G_\Delta$.  If $g\notin G_I$, then there is an $i\notin I$ 
such that $\rho_i(g)\not=1$; so 
\[
\iota_I^*\Big(\frac{x_i}{1-\rho_i(g)e^{-x_i}}\Big)=0
\]
for such $i$ because the Maclaurin expansion of 
$x_i/(1-\rho_i(g)e^{-x_i})$ has no constant term and 
$\iota_I^*(x_i)=0$.  Therefore, only elements $g$ in $G_I$ contribute 
to $\iota_I^*(\mathcal T^T(\Delta))$.  Now suppose $g\in G_I$.  
Then $\rho_i(g)=1$ for $i\notin I$, so 
\[
\iota_I^*\Big(\frac{x_i}{1-\rho_i(g)e^{-x_i}}\Big)=1
\]
for such $i$ because the Maclaurin expansion of 
$x_i/(1-\rho_i(g)e^{-x_i})$ has the constant term $1$ and 
$\iota_I^*(x_i)=0$.  Finally, since $\iota_I^*(x_i)=\u$ for $i\in I$, 
we thus have  
\[
\iota_I^*(\mathcal T^T(\Delta))=\sum_{g\in G_I}\prod_{i\in I}
\frac{\u}{1-\rho_i(g)e^{-\u}}.
\]
This together with the definition of $\pi_!$ and 
Corollary~\ref{coro:rho} shows that 
\begin{align*}
\pi_!(e^{c_1^T(\PO)}\mathcal T^T(\Delta))=&
\pi_!\Big(e^{c_1^T(\PO)}\sum_{g\in G_\Delta}\prod_{i=1}^d\frac{x_i}
{1-\rho_i(g)e^{-x_i}}\Big)\\
=&\sum_{I\in\Sigma^{(n)}}\frac{w(I)e^{\iota_I^*(c_1^T(\PO))}}{|G_I|}
\sum_{g\in G_I}\frac{1}{\prod_{i\in I}(1-\rho_i(g)e^{-\u})}\\
=&\sum_{u\in H^2(BT)}\DHF_{\PO_+}(u)e^u.
\end{align*}
This implies that 
\[ \int_\Delta e^{c_1(\PO)}\mathcal T(\Delta)=\sum_{u\in H^2(BT)}
\DHF_{\PO_+}(u)=\sharp(\PO).
\]
\end{proof}

\begin{rema} The argument developed above in this section is purely 
combinatorial, but it is possible to take a topological approach. 
Namely, associated with a complete simplicial multi-fan $\Delta$, 
one can construct a torus space $M_\Delta$ with $H^*_T(M_\Delta;\Q)=
H^*_T(\Delta)\otimes\Q$ (see \cite{DJ}). It is not 
necessarily a manifold but has a fundamental class so that the 
equivariant Gysin homomorphism 
$\pi_!\colon H^*_T(M_\Delta;\Q)=H^*_T(\Delta)\otimes\Q\to 
H^{*-2n}_T(pt;\Q)=
H^{*-2n}(BT;\Q)$, that is, the index map, can be defined. 
\end{rema}

As an application of the theorem above, we shall show that 
Khovanskii-Pukhlikov formula, which relates a certain variation 
of the volume of a simple convex lattice 
polytope to the number of lattice 
points in it, can be generalized to simple multi-polytopes. 
We begin with 

\begin{lemm} \label{lemm:volume}
$\displaystyle{\vol(\PO)=\frac{1}{n!}\int_\Delta c_1(\PO)^n
=\int_\Delta e^{c_1(\PO)}}$ for 
a simple multi-polytope $\PO$. 
\end{lemm}

\begin{proof} The latter equality is obvious because only elements of 
degree $2n$ in $H^*(\Delta)\otimes\R$ survive through the map 
$\int_\Delta$.  We shall prove the former equality. 

\emph{Step 1.} If $\PO$ is a lattice multi-polytope, then 
Theorem~\ref{theo:integral_formula} applied to $\m\PO$ for any integer 
$\m$ implies 
\[ \int_\Delta e^{c_1(\m\PO)}\mathcal T(\Delta)=\sharp(\m\PO).\]
We compare the coefficients of $\m^n$ at the both sides above.  Since 
$c_1(\m\PO)=\m c_1(\PO)$, the coefficient of $\m^n$ at the left-hand side 
is $\frac{1}{n!}\int_\Delta c_1(\PO)^n$, 
while the one at the right-hand side is $\vol(\PO)$ by 
Theorem~\ref{theo:Ehrhart_polynomial_for_PO} (3).  Therefore the lemma 
is proven for a lattice multi-polytope $\PO$. 

\emph{Step 2.} If $\PO$ is \emph{rational}, by which 
we mean that there is a nonzero integer $m$ such that $m\PO$ is a lattice 
multi-polytope, then $\vol(m\PO)=\frac{1}{n!}\int_\Delta c_1(m\PO)^n$ by 
Step 1. 
Since $\vol(m\PO)=m^n\vol(\PO)$ and $c_1(m\PO)=mc_1(\PO)$, the lemma 
is proven for a rational multi-polytope $\PO$.  

\emph{Step 3.} The functions $\vol(\cdot)$ and $\int_\Delta c_1(\cdot)^n$ 
are defined on the vector space $H^2_T(\Delta)\otimes\R$ through the 
equivariant first Chern class, and they are obviously continuous.  
By Step 2 they agree on all rational multi-polytopes which form a 
dense subset of the vector space, so they must agree on the entire 
vector space by continuity.  This completes the proof of the lemma. 
\end{proof}
Multi-polytopes defined on $\Delta$ 
form a vector space isomorphic to $H^2_T(\Delta)\otimes\R$ 
through the equivariant first Chern class, 
and Lemma~\ref{lemm:volume} implies that the 
volume function is a homogeneous polynomial function of degree $n$.  
In fact, if one writes 
$c_1^T(\PO)=\sum_{i=1}^d c_ix_i$, then 
$\vol(\PO)$ is a homogeneous polynomial in $c_1,\dots,c_d$ of degree $n$. 

For $h=(h_1,\dots,h_d)\in\R^d$, we denote by $\PO_h$ a multi-polytope 
with $c_1^T(\PO_h)=\sum_{i=1}^d (c_i+h_i)x_i$.  
Since $c_1(\PO_h)=\sum_{i=1}^d(c_i+h_i)\bar x_i$, 
Lemma~\ref{lemm:volume} 
applied to $\PO_h$ implies that $\vol(\PO_h)$ is a polynomial 
in $h_1,\dots,h_d$ (of total degree $n$).  We define 
the \emph{Todd operator} as follows: 
\[ \mathcal T(\partial/\partial h):=\sum_{g\in G_\Delta}
\prod_{i=1}^d \frac{\partial/
\partial h_i}{1-\rho_i(g)e^{-\partial/\partial h_i}}.\]
Although the Todd operator is of infinite order, 
its operation on $\vol(\PO_h)$ converges because $\vol(\PO_h)$ is a 
polynomial in $h_1,\dots,h_d$.  The following theorem extends 
the Khovanskii-Pukhlikov formula 
to simple lattice multi-polytopes. 

\begin{theo} \label{theo:KP_formula} If $\PO$ is a simple lattice 
multi-polytope, then 
\[
\mathcal T(\partial/\partial h)\vol(\PO_h)|_{h=0}=
\sharp(\PO).
\]
\end{theo}

\begin{proof} An elementary computation shows that 
\[ \frac{\partial/\partial h_i}{1-\rho_i(g)e^{-\partial/\partial h_i}}
e^{(c_i+h_i)\bar x_i}|_{h_i=0}=e^{c_i\bar x_i}\frac{\bar x_i}
{1-\rho_i(g)e^{-\bar x_i}}.\]
Therefore, it follows from Lemma~\ref{lemm:volume} and 
Theorem~\ref{theo:integral_formula} that 
\begin{align*}
\mathcal T(\partial/\partial h)\vol(\PO_h)|_{h=0}&=\mathcal T(\partial/
\partial h)\int_\Delta e^{c_1(\PO_h)}|_{h=0}
=\mathcal \int_\Delta \mathcal T(\partial/
\partial h) e^{c_1(\PO_h)}|_{h=0}\\
&=\int_\Delta\sum_{g\in G_\Delta}\prod_{i=1}^d \frac{\partial/\partial h_i}
{1-\rho_i(g)e^{-\partial/\partial h_i}} e^{(c_i+h_i)\bar x_i}|_{h_i=0}\\
&=\int_\Delta\sum_{g\in G_\Delta}\prod_{i=1}^d e^{c_i\bar x_i}
\frac{\bar x_i}{1-\rho_i(g)e^{-\bar x_i}}\\
&=\int_\Delta e^{c_1(\PO)}\mathcal T(\Delta)=\sharp(\PO), 
\end{align*}
proving the theorem. 
\end{proof}

\begin{rema}
One can reformulate the Khovanskii-Pukhlikov formula as follows. 
As remarked above, the volume function $\vol$ is a polynomial in 
$c_1,\dots,c_d$, so one can apply the Todd operator 
$\mathcal T(\partial/\partial c)$ (with the variables $c=(c_1,\dots,c_d)$ 
instead of $h=(h_1,\dots,h_d)$) 
to the volume function $\vol$ and evaluate at a simple lattice multi-polytope 
$\PO$.  The same argument as in the proof of 
Theorem~\ref{theo:KP_formula} shows that 
the evaluated value agrees with $\sharp(\PO)$.  
\end{rema}

\bigskip

\section{Multi-fan of a torus manifold} \label{sec:torus_manifold}

In this section we introduce the notion of a torus manifold and 
associate a complete non-singular multi-fan with it. A compact 
non-singular toric variety provides an example of a torus manifold, 
but the class of torus manifolds is much wider than that of compact 
non-singular toric varieties, (apparently, even wider than that of 
unitary toric manifolds introduced in \cite{Masuda}).
The basic theory of toric varieties says that there is a one-to-one 
correspondence between compact non-singular toric varieties and 
complete non-singular fans.   This correspondence is extended 
in one direction, namely from torus 
manifolds to complete non-singular multi-fans.
But the usual way to associate a fan with a toric variety (see 
\cite[Section 2.3]{Fulton}) does not work in our extended category.  
However, when a toric variety is compact and non-singular,  
the corresponding (complete and non-singular) fan can be reproduced using 
equivariant cohomology and this argument works even for torus manifolds. 
The idea is essentially same as in \cite{Masuda}.

We begin with the definition of a torus manifold.  
An elementary representation theory of a torus group tells us that 
if an $m$-dimensional torus $(S^1)^m$ acts effectively and smoothly on 
a connected smooth manifold of dimension $2n$ with non-empty fixed point set, 
then $m\le n$ and the dimension of the fixed point set is at most 
$2(n-m)$.  We are interested in an extreme case $m=n$.  
Let $M$ be a closed, connected, smooth manifold of dimension 
$2n$ with an effective smooth action of an $n$-dimensional torus group 
$T=(S^1)^n$ such that the fixed point set $M^T$ is non-empty.  
Then $M^T$ is necessarily isolated. 
A closed, connected, codimension two submanifold of $M$ is called 
\emph{characteristic} if it is a connected component of the set 
fixed pointwise by a certain circle subgroup of $T$ and contains 
at least one $T$-fixed point.  Since $M$ is compact, there are only 
finitely many characteristic submanifolds.  We denote them by $M_i$ 
$(i=1,\dots,d)$.  They are orientable if $M$ is orientable.  

\begin{defi} Let $M$ be a closed, connected, oriented, smooth manifold $M$ of 
dimension $2n$ with an effective smooth action of an $n$-dimensional 
torus group $T$ with non-empty fixed point set $M^T$. $M$ will be 
called a \emph{torus manifold} if a prefered orientation is 
given for each characteristic submanifold $M_i$. 
\end{defi}

A toric variety $X$ (of dimension $n$) is a normal complex 
algebraic variety of complex dimension $n$ with an effective algebraic action 
of $(\C^*)^n$ having a dense orbit.  
If $X$ is compact and non-singular, then $X$ with 
the restricted action of $T$ $(\subset (\C^*)^n)$ provides an 
example of a torus manifold of dimension $2n$.  
In this case, characteristic 
submanifolds are $(\C^*)^n$-invariant divisors. They 
have canonical orientations since they are complex manifolds. 
Similarly, when a torus manifold is equipped with a $T$-invariant unitary 
structure, characteristic submanifolds have canonical orientations. 
With these orientations of characteristic submanifolds, the torus manifold 
will be called a \emph{unitary torus manifold} (also called a unitary 
toric manifold in \cite{Masuda}). 

\begin{exam} A complex projective space $\C P^n$ 
with an action of $(\C^*)^n$ given by 
\[ [z_0,z_1,\dots,z_n]\to [z_0,g_1z_1,\dots,g_nz_n], \]
where $[z_0,z_1,\dots,z_n]\in \C P^n$ and $(g_1,\dots,g_n)\in 
(\C^*)^n$, is a compact non-singular toric 
variety.  This with the restricted $T$-action is a torus 
manifold and there are $n+1$ characteristic submanifolds, 
that are respectively defined by $z_i=0$ for $i=0,1,\dots,n$.  
\end{exam}

There are many torus manifolds which do not arise from 
compact non-singular toric varieties, see \cite{DJ}, \cite{Masuda}, 
\cite{OR}.  

Henceforth $M$ will denote a torus manifold of dimension $2n$. 
Let $p\in M^T$.  Since $M^T$ is isolated, the tangential 
$T$-representation $\tau_pM$ of $M$ at $p$ has no trivial factor, 
so it decomposes into the direct sum of $n$ irreducible real 
two-dimensional $T$-representations.  This implies that there are exactly 
$n$ characteristic submanifolds which contain $p$.  
In fact, an irreducible factor in $\tau_pM$ corresponds to the 
normal direction to a characteristic submanifold at $p$.  We set 
\[ \SigmaM:=\{ I\subset \{1,\dots,d\}\mid (\cap_{i\in I}M_i)^T
\not=\phi\}.\]
We add an empty set to $\SigmaM$ as a member, so that $\SigmaM$ 
becomes an \augmented.  The 
observation above implies that the cardinality of an element in 
$\SigmaM$ is at most $n$ and there is an element in $\SigmaM$ 
with cardinality $n$.

The \augmented\ $\SigmaM$ is closely related to the ring 
structure of the equivariant cohomology $H_T^*(M)$ of $M$ with integer 
coefficients. 
Let us explain this briefly.  Since $M_i$ and $M$ are 
oriented closed $T$-manifolds and the codimension of $M_i$ is two, 
the inclusion map from $M_i$ to $M$ induces a Gysin homomorphism 
$H_T^*(M_i)\to H_T^{*+2}(M)$ in equivariant cohomology which raises 
dgrees by two (see 
\cite{Kawakubo} for example).  Denote by $\xi_i\in H_T^2(M)$ 
the image of the identity element in $H_T^0(M_i)$.  
We may think of $\xi_i$ as the Poincar\'e dual of $M_i$ 
(considered as a cycle in $M$) in equivariant cohomology. 
If the orientation on $M$ or $M_i$ is reversed, then $\xi_i$ turns into 
$-\xi_i$.  

We take a polynomial ring 
$\Z[x_1,\dots,x_d]$ in $d$-variables and consider a map 
\[ \varphi\colon \Z[x_1,\dots,x_d]\to H_T^*(M)\]
which sends $x_i$ to $\xi_i$.  This map is often surjective. 
Here is a case. 

\begin{prop} {\rm (\cite{Masuda}, Proposition 3.4.)} 
\label{prop:ring_structure}
Suppose that $H^*(M)$ is generated by 
elements in $H^2(M)$ as a ring (this is the case when $M$ is a compact 
non-singular toric variety).  Then the map $\varphi$ is surjective and 
the kernel is the ideal generated by monomials $\prod_{i\in I}x_i$ 
for all subsets $I\subset \{1,\dots,d\}$ such that $I\notin \SigmaM$. 
In other words, $H_T^*(M)$ is isomorphic to the face ring (or 
Stanley-Reisner ring) of $\SigmaM$. 
\end{prop}

The equivariant cohomology $H_T^*(M)$ has a finer structure than 
the ring structure.  The map $\pi$ collapsing $M$ to a point induces 
a homomorphism $\pi^*\colon H^*_T(pt)=H^*(BT)\to H^*_T(M)$, so that 
$H^*_T(M)$ can be viewed as an algebra over $H^*(BT)$ through $\pi^*$.  
This algebra structure over $H^*(BT)$ cannot be determined by $\SigmaM$ 
and contains more information on the torus manifold $M$.  To see 
the algebra structure, it is enough to see the image of $H^2(BT)$ by 
$\pi^*$ because $H^*(BT)$ is a polynomial ring generated by elements 
in $H^2(BT)$.  

\begin{lemm} {\rm (\cite{Masuda}, Lemma 1.5.)} \label{lemm:edge_vector}
For each $i\in \{1,\dots,d\}$ there exists 
a unique element $v_i\in H_2(BT)$ such that 
\[ \pi^*(u)=\sum_{i=1}^d\langle u,v_i\rangle \xi_i 
\qquad\text{modulo $H^*(BT)$-torsions}
\]
for any $u\in H^2(BT)$, where $\langle\ ,\ 
\rangle$ denotes the usual pairing between cohomology and homology. 
\end{lemm}

\begin{proof} The proof is given in \cite{Masuda}, but we shall give 
a simple proof for the reader's convenience when $M$ is as in 
Proposition~\ref{prop:ring_structure}.  Since $H_T^2(M)$ is additively 
generated by $\xi_i$'s, one can express 
\[ \pi^*(u)=\sum_{i=1}^d v_i(u)\xi_i\]
with a unique integer $v_i(u)$ depending on $u$ for each $i$. We view 
$v_i(u)$ as a function of $u\in H^2(BT)$.  Since it is linear, 
it defines an element $v_i$ of $\Hom(H^2(BT),\Z)=H_2(BT)$ such that 
$v_i(u)=\langle u,v_i\rangle$. 
\end{proof}

\medskip
\noindent
{\bf Note.}   A geometrical interpretation of the vectors $v_i$ will 
be given in Section~\ref{sec:orbifolds}.
\medskip

In order to introduce a multi-fan, we adopt $H_2(BT)$ as the lattice 
$N$ and identify $H_2(BT;\R)$ with the vector space $N_\R=N\otimes \R$. 
Then we define a map 
\[ \LambdaM\colon \SigmaM\to \Cone(N)\]
by sending $I\in\SigmaM$ to the cone in $H_2(BT;\R)$ 
spanned by $v_i$'s $(i\in I)$ (and the empty set to $\{0\}$). 

Finally we shall define a pair of weight functions on maximal cones 
of dimension $n$.  Remember that a characteristic submanifold $M_i$ is 
a connected component of the set fixed pointwise by a certain circle 
subgroup, say $T_i$, of $T$.  
It turns out that $T_i$ agrees with the circle subgroup determined by 
$v_i\in H_2(BT)$ through the natural identification $H_2(BT)\cong 
\Hom(S^1,T)$ (\cite{Masuda}, Lemma 1.10).  
Therefore $M_I:=\cap_{i\in I}M_i$ is 
fixed pointwise by a subtorus $T_I$ generated by $T_i$'s for $i\in I$. 

\begin{lemm} \label{lemm:non-singularity}
{\rm (\cite{Masuda}, Lemma 1.7.)}  Suppose $I\in\SigmaM^{(n)}$.  Then 
the set $\{v_i\mid i\in I\}$ forms a basis of $H_2(BT)$, so that 
$M_I$ is a subset of $M^T$ and the cone $\LambdaM(I)$ is of dimension $n$. 
\end{lemm} 

A fixed point $p\in M^T$ belongs to $M_I$ for some $I\in
\Sigma^{(n)}$, and the tangent space $\tau_pM$ at $p\in M_I$ 
naturally decomposes into 
\[\tau_pM\cong \bigoplus_{i\in I}(\tau_pM/\tau_pM_i).\]
The orientations on $M$ and $M_i$ determine an orientation on 
$\tau_pM/\tau_pM_i$ for each $i$, and then an orientation on $\tau_pM$ 
through the above isomorphism.  
On the other hand, $\tau_pM$ 
has a given orientation since $M$ is oriented.  These two orientations 
on $\tau_pM$ may disagree.  We define the sign $\epsilon_p$ at $p$ to be
$+1$ or $-1$ according as the two orientations agree or disagree,
and set
\begin{align*} 
\wM^+(I):=&\text{ the number of $\{ p\in M_I\mid \epsilon_p=+1\}$},\\
\wM^-(I):=&\text{ the number of $\{ p\in M_I\mid \epsilon_p=-1\}$}.
\end{align*}
Note that $\wM^+(I)=1$ and $\wM^-(I)=0$ for all $I\in\Sigma^{(n)}$ if 
$M$ is a compact non-singular toric variety. 

\begin{defi} We call the triple $\DeltaM:=(\SigmaM,\LambdaM,\wM^\pm)$ 
the multi-fan of $M$.  
\end{defi}

A characteristic submanifold 
of $M_i$ is a connected component of $M_i\cap M_j$ for some $j$ containing 
a $T$-fixed point. We give it the orientation induced from those 
on $M_i$ and $M_j$.  
With these orientations equipped, $M_i$, on which $T/T_i$ acts 
effectively, is considered as a torus manifold.  If $M_i\cap M_j$ is 
connected for any $j\in\SigmaMsubi^{(1)}$ (this is the case when $M$ is a 
compact non-singular toric variety), then the multi-fan $\DeltasubMi$ of 
$M_i$ agrees with the projected multi-fan $\DeltaMsubi$ with respect to 
$\{i\}\in\SigmaM^{(1)}$. They are different otherwise but there is a 
natural surjective map from $\SigmasubMi$ to $\SigmaMsubi$.

Similarly, a connected component of $M_K$ for $K\in\SigmaM$ containing 
a $T$-fixed point is considered as 
a torus manifold, and $\DeltasubMK$ agrees with 
$\DeltaMsubK$ if $M_K$ and $M_K\cap M_j$ are connected for all 
$j\in\SigmaMsubK^{(1)}$, but otherwise they are different although 
there is a natural surjective map 
from $\SigmasubMK$ to $\SigmaMsubK$, where $\SigmasubMK$ is an 
\augmented\ obtained from the union of the simplicial sets associated 
with the connected components of $M_K$.  

The multi-fan $\DeltaM$ is non-singular by Lemma~\ref{lemm:non-singularity}.  
We shall show that it is complete. 

\begin{lemm} \label{lemm:complete}
$\DeltaM$ is complete. 
\end{lemm}

\begin{proof} As we remarked in Section~\ref{sec:multi-fan}
after the definition of the completeness of a multi-fan, 
it suffices to prove the pre-completeness of $\DeltaMsubJ$ 
for any $J\in\SigmaM^{(n-1)}$. Choose a generic vector $v$ from 
$N=H_2(BT)$.  The sign $(-1)^{\{i\}}$ for $i\in \SigmaMsubJ^{(1)}$ 
is defined as in Section~\ref{sec:DH} with respect to the projection 
image of $v$ on the quotient lattice of $N$ by the sublattice generated 
by $\LambdaM(J)\cap N$.  The pre-completeness of $\DeltaMsubJ$ 
is equivalent to the equality:
\[ \sum_{\{i\}\in\SigmaMsubJ^{(1)}}(-1)^{\{i\}}\wMsubJ(\{i\})=0,\]
which we will verify in the following.  Since $|J|=n-1$, a connected 
component of $M_J$ containing a $T$-fixed point is a $2$-dimensional 
sphere on which $T^J:=T/T_J$ acts effectively. We denote those connected 
components by $\Stwoalpha$'s.  They are torus manifolds equipped 
with the orientations 
discussed before this lemma.  Since $\Stwoalpha$ has two $T^J$-fixed points, 
$\SigmaS^{(1)}$ consists of two elements, denoted by $\alpha_\pm$, 
corresponding to the $T^J$-fixed points.  One easily checks that the 
multi-fan $\DeltaS$ of $\Stwoalpha$ is complete, which is equivalent to 
the equality:
\begin{equation}
(-1)^{\alpha_+}\wS(\alpha_+)+(-1)^{\alpha_-}\wS(\alpha_-)=0.\tag{9.1}
\end{equation}

As discussed before this lemma, 
we have a natural map $\pi_J\colon \SigmasubMJ\to 
\SigmaMsubJ$.  Note that if $\pi_J(\alpha_\epsilon)=\{i\}$ where 
$\epsilon$ stands for $+$ or $-$, then $(-1)^{\alpha_\epsilon}=
(-1)^{\{i\}}$.  On the other hand, we have 
\[ \wMsubJ(\{i\})=\sum_{\pi_J(\alpha_\epsilon)=\{i\}}\wS(\alpha_\epsilon).
\]
Therefore
\[ \sum_{\{i\}\in\SigmaMsubJ^{(1)}}(-1)^{\{i\}}\wMsubJ(\{i\})=
\sum_{\alpha_\epsilon}(-1)^{\alpha_\epsilon}\wS(\alpha_\epsilon),\]
which vanishes by (9.1), proving the lemma. 
\end{proof}

We make a remark on orientations at this point.  Choose an orientation on $T$ 
and fix it.  It induces an orientation on $H_2(BT;\R)$,  
so that $[\DeltaM]\in H_{n-1}(\SigmaM)$ is defined.  If the orientation 
on $T$ or $M$ is reversed, then $[\DeltaM]$ turns into $-[\DeltaM]$. 
But we have 

\begin{lemm} $[\DeltaM]$ does not depend on the 
orientations on $M_i$'s. 
\end{lemm}

\begin{proof} Recall that the cycle which defines $[\DeltaM]$ is 
$\sum_{I\in \SigmaM^{(n)}}\wM(I)\langle I\rangle$.  
We reverse the orientation on $M_i$.  Obviously, $\wM(I)$ and 
$\langle I\rangle$ remain unchanged unless $i\in I$.  Suppose $i\in I$. 
Then, since the orientation on $\tau_p M/\tau_p M_i$ is reversed, 
$\wM^+(I)$ and $\wM^-(I)$ will be interchanged, so that $\wM(I)$ 
turns into $-\wM(I)$.  As for $\langle I\rangle$, $\xi_i$ turns into 
$-\xi_i$ as remarked before and hence so does $v_i$ by 
Lemma~\ref{lemm:edge_vector}.  Thus, $\langle I\rangle$ turns into 
$-\langle I\rangle$ if $i\in I$.  After all, $\wM(I)\langle 
I\rangle$ does not depend on the orientations on $M_i$'s for 
any $I\in\SigmaM^{(n)}$. 
\end{proof}

Remember that there is a canonical isomorphism $\Hom(T,S^1)\cong H^2(BT)$. 
We denote by $t^u$ the element in $\Hom(T,S^1)$ 
corresponding to $u\in H^2(BT)$. 
Elements of $\Hom(T,S^1)$ are complex one-dimensional representations 
of $T$ and they generate the representation ring $R(T)$ of $T$ which 
is identified with the group ring of $H^2(BT)$. 
Since $\xi_i$ is the image of $1\in H^0_T(M_i)$ by the equivariant 
Gysin map from $M_i$ to $M$, its restriction to a $T$-fixed point 
$p$ in $M_i$, denoted by $\xi_i|_p$, gives the equivariant Euler class of 
the $T$-representation $\tau_pM/\tau_p M_i$; so $\tau_pM/\tau_p M_i
=t^{\xi_i|_p}$.  On the other hand, the equality in 
Lemma~\ref{lemm:edge_vector} restricted to 
$p$ shows that $\{\xi_i|_p\mid i\in I\}$ is the dual basis of 
$\{v_i\mid i\in I\}$, so $\xi_i|_p$ is independent of the choice of 
$p\in M_I$ and $\xi_i|_p=\u$ in the notation of 
Section~\ref{sec:Ehrhart_polynomial}.  Therefore we have
\[ \tau_pM=\bigoplus_{i\in I}t^{\u}\]
as a $T$-representation whenever $p\in M_I$.  

The elements $\xi_i$'s $(i=1,\dots,d)$ generate $H^2_T(M)$ additively 
modulo $H^*(BT)$-torsions (\cite[Lemma 3.2]{Masuda}) and the 
torsion elements vanish when restricted to the fixed point set 
$M^T$ because $H^*_T(M^T)$ is a free $H^*(BT)$-module.  Since the restriction 
$\xi_i|_p$ $(p\in M_I)$ depends on only $I$, we shall denote an element 
$\xi\in H^2_T(M)$ restricted to a point in $M_I$ by $\xi|_I$.  Note that 
\begin{equation*}
\xi_i|_I=\begin{cases}\u \quad&\text{if $i\in I$,}\\
0 &\text{otherwise.}\end{cases} \tag{9.2}
\end{equation*}
\begin{lemm} \label{lemm:index} For any $\xi\in H^2_T(M)$, 
\begin{equation*}
\sum_{I\in\SigmaM^{(n)}}\frac{\wM(I)t^{\xi|_I}}{\prod_{i\in I}
(1-t^{-\u})}
\end{equation*}
is an element of $R(T)$ when $M$ is a torus manifold. 
\end{lemm}

\begin{proof} Since $\xi_i$'s generate $H^2_T(M)$ additively modulo 
$H^*(BT)$-torsions, 
$\xi=\sum_{i=1}^d c_i\xi_i$ modulo $H^*(BT)$-torsions
with some integers $c_i$'s.  We define a map $\mathcal F_\xi\colon 
\SigmaM^{(1)}\to \HP(H^2(BT;\R))$ by 
\[ \mathcal F_\xi(\{i\}):=\{u\in H^2(BT;\R)\mid \langle u,v_i\rangle 
=c_i\}.
\]
The pair $(\DeltaM,\mathcal F_\xi)$ is a lattice multi-polytope, and 
$\cap_{i\in I}\mathcal F_\xi(\{i\})=\xi|_I$ for 
$I\in\SigmaM^{(n)}$ which follows from (9.2).   
Since $\DeltaM$ is non-singular by Lemma~\ref{lemm:non-singularity} and 
complete by Lemma~\ref{lemm:complete}, 
the lemma follows from Corollary~\ref{coro:key_corollary} applied to 
the multi-polytope $(\DeltaM,\mathcal F_\xi)$.
\end{proof}

\bigskip

\section{$T_y$-genus of a torus manifold} 
\label{sec:T_y_genus_of_a_torus_manifold}

When $M$ is a unitary torus manifold, the localization formula of the 
$T_y$-genus $T_y[M]$ of $M$ tells us that 
\begin{equation} T_y[M]=
\sum_{I\in\Sigma^{(n)}}\wM(I)\frac{\prod_{i\in I}(1+yt^{-\u})}
{\prod_{i\in I}(1-t^{-\u})}\tag{10.1}
\end{equation}
and this is actually a polynomial in $y$ with constant coefficients.  
As is well known, $T_0[M]$ agrees with the Todd genus of $M$ and 
$T_1[M]$ agrees with the signature of $M$, see \cite{Hirzebruch}. 
The $T_y$-genus is a genus for unitary manifolds 
and it is not defined for general torus manifolds.
But the right-hand side of (10.1) makes sense even for a torus 
manifold, and we take it 
as the definition of the $T_y$-genus $T_y[M]$ of $M$ and define 
the Todd genus of $M$ to be $T_0[M]$.  Note that the signature 
of $M$ is already defined for a torus manifold $M$ because 
$M$ is an oriented closed manifold, and that it agrees with 
$T_1[M]$ which follows from the Atiyah-Singer $G$-signature 
theorem.  

\begin{theo} \label{theo:T_y_genus of a torus manifold}
Let $M$ be a torus manifold of dimension $2n$.  Then 
\[ T_y[M]=T_y[\DeltaM]=\sum_{m=0}^n \e_{n-m}(\DeltaM)(-1-y)^m.\]
$($See Section~\ref{sec:T_y_genus} for $\e_q(\DeltaM)$.$)$
In particular, the Todd genus $T_0[M]$ of $M$ equals $\deg(\Delta)$. 
\end{theo}

\begin{proof} Look at the expansion of the right-hand 
side of (10.1) with respect to $y$. 
It follows from (9.2) and Lemma~\ref{lemm:index} that 
all coefficients of powers of $y$ in (10.1) are elements of $R(T)$.  
Take a generic vector $v\in H_2(BT)$ and evaluate the right-hand side 
of (10.1) on $v$.  
Then we get the following polynomial in $y$ whose coefficients are
Laurent polynomials in $z$:
\begin{equation}
\sum_{I\in\Sigma^{(n)}}\wM(I)\frac{\prod_{i\in I}(1+yz^{-\langle\u,v
\rangle})}{\prod_{i\in I}(1-z^{-\langle \u,v\rangle})}\tag{10.2}
\end{equation}
It is easily seen that (10.2) approaches to a polynomial in $y$ with 
constant coefficents if $z$ tends either to $0$ or to $\infty$. 
This means that (10.2) itself is a polynomial with constant coefficients. 
Since $v$ is generic, this implies that (10.1), that is $T_y[M]$, 
is actually a polynomial with constant coefficients equal to (10.2).
Then, by letting $z$ tend to $0$, we obtain 
\[ T_y[M]=\sum_{I\in \Sigma^{(n)}}\wM(I)(-y)^{\mu(I)},\]
where $\mu(I)=\sharp\{i\in I\mid \langle\u,v\rangle>0\}$. 
This $\mu(I)$ agrees with the $\mu(I)$ in 
Section~\ref{sec:T_y_genus} because $\{\u\mid i\in I\}$ is the dual 
basis of $\{ v_i\mid i\in I\}$. 
Hence $T_y[M]=T_y[\Delta(M)]$, proving the former equality in the theorem.
The latter follows from Corollary~\ref{coro:h=e}. 

As noted in the definition of $T_y[\Delta]$ in 
Section~\ref{sec:T_y_genus}, $T_0[\DeltaM]=\deg(\DeltaM)$. 
Since $T_0[M]=T_0[\DeltaM]$, the last statement in the theorem follows. 
\end{proof}

\begin{coro} \label{coro:signature} 
The signature $\Sign(M)$ of a torus manifold $M$ is given by 
\[ \Sign(M)=\sum_{m=0}^n(-2)^m e_{n-m}(\DeltaM).\]
If $T[M]=1$ and $w(M)(I)=1$ for all $I\in\Sigma(M)^{(n)}$, then 
$e_q(\DeltaM)$ agrees with the number of cones of dimension 
$q$ in $\DeltaM$.  
\end{coro}

\begin{proof} Since $\Sign(M)$ equals $T_1[M]$, 
the former statement follows from 
Theorem~\ref{theo:T_y_genus of a torus manifold}.  
The latter statement is 
noted in the definition of $\e_q(\Delta)$ in Section~\ref{sec:T_y_genus}. 
\end{proof}

\begin{rema} If $M$ is a compact non-singular toric variety, then 
$T[M]=1$ and $w(M)(I)=1$ for all $I\in\Sigma(M)^{(n)}$, and the 
formula above is already known in that case (\cite[Theorem 3.12(3)]{Oda}).
\end{rema}

\bigskip

\section{Equivariant index of a torus manifold} 
\label{sec:equivariant_index}

If $M$ is a unitary torus manifold, 
then the map $\pi$ collapsing $M$ to a point induces, in equivariant 
K-theory, the equivariant Gysin homomorphism 
\[ \pi_!\colon K_T(M)\to K_T(pt)=R(T).\]
If $E$ is a complex $T$-vector bundle over $M$, then $\pi_!(E)$ equals 
the index of a Dirac operator twisted by $E$.  It is sometimes called 
the equivariant Riemann-Roch number.  
The Todd genus of $M$ is equal to $\pi_!(1)$.  

Let $L$ be a complex $T$-line bundle over a unitary torus manifold $M$. 
Since $\pi_!(L)$ is an element of $R(T)$, one can express 
\begin{equation}
 \pi_!(L)=\sum_{u\in H^2(BT)}m_L(u)t^u
\tag{11.1}
\end{equation}
with integers $m_L(u)$ which are zero for all but finitely many elements $u$. 
In this section we describe the 
multiplicity $m_L(u)$ of $t^u$ in terms of the 
(shifted) moment map associated with $L$ when $M$ is a torus manifold.  
For that, we need to define $\pi_!(L)$ when $M$ is a torus manifold. 
This is done as follows.  When $M$ is a unitary torus manifold, 
the localization formula applied to $\pi_!(L)$ tells us that 
\begin{equation}
\pi_!(L)=\sum_{I\in\SigmaM^{(n)}}\frac{\wM(I)t^{c_1^T(L)|_I}}{\prod_{i\in I}
(1-t^{-\u})}\tag{11.2}
\end{equation}
where $c_1^T(L)\in H^2_T(M)$ denotes the equivariant first Chern class of $L$. 
(Note that $t^{c_1^T(L)|_I}$ is nothing but 
the complex one-dimensional $T$-representation 
obtained by restricting $L$ to a point in $M_I$.) 
The right-hand side of (11.2) is an element of $R(T)$ 
by Lemma~\ref{lemm:index} whenever 
$M$ is a torus manifold although $\pi_!$ may not be defined.  
Thus we define $\pi_!(L)$ as the right-hand side of (11.2) when 
$M$ is a torus manifold, and then define $m_L(u)$ as before 
using (11.1).  

In the following, we will make the following assumption on a torus 
manifold $M$, which is satisfied for 
compact non-singular toric varieties with restricted $T$-actions: 
\emph{all isotropy subgroups of $M$ are subtori of $T$ and each connected 
component fixed pointwise by a subtorus contains at least one 
$T$-fixed point}.  Then the union $\cup_{i=1}^d M_i$ is the set of 
points with nontrivial isotropy subgroups, and it follows from the 
slice theorem that the orbit space $M/T$ is a compact connected 
smooth manifold 
of dimension $n$ with $\cup_{i=1}^d M_i/T$ as boundary (after 
the corners are rounded).  

We make a further remark on orientations. 
The orbit space $M/T$ is orientable (see \cite{Masuda}, Lemma 6.7) 
and we orient it in such a way that the orientation on $T$ followed by 
that of $M/T$ agrees with that of $M$ times $(-1)^{n(n-1)/2}$. 
This determines a fundamental class in $H_n(M/T,\partial(M/T))$ and 
hence in $H_{n-1}(\partial(M/T))$, denoted by $[\partial(M/T)]$, 
through the boundary homomorphism from $H_n(M/T,\partial(M/T))$ to 
$H_{n-1}(\partial(M/T))$.  

Since $H^2_T(M)$ is additively generated by $\xi_i$'s $(i=1,\dots,d)$ 
modulo $H^*(BT)$-torsions, $c_1^T(L)=\sum_i c_i\xi_i$ 
modulo $H^*(BT)$-torsions with some integers $c_i$'s. 
Associated with $L$, there is defined 
the moment map $\Phi_L\colon M\to H^2(BT;\R)=L(T)^*$. 
It maps $M_i$ into an affine hyperplane 
$\{ u\in H^2(BT;\R)\mid \langle u,v_i\rangle=c_i\}$ for each $i$ 
(see \cite{Masuda}, Lemma 6.5).  
We slightly shift $\Phi_L$ so that the shifted map $\Phi_L'$ maps 
$M_i$ into 
\[ \FL'(\{i\}):=\{u\in H^2(BT;\R)\mid \langle u,v_i\rangle =c_i+
\frac{1}{2}\}\]
for each $i$.  In fact, $\Phi_L'$ is defined as follows.  
Let $K$ be a complex $T$-line bundle over $M$ with $c_1^T(K)=
-\sum_{i=1}^d\xi_i$.  Such $K$ exists (\cite{HY}).  When $M$ is 
a compact non-singular toric variety, $K$ is the canonical line 
bundle of $M$.  Using the moment map $\Phi_K\colon M\to H^2(BT;\R)$ 
associated with $K$, we define 
\[ \Phi_L':=\Phi_L-\frac{1}{2}\Phi_K.\]
The moment maps $\Phi_L$ and $\Phi_K$ are equivariant, the $T$-action 
on the target $H^2(BT;\R)$ being trivial; so $\Phi_L'$ induces a map 
\[ \bar\Phi_L'\colon M/T\to H^2(BT;\R).\]
The shifted affine hyperplanes $\FL'(\{i\})$'s miss the lattice 
$H^2(BT)$.  Since $\partial(M/T)=\cup_i(M_i/T)$ and $\bar\Phi_L'$ maps 
$M_i/T$ to $\FL'(\{i\})$ for each $i$, $\bar\Phi_L'$ induces a 
homomorphism 
\[ (\bar\Phi_L')_*\colon H_{n-1}(\partial(M/T))\to H_{n-1}(H^2(BT;\R)
\backslash\{u\})\]
for each lattice point $u\in H^2(BT)$.  We define 
\[ d_L'(u):=\text{ the mapping degree of $(\bar\Phi_L')_*$}\]
where the orientation on $H^2(BT;\R)$ is determined by that on $T$. 
Our main theorem in this section is the following. 

\begin{theo} \label{theo:multiplicity}
Let $M$ be a torus manifold.  Then $m_L=d_L'$ on $H^2(BT)$. 
\end{theo}

\begin{rema} This theorem was first established by Karshon-Tolman \cite{KT} 
when $M$ is a compact non-singular toric variety, and then extended to 
Spin$^c$ manifolds with 
torus actions by Grossberg-Karshon 
\cite{GK} and to a unitary torus manifold by the second named author \cite{Masuda}.
The family of torus manifolds contains these manifolds. 
\end{rema}

Let $\SM$ be the realization of the first barycentric subdivision of 
$\SigmaM$ and let $\SMsubi$ be the union of simplicies in $\SM$ which 
contain the vertex $\{i\}$ as in Section~\ref{sec:WN}.  
Since $\SMsubI=\cap_{i\in I}\SMsubi$ is 
contractible for any non-empty set $I\in {\SigmaM}$ and $\partial(M/T)=
\cup_{i=1}^d(M_i/T)$, it follows from Lemma~\ref{lemm:Psi} that there is 
a continuous map 
\[ \PsiM\colon \partial(M/T)\to \SM\]
sending $\cap_{i\in I}(M_i/T)$ to $\SMsubI$ for each $I\in \SigmaM$, 
and that such a map is unique up to homotopy preserving the 
stratifications, where the stratifications on $\partial(M/T)$ 
and $\SM$ mean subspaces $\cap_{i\in I}\partial(M_i/T)$ and $\SMsubI$ 
indexed by elements $I$'s in $\SigmaM$. 

If the orientation on $T$ or $M$ is reversed, then $[\partial(M/T)]$ 
and $[\DeltaM]$ will be multiplied by $-1$ simultaneously; so the 
following lemma makes sense.  

\begin{lemm} \label{lemm:rho}
$\PsiM_*([\partial(M/T)])=[\DeltaM]$. 
\end{lemm}

\begin{proof} We prove the lemma by induction on the dimension $n=
\dim(M/T)$.  When $n=1$, $M$ is $S^2$ with a nontrivial smooth 
$S^1$-action.  In this case, it is not difficult to check the lemma, 
which we leave to the reader.  

Assume that $n>1$.
Since a characteristic submanifold of $M_i$ is a connected component 
of $M_i\cap M_j$ for some $j$ and such $j$ is uniquely determined by 
the characteristic submanifold of $M_i$, there is a natural map 
$\pi_i\colon\SigmasubMi\to\SigmaMsubi$. 
This map is an isomorphism if $M_i\cap M_j$ is connected for any $j$, 
but otherwise it is only surjective.  
As we did in Lemma~\ref{lemm:Delta}, we identify 
the realization of $\SigmaMsubi$ with $\partial(\SMsubi)$. 
One sees that 
\begin{equation}
{\pi_i}_*([\DeltasubMi])=\sum_{i\in I\in \SigmaM^{(n)}}\wM(I)
\langle I\backslash\{i\}\rangle\in H_{n-2}(\partial(\SMsubi))
=H_{n-2}(\SigmaMsubi).
\tag{11.3}
\end{equation}

Since $M_i$ is itself a torus manifold, the spaces $\partial(M_i/T)$ 
and $\SsubMi$ have stratifications like for $M$, and hence 
we have a map $\PsiMi\colon \partial(M_i/T)\to \SsubMi$ 
preserving the stratifications. By the induction assumption 
\begin{equation}
\PsiMi_*([\partial(M_i/T)])=[\DeltasubMi] \in H_{n-1}(\SsubMi)=
H_{n-1}(\SigmasubMi). 
\tag{11.4}
\end{equation}

On the other hand, $\partial(\SMsubi)$ has a stratification 
induced from $\SM$ and each stratum is contractible.  
Since $\PsiM$ restricted to $\partial(M_i/T)$ is a map from 
$\partial(M_i/T)$ to $\partial(\SMsubi)$ preserving the stratifications and 
so is $\pi_i\circ\PsiMi$ as well, they are homotopic 
preserving the stratifications by Lemma~\ref{lemm:Psi}. 
Therefore, we have the following commutative diagram:
\[
\begin{CD}
H_{n-1}(\partial(M/T))@>\text{injective}>>\bigoplus_i H_{n-1}
(M_i/T,\partial(M_i/T))@>\cong>>\bigoplus_i H_{n-2}(\partial(M_i/T))\\
@V \PsiM_* VV @VVV @V \oplus{\pi_i}_*\PsiMi_* VV\\
H_{n-1}(\SM)@>\text{injective}>>\bigoplus_i H_{n-1}(\SMsubi,\partial
(\SMsubi))@>\cong>> \bigoplus_i H_{n-2}(\partial(\SMsubi))
\end{CD}
\]
where the left horizontal maps are natural ones. 
Tracing the upper horizontal sequence from the 
left to the right, $[\partial(M/T)]\in H_{n-1}(\partial
(M/T))$ maps to $\bigoplus_i [\partial(M_i/T)]$, 
and down to $\sum_{i\in I\in \SigmaM^{(n)}}\wM(I)\langle I\backslash 
\{i\}\rangle\in \bigoplus_i H_{n-2}(\partial(\SMsubi))$ by 
(11.3) and (11.4), while  
$[\DeltaM]\in H_{n-1}(\SM)$ maps through the lower horizontal 
sequence to the same element as observed in Lemma~\ref{lemm:Delta}.
Since the horizontal sequences above are injective, the lemma follows. 
\end{proof}

\begin{proof}[Proof of Theorem~\ref{theo:multiplicity}]  
By Lemma~\ref{lemm:Psi} we have a map $\SM\to H^2(BT;\R)$ associated with 
the multi-polytope $\PO_L':=(\DeltaM,\FL')$.  We denote the map by 
$\Psi_L'$.  The composition $\Psi_L'\circ \rho_M$ is a map from 
$\partial(M/T)$ to $H^2(BT;\R)$ sending $\cap_{i\in I}(M_i/T)$ to 
$\cap_{i\in I}\FL'(\{i\})$ for any $I\in\SigmaM$, and so is 
$\bar\Phi_L'$ as well.  Therefore, $\Psi_L'\circ \rho_M$ and $\bar\Phi_L'$ 
are homotopic preserving the stratifications by Lemma~\ref{lemm:Psi}.  
It follows from Lemma~\ref{lemm:rho} that 
\begin{align*}
d_L'(u)=&\text{ the mapping degree of }(\bar\Phi_L')_*: 
H_{n-1}(\partial(M/T))\to H_{n-1}(H^2(BT;\R)\backslash\{u\})\\
=&\text{ the mapping degree of }(\Psi_L'\circ\rho_M)_*: 
H_{n-1}(\partial(M/T))\to H_{n-1}(H^2(BT;\R)\backslash\{u\})\\
=&\text{ the mapping degree of }(\Psi_L')_*: 
H_{n-1}(\SM)\to H_{n-1}(H^2(BT;\R)\backslash\{u\})\\
=&\WN_{\PO_L'}(u)=\DHF_{\PO_L'}(u)=\DHF_{(\PO_L)_+}(u).
\end{align*}
This together with Corollary~\ref{coro:key_corollary} and the definition 
of $m_L$ (i.e., (11.1) and (11.2)) proves the theorem. 
\end{proof}

\bigskip


\section{Torus orbifolds}
\label{sec:orbifolds}

The aim of this section is to give the definition of a torus orbifold
and provide its basic properties for generalizing the results of 
Sections 10 and 11.
We first recall basic definitions concerning orbifolds. We
refer to \cite{Satake}, \cite{Kawasaki} or \cite{Duistermaat}
for details.
The reference \cite{LT} will be also useful; it deals with torus 
actions on symplectic orbifolds. 
If $M$ is an orbifold of dimension $n$, then there is a family
$\{(U_\alpha ,V_\alpha ,H_\alpha ,p_\alpha)\}$ of orbifold charts,
where $\{U_\alpha\}$ is an open covering of $M$, $V_\alpha$ is an
$n$-dimensional manifold, $H_\alpha$ is a finite group acting on
$V_\alpha$  and $p_\alpha :V_\alpha \to U_\alpha$ is a map which induces
a homeomorphism from $V_\alpha/H_\alpha$ onto $U_\alpha$. If $U_\alpha$
and $U_\beta$ intersect each other, then the charts
$(U_\alpha ,V_\alpha ,H_\alpha ,p_\alpha)$ and
$(U_\beta ,V_\beta ,H_\beta ,p_\beta)$ satisfy suitable compatibility
conditions. Such a family $\{(U_\alpha ,V_\alpha ,H_\alpha ,p_\alpha)\}$ is
called an orbifold atlas. For any point $x$ in $M$, there exists a
special type of orbifold chart 
$(U_x,V_x,H_x,p_x)$ with the property that
$p_x^{-1}(x)$ is a single point $\tilde{x}\in V_x$. The isomorphism
class of the group $H_x$ is uniquely determined by $x$ and is called
the isotropy group of $x$. The order of $H_x$, denoted by $d_x$, is
called the multiplicity of the point $x$. Such an orbifold chart
will be called a reduced orbifold chart. When $M$ is connected, the
minimum of the multiplicities is called the multiplicity of the
orbifold $M$ and is denoted by $d(M)$. The set
$\{x\in M\mid d_x=d(M)\}$ is open and dense in $M$. It is a manifold.
This set is called the principal stratum of the orbifold $M$.
We have $d(M)=1$ if and only if the actions of all the
isotropy groups are effective.   
\par
A map $f:M\to M'$ from an orbifold $M$ into another
orbifold $M'$ is called smooth if, near every point $x$ in $M$,
there is a homomorphism $\rho_\alpha:H_\alpha \to H'_\alpha$ and
a $\rho_\alpha$-equivariant smooth map $f_\alpha :V_\alpha \to V'_\alpha$ 
for suitable orbifold charts $(U_\alpha ,V_\alpha ,H_\alpha ,p_\alpha)$
for $M$ around $x$ and $(U'_\alpha ,V'_\alpha ,H'_\alpha ,p'_\alpha)$
for $M'$ around $f(x)$ satisfying the commutativity relation
$p'_\alpha \circ f_\alpha = f\circ p_\alpha$. A subset $M$ of an
orbifold $M'$ is called a suborbifold if, for each orbifold chart
$(U'_\alpha ,V'_\alpha ,H'_\alpha ,p'_\alpha)$ of $M'$,
$V_\alpha ={p'_\alpha}^{-1}(M\cap U'_\alpha)$ is an $H'_\alpha$-invariant
submanifold of $V'_\alpha$. If this is the case, $M$ becomes an
orbifold with orbifold charts $(U_\alpha ,V_\alpha ,H'_\alpha ,p'_\alpha)$
where $U_\alpha =M\cap U'_\alpha$, and the inclusion $M\to M'$ becomes
a smooth map. It may happen that $d(M)>d(M')$ ($M$ and $M'$ are assumed
connected).
The integer 
$d(M|M')=d(M)/d(M')$ will be called the relative multiplicity of
the pair $(M,M')$.
\par
Orbifold vector bundles
are also defined. Typical examples are the tangent bundle of an
orbifold and the normal bundle of a suborbifold. An orbifold is
orientable if its tangent bundle is orientable. If $E\to M$ is an
orbifold vector bundle over a connected orbifold, then the relative
multiplicity of the orbifold vector bundle $E$ is defined to be
$d(M|E)$ where $M$ is identified with the zero-section and is considered
as a suborbifold of $E$. If $M$ is a suborbifold of $M'$ and $\nu$ is
the normal bundle of $M$ in $M'$, then $d(M|\nu)$ equals $d(M|M')$.
\par 
Let $G$ be a Lie group. An action of $G$ on an orbifold $M$ is a
smooth map $\psi :G\times M \to M$ satisfying the usual rule of 
group action. Suppose that $G$ is connected. If $x\in M$ is a fixed
point of the action, and 
$(U_x,V_x,H_x,p_x)$ is a reduced orbifold chart around $x$ such
that $U_x$ is invariant under the action of $G$, then there is a
finite covering group $\tilde{G}_x$ of $G$ and an action of
$\tilde{G}_x$  on $V_x$ which covers the action of $G$ on $U_x$.
If $G$ is compact, the fixed point set of the action is a suborbifold.
\par
Now let $M$ be an oriented, closed orbifold of dimension $2n$ with
an effective action of an $n$-dimensional torus $T$.
A connected component of the fixed point set by a circle subgroup is
a suborbifold. A suborbifold of this type which has codimension two
and contains at least one fixed point of the $T$-action will be
called a \emph{characteristic suborbifold}. Let $M_i$ be a characteristic
suborbifold and $x\in M_i$. We take, as we may, a reduced orbifold
chart $(U_x,V_x,H_x,p_x)$ around $x$ such that $V_x$ is an open
disk in $\mathbb{R}^{2n}$ and the action of $H_x$ on $V_x$ is linear.
We denote by the same symbol $V_x$ the tangent space to $V_x$
at the point
$\tilde{x}=p_x^{-1}(x)$. Then the vector space $V_x$ decomposes
into a direct sum $V_{ix}\oplus V_{ix}^\perp$ where 
$V_{ix}^\perp$ is tangent to $p_x^{-1}(U_x\cap M_i)$, and 
the vector space $V_{ix}$ represents the fiber direction of
the normal bundle of $M_i$ in $M$.
The isotropy group $H_x$ acts on $V_{ix}$.
\begin{lemm}\label{lemm:Si}
Let $M$ be an oriented closed orbifold
as above and $M_i$
a characteristic suborbifold. Let $S_i$ denote the circle subgroup
which fixes the points of $M_i$. Then there exists a finite covering
group $\tilde{S}_i$ of $S_i$ and a lifting of the action of $S_i$
to the action of $\tilde{S}_i$ on $V_x$ for any point $x\in M_i$.
The lifted action of $\tilde{S}_i$ preserves $V_{ix}$.
\end{lemm}
\begin{proof}
To $x\in M_i$ we correspond the degree of the minimal finite covering
$\tilde{S}_{ix}$ of $S_i$ such that there is a lifting of the action
to
$\tilde{S}_{ix}$. The lifted action necessarily preserves $V_{ix}$.
It is not difficult to see that the correspondence is locally
constant. 
Since $M_i$ is connected the correspondence must be constant.
\end{proof}
Hereafter we denote by $\rho_i:\tilde{S}_i\to S_i$ the minimal
finite covering of
$S_i$ with the above property. $\tilde{S}_i$ acts effectively on
$V_x$.
\par
An oriented, closed orbifold $M$ of dimension $2n$ with an effective
action of a torus $T$ of dimension $n$ with non-empty fixed
point set $M^T$ equipped with a preferred
orientation of the normal bundle of each characteristic
suborbifold will be called a \emph{torus
orbifold} if, for each $M_i$ and at each point $x\in M_i$, the
action  of
$H_x$ preserves the orientation of each $V_{ix}$. 
Note that choosing an orientation of a characteristic submanifold 
is equivalent to choosing an orientation of its normal bundle. 
Thus a torus manifold is a torus orbifold in the above sense. 
Another example is a unitary torus orbifold. 
A unitary torus
orbifold is a torus orbifold such that $V_\alpha$ is a unitary
manifold, the action of $H_\alpha$ preserves the unitary structure
of $V_\alpha$ for each orbifold chart
$(U_\alpha ,V_\alpha ,H_\alpha ,p_\alpha)$ and the action of $T$ on
$M$ also preserves the unitary structure of $V_\alpha '$s.
\par
Let $M$ be a torus orbifold. The preferred orientation of
the normal bundle
$\nu_i$ of $M_i$ makes it a complex orbifold line bundle. Then
there is a
unique isomorphism $\varphi_i:S^1\to \tilde{S}_i$ such that
$\varphi(z)$ acts by the complex multiplication of
$z$ on each $V_{ix}$. We identify $\tilde{S}_i$ with $S^1$ via
$\varphi_i$. The homomorphism $\rho_i:S^1=\tilde{S}_i\to T$ defines
an element $v_i\in \Hom(S^1,T)=H_2(BT;\Z)$. 
We are now ready to define the multi-fan
$\DeltaM=(\SigmaM,\Lambda(M),w(M)^\pm)$ associated with a torus
orbifold $M$ in an entirely similar way to the case of torus
manifolds. Specifically 
\[ \SigmaM=\{I\mid (\cap_{i\in I}M_i)^T\not= \emptyset\}, \]
and $\Lambda(M)(I)$ is the cone in $H_2(BT;\R)$ with apex at $0$
and spanned by $\{v_i\mid i\in I\}$. Furthermore
$w(M)^\pm(I)=\#\{x\in M_I\mid \epsilon_x=\pm 1\}$ for
$I \in \SigmaM^{(n)}$, where
$\epsilon_x$ is defined to be the ratio of two orientations at $x$,
one which is given by the orientation of $M$ and the other by
that of the oriented vector space $V_x=\oplus_{i\in I} V_{ix}$.
\par
We set $\tilT_I=\prod_{i\in I}\tilSi$ for $I\in \SigmaM^{(k)}$
and $\rho_I=\prod_{i\in I}\rho_i:\tilT_I\to T$. The image of $\rho_I$
is denoted by $T_I$. $\rho_I:\tilT_I\to T_I$ is a finite covering. 
$T_I$ fixes the points of $M_I=\bigcap_{i\in I}M_i$. If $I\in
\SigmaM^{(n)}$, then $T_I=T$.  
Let $x$ be a fixed point of the action of $T$ on $M$. Then there is
a unique $I\in \SigmaM^{(n)}$ such that $x$ belongs to 
$M_I$. The inclusion
$S^1=\tilde{S}_i\to \tilde{T}_I$ defines an element
$\tilde{v}_i\in \Hom(S^1,\tilde{T}_I)=H_2(B\tilde{T}_I;\Z)$, and
we have $\rho_{I*}(\tilde{v}_i)=v_i$. 
$V_x$ and $V_{ix},\ i\in I,$ are complex
$\tilde{T}_I$-modules, and the decomposition 
$ V_x = \oplus_{i\in I} V_{ix}$ is compatible with the action
of $\tilde{T}_I$.
The effectiveness of the $T$-action on $M$ implies that
$\tilde{T}_I$ effectively acts on $V_x$; 
equivalently, it implies that $\{\tilde{v}_i\mid i\in I\}$
is a basis of $H_2(B\tilde{T}_I;\Z)$. Since $\rho_{I*}:
H_2(B\tilde{T}_I;\Z)\to H_2(BT;\Z)$ is injective,
the $v_i,\ i\in I,$ are linearly independent
in $H_2(BT;\R)$.
\begin{lemm}\label{lemm:complete_for_orbifold}
$\DeltaM$ is a complete multi-fan.
\end{lemm}
\begin{proof}
The argument is almost similar to the case of torus manifolds.
One has only to observe that the characteristic suborbifolds and
their intersections are torus orbifolds and a 2-dimensional torus
orbifold is topologically a 2-sphere acted on by a circle group with
exactly two fixed points.
\end{proof}
\par

\begin{lemm}\label{lemm:isotropy}
Suppose $d(M)=1$. Let $I\in \SigmaM^{(k)}$, and let $x$ be a point
in the principal stratum (as an orbifold) of $M_I$. Then the
isotropy group $H_x$ of $x$ is isomorphic to the kernel of
$\rho_I:\tilT_I \to T$. 
\end{lemm}

\begin{proof}
Let $(U_x,V_x,H_x,p_x)$ be an orbifold chart around $x$. We may
regard $V_x$ as an $n$ dimensional $\tilT_I$-module as before.
As such, $V_x$ is decomposed as a direct sum of $\tilT_I$-modules
\[ V_x=\left(\oplus_{i\in I}V_{ix}\right)\oplus V' \]
where $V'$ is projected into $M_I$ by $p_x$. 
$\tilT_I=\prod_{i\in I}\tilSi$ can be
regarded as embedded in the general linear group of 
$\oplus_{i\in I}V_{ix}$. Since $H_x$ acts on each $V_{ix}$
preserving its orientation, there is a homomorphism $H_x\to \tilT_I$.
The action of $H_x$ on $V'$ is trivial. Moreover the action of $H_x$
on $V_x$ is effective because $d(M)=1$. It follows that the 
homomorphism above embeds $H_x$ into $\tilT_I$. Since the kernel 
of $\rho_I$ is equal to the intersection of $\tilT_I$ with the image 
of $H_x$, it is isomorphic to $H_x$.  
\end{proof} 

It is known that a closed oriented orbifold $M$ of dimension $n$ has
the fundamental
class $[M]\in H_n(M;\Z)$, and that the Poincar\'{e} duality holds, i.e.,
the operation $ \vartheta=[M]\cap \ :H^q(M;\Q)
\to H_{n-q}(M;\Q)$ is an isomorphism. If $f:M\to M'$ is a smooth map from
an oriented close orbifold $M$ to another such $M'$, then the
Gysin homomorphism $f_!:H^q(M;\Q)\to H^{q+n-n'}(M';\Q)$ is defined
to be the
compostion $\vartheta^{-1}\circ f_*\circ \vartheta$, where $n'$ is
the dimension of $M'$. If a compact Lie group $G$ acts
on $M$ and $M'$, and $f$ is equivariant, then the equivariant Gysin
homomorphism $f_!:H_G^q(M;\Q)\to H_G^{q+n-n'}(M';\Q)$ is also defined.
\par
Henceforth $M$ will be a torus orbifold. For each $i\in \SigmaM^{(1)}$,
we set
\[ \xi_i = (f_i)_!(1)\in H_T^2(M;\Q), \]
where $f_i:M_i\to M$ is the inclusion.
\begin{lemm}\label{lemm:xi}
Let $c_1^T(\nu_i)$ be the equivariant first Chern class of the
normal bundle $\nu_i$. Then we have
\[ c_1^T(\nu_i)=f_i^*(\xi_i). \]
\end{lemm}
\begin{proof}
We may assume that $d(M)=1$. 
Take an equivariant Thom form $\phi$ for the equivariant
orbifold bundle
$\nu_i$ (we refer to \cite{Bott} for Thom form and Chern form, cf. also \cite{Duistermaat}).
Let $x$ be a point in the principal stratum of $M_i$,
and $(U_x,V_x,H_x,p_x)$ a reduced orbifold chart around $x$. 
The restriction of $\phi$ to $V_x$ is invariant under the action
of $H_x$ and its
support is contained in a tubular neighborhood $W_i$ of 
$V_i=p_x^{-1}(U_i)$,
where $U_i=U_x\cap M_i$.
Moreover, with respect to the fibering 
$\tilde{\pi}_i:W_i\to V_i$, we have
$|H_x|^{-1}(\tilde{\pi}_i)_*(\phi)=1$,
where $(\tilde{\pi}_i)_*$ is the
integration along the fiber of $\tilde{\pi}_i$.
Note that the fiber is $V_{ix}$,
and that the action of $H_x$ preserves the orientation of $V_{ix}$.
The equivariant Chern class $c_1^T(\nu_i)$ is the restriction
to $M_i$ of the cohomology class $[\phi]$ of $\phi$. Here
$[\phi]$ is considered as a relative class
in $H_T^2(W,W\setminus M_i;\R)$ where $W$ is a tubular neighborhood
of $M_i$.
\par
On the other hand, $\xi_i$ is the restriction of a cohomology class
$\psi\in H_T^2(W,W\setminus M_i;\R)$ such that 
\[ \pi_*(\psi) =1\in H_T^0(W;\R)=H_T^0(M_i;\R), \]
where $\pi:W\to M_i$
denotes the projection of the fibration. Note that the fiber of $\pi$
is $U_{ix}=V_{ix}/H_x$, where $H_x$ acts effectively
on $V_{ix}$. We have
\[
\pi_*([\phi])=|H_x|^{-1}(\tilde{\pi}_i)_*([\phi])
             =1
             =\pi_*(\psi).
\]
But $\pi_*$ is an isomorphism (Thom isomorphism). Hence we have
$ [\phi]=\psi $, and consequently
\[ c_1^T(\nu_i)=[\phi]|M_i=\psi|M_i=f_i^*(\xi_i). \] 
\end{proof}
\par
We noticed that, for $I\in\SigmaM^{(n)}$, 
$\{v_i\mid i\in I\}$ was a basis of $H_2(BT;\R)$.
Let $\{u_i^I\}$ be the dual basis in $H^2(BT;\R)$. This can be
interpreted in the following way. Let $\{\tilde{u}_i\mid i\in I\}$
be the basis of $H^2(B\tilde{T}_I;\Z)$ dual to $\{\tilde{v}_i\mid i\in I\}$.
We have $\rho_I^*(u_i^I)=\tilde{u}_i$, since
$\rho_{I*}(\tilde{v}_i)=v_i$. We identify $H^2(B\tilde{T}_I;\R)$
with $H^2(BT;\R)$ by the isomorphism $\rho_I^*$. Then
$H^2(B\tilde{T}_I;\Z)$ can be considered as embedded in
$H^2(BT;\R)$. With this convention we have $u_i^I =\tilde{u}_i$.
\par
Let $x\in M^T$ be a fixed point of the $T$-action. In the sequel
we identify $H_T^2(x;\R)$ with $H^2(BT;\R)$.
\begin{lemm}\label{lemm:ui}
Let $I\in \SigmaM^{(n)}$ and $x\in M_I$. Then 
$\xi_i|x=u_i^I\in H^2(BT;\R)$ for $i\in I$. If $j\notin I$, then
$\xi_j|x=0$.
\end{lemm}
\begin{proof}
By Lemma \ref{lemm:xi} we have 
\[ \xi_i|x=c_1^T(\nu_i|x) .\] 
But $\nu_i|x$  viewed as $\tilde{T}_I$-module is $V_{ix}$. It follows
that $c_1^{\tilde{T}_I}(\nu_i|x) =\tilde{u}_i$. Hence
\[ c_1^T(\nu_i|x) =u_i^I. \]
If $j\notin I$, then $x\notin M_j$. Therefore $\xi_j|x=0$.
\end{proof}
\par
If we consider $u_i^I=\tilde{u}_i$ as an element of 
$\Hom(\tilT_I,S^1)=H^2(B\tilT_I;\Z)$, then Lemmas~\ref{lemm:ui} and 
\ref{lemm:orbiedge_vector} imply that 
$u_i^I$ is nothing but the $\tilT_I$-module $V_{ix}$.
The following Lemma describes the algebra structure of $H_T^*(M;\R)$
over $H^*(BT;\R)$ modulo $H^*(BT;\R)$-torsion as in the case of
torus manifolds (Lemma \ref{lemm:edge_vector}).
\begin{lemm}\label{lemm:orbiedge_vector}
The following equality holds for any $u\in H^2(BT;\R)$:
\[ \pi ^*(u)=\sum_{i\in \SigmaM^{(1)}}\langle u,v_i\rangle \xi_i\quad
\text{modulo $H^*(BT;\R)$-torsion} .\]
\end{lemm}
\begin{proof}
Let $x\in M_I\subset M^T$ be a fixed point of the $T$-action. We restrict
both sides of the equality in Lemma \ref{lemm:orbiedge_vector} to $x$.
On the left hand side we get $u$. On the right hand side the result is
\[ \sum_{i\in I}\langle u, v_i\rangle u_i^I \]
by virtue of Lemma \ref{lemm:ui}. But this is equal to $u$ by the
definition of the $u_i^I$. Thus both sides coincide after the restriction
to each $x\in M^T$. Since the restriction homomorphism
$\pi^* :H_T^*(M;\R) \to H_T^*(M^T;\R)$ is injective modulo
$H^*(BT;\R)$-torsion, the equality is confirmed.
\end{proof}
\begin{rema}
The equality in Lemma \ref{lemm:orbiedge_vector} characterizes the vectors
$v_i$ in terms of the $\xi_i$ as in Lemma \ref{lemm:edge_vector}.
\end{rema}
\par
We set $N=H_2(BT;\Z)$ and define $N_I$ for
$I\in \SigmaM^{(n)}$ to be the lattice generated by the
$v_i,\ i\in I$.
\begin{lemm}\label{lemm:N/NI}
Assume that $d(M)=1$. Let $x\in M_I$ with $I\in \SigmaM^{(n)}$.
Then $H_x$ is isomorphic to $\Ker \rho_I$. Moreover $\Ker \rho_I$
is isomorphic to $N/N_I$.
\end{lemm}
\begin{proof}
We have
already shown that $H_x$ is isomorphic to the kernel of $\rho_I$
in Lemma \ref{lemm:isotropy}.
For the second part it suffices to note that $N$ and $N_I$ can be
identified with the fundamental group of $T$ and $\tilde{T}_I$.
Therefore the kernel of $\rho_I$ is isomorphic to $N/N_I$.
\end{proof}
\begin{rema}
Hereafter we identify $H_x$ and $N/N_I$ with 
$\Ker \rho_I\subset \tilde{T}_I$ through the isomorphisms
given in Lemma \ref{lemm:N/NI}. 
We put $\chi_I(u,v)=\exp(2\pi\sqrt{-1}\langle u,v\rangle)$
for $u\in H^2(B\tilde{T}_I;\Z)$ and $v\in H_2(BT;\R)$. If $u$
is fixed, then the value $\chi_I(u,v)$ depends only on the
equivalence class of $v$ modulo $N_I$.
Hence, if we identify $\tilde{S}_i$ with
$S^1$ via $\varphi_i$ as before and $\tilde{T}_I$
with $\prod_{i\in I}S^1$ via $\prod_{i\in I}\varphi_i$, then the map 
$\exp :H_2(BT;\R)\to \tilde{T}_I$ defined by
$\exp(v)=\prod_{i\in I}\exp(2\pi\sqrt{-1}\langle u_i^I,v\rangle)$
is a universal covering map and its kernel is $N_I$. It induces
an isomorphism from $H_x=N/N_I$ onto $\Ker \rho_I$. We shall
write $\chi_I(u,g)$ instead of $\chi_I(u,v)$ for
$g=\exp(v)\in \tilde{T}_I$ as in Section \ref{sec:Ehrhart_polynomial}.
Let $V$ be a one dimensional $\tilde{T}_I$-module. It defines
an element $u\in \Hom(\tilde{T}_I ,S^1)=H^2(B\tilde{T}_I;\Z)$.
Then the action of $g\in \tilde{T}_I$
on $V$ is given by the complex multiplication by
$\chi_I(u,g)$.  
\end{rema}
\par
Suppose that $M$ is a unitary torus orbifold such that $d(M)=1$.
Let $L$ be a
$T$-invariant complex line bundle over $M$. By using the
hermitian connection of $M$ and a hermitian
connection of $L$, a Dirac operator twisted by $L$ is defined as in 
the case of torus manifolds. Its
index is a $T$-module. It is called the equivariant Riemann-Roch
number with coefficient in $L$, and is denoted by $RR^T(M,L)\in R(T)$.
It can be expressed by the fixed point
formula due to Vergne \cite{Vergne}; cf. also \cite{Duistermaat}.
The formula
is particularly simple when all the fixed points are isolated. 
It is convenient to write down the image of $RR^T(M,L)$ by
$\ch :R(T)\to H^{**}(BT;\R)$; the result is 
\begin{lemm}\label{lemm:RR}
Let $\xi=c_1^T(L)$ be the equivariant Chern class of $L$.
Then 
\[ \ch (RR^T(M,L))=\sum_{x\in M^T}\frac{\epsilon_x e^{\xi|x}}{|H_x|}
\sum_{g\in H_x}\frac{\chi_{I_x}(\xi|x,g)}{\prod_{i\in I_x}
(1-\chi_{I_x}(u_i^{I_x},g)^{-1}e^{-u_i^{I_x}}) }\ , \]
where $I_x\in \SigmaM^{(n)}$ is such that $x\in M_{I_x}$.
\end{lemm}
It can be shown that, if $x$ and $y$ both lie in the same $M_I$,
then $\xi|x=\xi|y$ for $\xi=c_1^T(L)$. The proof is same as in the
case of torus manifolds as was given in \cite{Masuda}. We shall
write $u_I=c_1^T(L)|x$ for $x\in M_I$. Taking Remark below Lemma
\ref{lemm:N/NI} in account, we get
\begin{prop}\label{prop:RR}
\[ \ch (RR^T(M,L))=\sum_{I\in \SigmaM^{(n)}}\frac{w(M)(I)e^{u_I}}{|N/N_I|}
\sum_{g\in N/N_I}\frac{\chi_I(u_I,g)}{\prod_{i\in I}
(1-\chi_I(u_i^I,g)^{-1}e^{-u_i^I} )}\ . \]
\end{prop}
Since $\ch:R(T)\to H^{**}(BT;\R)$ is injective, the formula in 
Proposition \ref{prop:RR} characterizes $RR^T(M,L)$. Using the
notation in Section \ref{sec:Ehrhart_polynomial}, we obtain
\begin{coro}\label{coro:RR}
\[ RR^T(M,L)=\sum_{I\in \SigmaM^{(n)}}\frac{w(M)(I)t^{u_I}}{|N/N_I|}
\sum_{g\in N/N_I}\frac{\chi_I(u_I,g)}{\prod_{i\in I}
(1-\chi_I(u_i^I,g)^{-1}t^{-u_i^I} )}\ . \]
\end{coro}
When $u_I=c_1^T(L)|x,\ x\in M_I,$ lies in $N^*=H^2(BT;\Z)$, 
then $\chi_I(u_I,g)=1$ for all $g\in N/N_I$. Therefore, if
$u_I\in N^*$ for all $I\in \SigmaM^{(n)}$, then 
\[ RR^T(M,L)=\sum_{I\in \SigmaM^{(n)}}\frac{w(M)(I)t^{u_I}}{|N/N_I|}
\sum_{g\in N/N_I}\frac{1}{\prod_{i\in I}
(1-\chi_I(u_i^I,g)^{-1}t^{-u_i^I} )}\ . \]
By observing that $g\mapsto \chi_I(u,g)$ is a character of $N/N_I$ for
any $u\in H^2(B\tilT_I,\Z)=\Hom(\tilT_I,S^1)$, 
the formula above can be rewritten
in the following form:
\begin{equation*}
RR^T(M,L)=\sum_{I\in \SigmaM^{(n)}}\frac{w(M)(I)t^{u_I}}{|N/N_I|}
\sum_{g\in N/N_I}\frac{1}{\prod_{i\in I}
(1-\chi_I(u_i^I,g)t^{-u_i^I} )}\ .\tag{12.1}
\end{equation*}
The right hand side of this formula 
(12.1) appeared
in Corollary \ref{coro:key_corollary}. There, it was related to a
lattice multi-polytope $\mathcal{P}$, in which $u_I$ is contained in
$N^*$ for all $I\in \Sigma^{(n)}$, and the Duistermaat-Heckman function
$DH_{\mathcal{P}_+}$. Suppose that $c_1^T(L)$ is of the form 
$c_1^T(L)=\sum_{i\in \SigmaM^{(1)}}c_i\xi_i\in H_T^2(M;\R)$. Then 
the above multi-polytope
$\mathcal{P}$
is nothing but the one whose first Chern class is 
$c_1(\mathcal{P})=\sum c_ix_i$ as in Section 8.
Note that $\PO$ is not always a lattice multi-polytope in this case. 
\begin{rema}
Corollary \ref{coro:key_corollary} shows that the right hand side of
the formula (12.1) 
depends only on $\DeltaM$ and
$\mathcal{P}$; namely, it does not depend on the choice of
generating vectors  $v_i\in H^2(BT;\R)$ in so far as they
lie in $N=H^2(BT;\Z)$ and $\{u_i^I\mid i\in I\}$ is interpreted as the
dual basis of $\{v_i\mid i\in I\}$.
\end{rema}
When 
$M$ is a torus manifold, the Duistermaat-Heckman function has
a geometric meaning coming from the moment map of the line bundle
$L$ as was explained in Section \ref{sec:equivariant_index}. 
There the role of the winding number was also explained. 
These notions are generalized to the case of
torus orbifolds and similar results hold in this case too. The
details can be worked out without much alteration and are left
to the reader.

The $T_y$-genus of a torus orbifold is also defined by using the 
fixed point formula due to Vergne in a similar way as in Section 10,
and the analogue of Theorem 10.1 holds. the details are left to
the reader.
\bigskip
\section{Realizing multi-fans by torus orbifolds} \label{sect:realization}

In the previous section, we associated a complete simplicial multi-fan 
of dimension $n$ with a torus orbifold of dimension $2n$.  
In this section, we consider the converse problem. 
If a multi-fan $\Delta$ is associated with a torus orbifold $M$, then 
we say that $\Delta$ is \emph{(geometrically) realized} by $M$, or 
$M$ \emph{realizes} $\Delta$. 

We recall how the multi-fan of $M$ changes when 
the orientations on $M$ or $M_i$ are reversed.  
If the orientation of $M$ is unchanged but that of $M_i$ is reversed, 
then the orientation of the normal bundle of $M_i$ is reversed and,
hence, $1$-dimensional cone $C(i)$ tunrs into the cone $-C(i)$,
and the pair $(w(M)^+(I),w(M)^-(I))$
turns into $(w(M)^-(I),w(M)^+(I))$ for $I\in\SigmaM^{(n)}$ containing 
$i$ while others remain unchanged.  
If the orientations of $M$ and of all the $M_i$'s are reversed, 
then all the cones $C(i)$'s remain unchanged but
$(w(M)^+(I),w(M)^-(I))$ turns into 
$(w(M)^-(I),w(M)^+(I))$ for any $I\in\SigmaM^{(n)}$ so that $w(M)(I)$ 
turns into $-w(M)(I)$ for any $I\in\SigmaM^{(n)}$.  The torus orbifold 
$M$ with the reversed orientations of $M$ and all the $M_i$'s will be denoted 
by $-M$. 

The underlying space of a torus orbifold of dimension $2$ is $S^2$ with 
the standard $S^1$-action.  In this case, there are two characteristic 
submanifolds.  They are $S^1$-fixed points.  Taking orientations 
on $S^2$ and its characteristic submanifolds into account, we easily have 
the following theorem.  

\begin{theo} \label{theo:2-dim torus orbifold}
A complete simplicial multi-fan $\Delta=(\Sigma, C, w^\pm)$ of dimension $1$ 
is geometrically realized if and only if $\Sigma$ is isomorphic to 
the argumented simplicial set obtained from the boundary of a $1$-simplex and 
$\{w^+(I),w^-(I)\}=\{1,0\}$ as a set for $I\in\Sigma^{(1)}$. 
\end{theo}

The analysis of a torus orbifold of dimension $4$ is more 
complicated.  In this case, each characteristic suborbifold is 
homeomorphic to $S^2$ and has two fixed points. 
Therefore, if two of the characteristic suborbifolds intersect, then 
they intersect at one point or 
two points, and if they intersect at two points, then they do not intersect 
at any other characteristic suborbifolds.  
We also note that a $T$-fixed point is an intersection of 
two characteristic suborbifolds. 
These facts imply the \lq\lq only if" 
part in the following theorem.  We will prove the \lq\lq if" part later. 

\begin{theo} \label{theo:4-dim torus orbifold}
A complete simplicial multi-fan $\Delta=(\Sigma, C, w^\pm)$ of dimension $2$ 
is geometrically realized if and only if the following two conditions are 
satisfied for each $I\in \Sigma^{(2)}$:
\begin{enumerate}
\item $\{w^+(I),w^-(I)\}=\{1,0\}$ or $\{1,1\}$,
\item when $\{w^+(I),w^-(I)\}=\{1,0\}$, 
there are exactly two elements, say $I'$ and $I''$, in $\Sigma^{(2)}$ 
such that $I\cap I'$ and $I\cap I''$ 
are in $\Sigma^{(1)}$ and $I\cap I' 
\cap I''=\emptyset$, and when 
$\{w^+(I),w^-(I)\}=\{1,1\}$, there is no element 
$I'\in \Sigma^{(2)}$ such that $I\cap I'\in\Sigma^{(1)}$.
\end{enumerate}
\end{theo}

In contrast to the low dimensional cases above, we have 

\begin{theo} \label{theo:higher dim torus orbifold}
Any complete simplicial multi-fan of dimension $\ge 3$ is geometrically 
realized.
\end{theo} 

In the following $\Delta=(\Sigma,\Lambda,w^\pm)$ will be a complete 
simplicial multi-fan of dimension $n\ge 2$ unless otherwise stated. 
Here is an outline of how to realize $\Delta$ by a torus orbifold. 
We choose and fix a generic (rational) $1$-dimensional cone 
in $N_\R$, and decompose $\Delta$ using it 
into a number of what we call \emph{minimal} multi-fans.  
Minimal multi-fans can essentially be realized by weighted projective
spaces. We paste them together by performing 
equivariant connected sum along characteristic suborbifolds 
and at $T$-fixed points to obtain a desired torus orbifold realizing 
the given $\Delta$.  

Equivariant connected sum is performed through two isomorphic orbifold
charts. In this way attention should be payed to orbifold structures.
So we make a remark on orbifold structures at this point.  
There are many choices of an orbifold structure on $M$ (e.g. 
$S^2$ with the standard $S^1$-action admits infinitely many orbifold 
structures), but the associated multi-fan does not depend on the choice 
of an orbifold structure.  In fact, the circle subgroup $\Si$
determined by 
the vector $v_i$ in the previous section is the one which fixes points 
in the 
characteristic suborbifold $M_i$, so the line generated by $v_i$ is 
independent of the orbifold structure.  Moreover the direction 
of $v_i$ is determined by the choice of orientations on $M$ and $M_i$, 
so the cone spanned by $v_i$ is independent of the orbifold 
structure.  What depends on the orbifold structure is the 
length of $v_i$ which is equal to the degree of the covering map
$\tilSi \to \Si$. In this way the vectors $v_i$ reflect the orbifold
structure related to the torus action. We shall call the vector $v_i$
the edge vector of the $1$-dimensional cone $C(i)$.

We shall use two types of 
equivariant connected sum; one is at 
$T$-fixed points and the other is along characteristic suborbifolds.  
Let us explain the former first.  Suppose that torus orbifolds $M$
and $M'$ with $d(M)=d(M')$ have $T$-fixed points $q$ and $q'$ 
respectively such that 
the $n$-dimensional cones and the edge vectors corresponding to them
are the same and the signs $\epsilon_q$ and $\epsilon_{q'}$ at
$q$ and $q'$
are opposites.  Then there are a finite covering $\tilT$ of $T$,
a finite subgroup $H$ of $\tilT$ and orbifold charts $(U,V,H,p)$
and $(U',V,H,p')$ around $q$ and $q'$ respectively such that
$V$ is an invariant open disk centered at the origin in a
$\tilT$-module. In particular a diffeomorphism (in the sense
of orbifold) $f$ from the closure of $U$ onto that of
$U'$ is induced. Moreover $f$ sends characteristic suborbifolds
that contain $q$ onto characteristic suborbifolds that contain $q'$.
It should be noticed that $f$ is orientation reversing on $U$ and
on all the characteristic suborbifolds.  
We remove $U$ and $U'$ from $M$ and $M'$ respectively and glue
their boundaries through the diffeomorphism $f$ restricted 
to the boundaries.  The resulting space is a torus orbifold with 
the orientations compatible with the torus orbifolds $M$ and $M'$. 

Let us explain the equivariant connected sum along characteristic 
suborbifolds. For the sake of simplicity we assume that $d(M)=1$. 
Let $M_i$ be a characteristic suborbifold,
$p$ a point in the principal stratum of the orbifold $M_i$.
We may assume that the isotropy subgroup at $p$ of the $T$-action 
is the circle group $S_i$. Let $\tilSi$ be the covering group
of $S_i$ corresponding to the edge vector $v_i$ as introduced in
the previous section. Denote by $V_i$ 
the standard complex $1$-dimensional $\tilSi$-module 
and by $D(V_i)$ the unit disk of $V_i$. 
Then it follows from the Slice Theorem and Lemma \ref{lemm:isotropy}
that the $T$-orbit of $p$ has a closed invariant tubular neighborhood 
$\bar{U}_i$ in $M$ equivariantly diffeomorphic to
\begin{equation*}
(T\times_{\tilSi} D(V_i))\times D^{n-1} \tag{13.1}
\end{equation*}
where 
$T\times_{\tilSi} D(V_i)$ denotes the orbit space of $T\times D(V_i)$ by 
the $\tilSi$-action defined by $s(t,x)=(t\rho_i(s)^{-1},sx)$ for 
$s\in \tilSi, t\in T$ and $x\in D(V_i)$.  

Suppose that there are characteristic suborbifolds $M_i$ and $M'_{i'}$ 
of torus orbifolds $M$ and $M'$ with $d(M)=d(M')=1$ respectively
such that the corresponding edge vectors coincide.
Then the corresponding circle subgroups $\tilSi$ and
$\tilde{S}'_{i'}$ 
agree and there is an equivariant diffeomorpism between $\bar{U}_i$ 
and $\bar{U}'_{i'}$ 
reversing the orientations induced from $M$, $M_i$, $M'$ and $M'_{i'}$ 
because both $\bar{U}_i$ and $\bar{U}'_{i'}$ are equivariantly 
diffeomorphic to the space in (13.1) and 
$D^{n-1}$ ($n\ge 2$) has an orientation reversing self-diffeomorphism.  
We remove the interior of $\bar{U}_i$ and $\bar{U}'_{i'}$ from $M$ and $M'$
and paste them together along the boundaries of $\bar{U}_i$ and $\bar{U}'_{i'}$
through the orientation reversing equivariant diffeomorphism
restricted to the boundaries, 
producing a new torus orbifold, say $M''$.  We call this procedure 
the equivariant connected sum of $M$ and $M'$ along $M_i$ and $M'_{i'}$. 
The codimension of the principal orbits in $M_i$ and $M'_{i'}$ 
is $n-1$, so when $n\ge 3$, 
$M_i$ and $M'_{i'}$ are pasted together to become one characteristic 
suborbifold in $M''$ and $\Delta(M'')$ 
is obtained from $\Delta(M)$ and $\Delta(M')$ by identifying
$i$ with $i'$. 
However, when $n=2$, the characteristic suborbifolds 
$M_i$ and $M'_{i'}$ are $S^2$ 
and the principal orbits in them are circles; so the orbits separate 
$M_i$ and $M'_{i'}$ into 
two connected components respectively and hence two characteristic
suborbifolds of $M''$ are produced.  

Let $I\in\Sigma(M)^{(n)}$ and $I'\in\Sigma(M')^{(n)}$ be such that 
$C(M)(I)=C(M)(I')$. Suppose that the corresponding edge vectors are
the same for $I$ and $I'$. Then one can make equivariant connected sum 
of $M$ and $M'$ along each pair of characteristic suborbifolds $M_i$
and $M'_{i'}$ 
such that $C(M)(i)=C(M')(i')$ 
for $i\in I$ and $i'\in I'$, and then elements in $I$ and $I'$ will be 
identified in pairs in the multi-fan of 
the resulting torus orbifold and the weights $w^{\pm}$ on the
identified $n$-dimensional cone  is the sum of those at $I$ and $I'$. 

We say that $\Delta$ is \emph{connected} 
if $\Sigma$ is connected.  
According to the decomposition of $\Sigma$ into connected components, 
the multi-fan $\Delta$ decomposes 
into connected multi-fans 
which are again complete simplicial and of dimension $n$. 

\begin{lemm} \label{lemm:connected}
Suppose $n\ge 2$. Then 
the multi-fan $\Delta$ is geometrically realized if all connected 
components of $\Delta$ are geometrically realized.
\end{lemm}

\begin{proof} Let $M$ be a torus orbifold of dimension $2n$ 
and let $p$ be a point in the principal stratum of $M$.
We may suppose that $d(M)=1$. A closed tubular neighborhood $\bar{U}$ of the 
orbit of $p$ is equivariantly diffeomorphic to $T\times D^n$ 
and the complement of $\bar{U}$ is connected because $M$ is connected and the orbit 
has codimension $n\ge 2$.  

Let $M'$ be 
another torus orbifold of dimension $2n$ with $d(M')=1$, and
let $\bar{U}'$ be a closed subset 
in $M'$ corresponding to $\bar{U}$ in $M$.  
Since both $\bar{U}$ and $\bar{U}'$ are equivariantly 
diffeomorphic to $T\times D^n$ and $D^n$ has an orientation reversing 
diffeomorphim, there is an orientation reversing equivariant diffeomorphism 
between $\bar{U}$ and $\bar{U}'$.  We remove the interior of $\bar{U}$ and $\bar{U}'$ 
from $M$ and $M'$ respectively and glue their boundaries through the 
diffeomorphism restricted to the boundaries and 
obtain a new torus orbifold $M''$.  The multi-fan $\Delta(M'')$ is the 
disjoint union of $\Delta(M)$ and $\Delta(M')$. (Precisely speaking, 
$\Sigma(M'')$ is the disjoint union of $\Sigma(M)$ and $\Sigma(M')$ 
with the empty sets in them identified.)  

If all connected 
components of $\Delta$ are geometrically realized, then we connect 
torus orbifolds that realize the connected components of $\Delta$ by the 
above method. Then the resulting torus orbifold realizes $\Delta$. 
\end{proof}

As is shown in the proof of Lemma~\ref{lemm:connected}, whenever we have 
more than two torus orbifolds of dimension $n\ge 2$, we can connect 
them and the multi-fan of the resulting torus orbifold is the 
disjoint union of the multi-fans of the torus orbifolds we had. 

\begin{defi} 
We say that a complete simplicial multi-fan $\Delta=(\Sigma,\Lambda,w^\pm)$ 
of dimension $n$ is \emph{minimal} if 
\begin{enumerate}
\item $\Sigma$ is isomorphic to the argumented simplicial set obtained 
from the boundary of an $n$-simplex, and 
\item the set $\{w^+(I),w^-(I)\}$ is independent of $I\in\Sigma^{(n)}$. 
\end{enumerate}
\end{defi}

Although the set $\{w^+(I),w^-(I)\}$ is independent of $I$ for a minimal 
multi-fan $\Delta$, the pair $(w^+(I),w^-(I))$ may not be independent of 
$I\in\Sigma^{(n)}$.  
But one can convert $\Delta$ into another minimal multi-fan $\bar\Delta
=(\Sigma,\bar C,\bar w^\pm)$ such that the pair 
$(\bar w^+(I),\bar w^-(I))$ is independent of $I$.  
The definition of $\bar\Delta$ is as follows. 
Since $\Delta$ is of dimension $n$ and the cardinality of $\Sigma^{(1)}$ 
is $n+1$, there is a relation $\sum_{i\in\Sigma^{(1)}}b_iv_i=0$ among 
the edge vectors $v_i$ with non-zero real numbers $b_i$.  We then define 
\[
\bar C(i):=\begin{cases} C(i) \quad&\text{if $b_i>0$,}\\
-C(i) \quad&\text{if $b_i<0$,}
\end{cases}
\]
and define $\bar C(K)$ for $K\in\Sigma^{(m)}$ with $m\ge 2$ to be the cone 
spanned by $\bar C(k)$'s for $k\in K$.  We also define 
\[
(\bar w^+(I),\bar w^-(I)):=\begin{cases}
(w^+(I),w^-(I))\quad&\text{if $\sharp\{i\in I\mid b_i<0\}$ is even,}\\
(w^-(I),w^+(I))\quad&\text{if $\sharp\{i\in I\mid b_i<0\}$ is odd,}
\end{cases}
\]
for $I\in\Sigma^{(n)}$.  

\begin{lemm} \label{lemm:Delta_bar}
$\bar\Delta$ is minimal and satisfies the following two conditions:
\begin{enumerate}
\item the $n$-dimensional cones $\bar C(I)$ $(I\in\Sigma^{(n)})$ do not 
overlap and their union covers the entire space $N_\R$, and
\item the pair $(\bar w^+(I),\bar w^-(I))$ is independent of 
$I\in\Sigma^{(n)}$.
\end{enumerate}
Moreover $\Delta$ is geometrically realized if and only if so is $\bar\Delta$. 
\end{lemm}

\begin{proof} Let $\bar v_i$ be a non-zero vector in the cone $\bar C(i)$. 
One may choose it to be $v_i$ if $b_i>0$ and $-v_i$ if $b_i<0$. 
Then one has a relation $\sum_{i\in\Sigma^{(1)}}\bar b_i\bar v_i=0$ with 
positive numbers $\bar b_i$.  This implies the statement (1) in the lemma. 

We shall prove the statement (2) in the lemma. 
Let $J\in\Sigma^{(n-1)}$.  Since the cardinality of $\Sigma^{(1)}$ is $n+1$, 
there are exactly two elements $i, i'\in\Sigma^{(1)}$ not contained in 
$J$, and $J\cup\{i\}$ and $J\cup\{i'\}$ are in $\Sigma^{(n)}$, in other words, 
the $(n-1)$-dimensional cone $C(J)$ is a facet of only two $n$-dimensional 
cones $C(J\cup\{i\})$ and $C(J\cup\{i'\})$.  We project them on 
$N_\R^{C(J)}$ (the quotient space 
of $N_\R$ by the subspace generated by $C(J)$). 
Then the vectors projected from $v_i$ and 
$v_{i'}$ are toward opposite directions if and only if $b_ib_{i'}>0$. 
It follows from the completeness of $\Delta$ that 
$w(J\cup\{i\})=\sign(b_ib_{i'})w(J\cup\{i'\})$.  This together with the 
definition of $\bar w^\pm$ shows that 
$\bar w(J\cup\{i\})=\bar w(J\cup\{i'\})$.  Since $J\in\Sigma^{(n-1)}$ is 
arbitrary, this proves the statement (2).  It also proves the 
completeness of $\bar\Delta$, so that $\bar\Delta$ is minimal. 

The procedure from $\Delta$ to $\bar\Delta$ corresponds 
to reversing orientations on characteristic suborbifolds $M_i$ with $b_i<0$, 
so the latter statement in the lemma is obvious.  
\end{proof}

\begin{lemm} \label{lemm:primitive}
Let $\Delta$ be a minimal multi-fan of dimension $n\ge 2$. 
If $n\ge 3$, then $\Delta$ is geometrically realized.  
If $n=2$, then $\Delta$ is geometrically realized if (and only if) 
$\{w^+(I),w^-(I)\}=\{1,0\}$ for any $I\in\Sigma^{(2)}$. In any case
we can take an orbifold structure on the realizing torus orbifold
such that the corresponding edge vectors $\{v_i\}$ are all primitive;
that is, if $v_i=a_iv'_i$ for some $v'_i\in N$ and $a_i\in \Z$, then
$a_i=\pm 1$. 
\end{lemm}

\begin{proof} 
By Lemma~\ref{lemm:Delta_bar}, we may assume that the union 
of cones $C(I)$ over $I\in\Sigma^{(n)}$ covers the entire space $N_\R$ 
and the pair $(w^+(I),w^-(I))$, which we denote by 
$(p,q)$, is independent of $I$.  When $(p,q)
=(1,0)$, $\Delta$ can be realized by a weighted projective space,
say $X$. There is an orbifold structure on a weighted projective
space such that the edge vectors are all primitive. We admit
these facts for a moment; the proof will be give in the appendix at the end
of this section. Then 
$-X$ realizes the case when $(p,q)=(0,1)$. 
This completes the proof when $n=2$. 

Suppose $n\ge 3$.  For a general value of 
$(p,q)$, we prepare $p$ copies of $X$ and $q$ copies of $-X$ and 
do equivariant connected sum along all $X_i$'s and $-X_i$'s 
for each $i\in\Sigma(X)$.  Then the resulting torus orbifold realizes $\Delta$.
The edge vectors are all primitive in this construction since
it is so for $X$.
\end{proof}

Now let $\Delta$ be an arbitrary complete simplicial multi-fan of 
dimension $n\ge 2$.  We decompose $\Delta$ into a number of 
minimal multi-fans as follows. 
We choose and fix a generic (rational) $1$-dimensional cone in $N_\R$, 
say $\ell$, which is not contained in any subspaces spanned by 
cones of dimension 
$\le n-1$ in $\Delta$.  We label $\ell$ as $\star$. 
To each $n$-dimensional cone $C(I)$ for $I\in\Sigma^{(n)}$, 
we form $n$ cones which are respectively spanned by 
$\ell$ and facets of $C(I)$.  These $n$ cones together 
with $C(I)$ determine a simplicial multi-fan 
$\Delta[I]=(\Sigma[I], \Lambda[I], w[I]^\pm)$, where $\Sigma[I]$ 
consists of all proper subsets of $I\cup\{\star\}$. 
The weight functions $w[I]^\pm$ are defined as follows.  
Let $v_i$ be a non-zero vector in $C(i)$ for each $i\in I$ and
$v_\star$ a non-zero vector in $\ell$.  Then there is a relation 
\begin{equation*}
v_\star+\sum_{i\in I}a_iv_i=0 \tag{13.2}
\end{equation*}
with non-zero real numbers $a_i$'s.  
Let $\mathcal I\in 
\Sigma[I]^{(n)}$.  Then $\mathcal I=I$ or $(I\backslash\{i\})\cup\{\star\}$ 
for $i\in I$.  We define 
\begin{equation*}
(w[I]^+(\mathcal I),w[I]^-(\mathcal I)):=\begin{cases} 
(w^+(I),w^-(I))\quad&\text{if $\mathcal I=I$ or }\\
&\text{ $\mathcal I=(I\backslash\{i\})\cup\{\star\}$ and $a_i>0$,}\\
(w^-(I),w^+(I))\quad&\text{if $\mathcal I=(I\backslash\{i\})\cup\{\star\}$ 
and $a_i<0$.}
\end{cases}\tag{13.3}
\end{equation*} 

\begin{lemm} 
$\Delta[I]$ is complete and hence minimal.  
\end{lemm}

\begin{proof} 
The proof is essentially the same as that of lemma~\ref{lemm:Delta_bar}. 
As remarked in Section 2, it suffices to show that, when a generic vector $v$ 
gets across an $(n-1)$-dimensional cone, the integer $d_v$ in Section 2 
remains unchanged. 
Let $\mathcal J$ be an element of $\Sigma[I]^{(n-1)}$ and let 
$i$ and $i'$ be 
the two elements in $(I\cup\{\star\})\backslash \mathcal J$.  Then 
$\mathcal I:=
\mathcal J\cup\{i\}$ and $\mathcal I':=\mathcal J\cup\{i'\}$ 
are the elements in 
$\Sigma[I]^{(n)}$ which contain $\mathcal J$.  
We project cones $C[I](\mathcal I)$ and $C[I](\mathcal I')$ on 
$N_\R^{C[I](\mathcal J)}$.  Then it follows from (13.2) that the vectors 
projected from $v_i$ and $v_{i'}$ are toward opposite directions if and 
only if $a_{i}a_{i'}>0$, 
where $a_\star$ is understood to be $1$.  This together with the definition 
(13.3) of $w[I]^\pm$ 
implies that $d_v$ remains unchanged regardless of the sign of 
$a_{i}a_{i'}$ when $v$ gets across the $(n-1)$-dimensional cone 
$C[I](\mathcal J)$. 
\end{proof}

Let $J\in\Sigma^{(n-1)}$ and let $I_1, \dots, I_r$ be the elements in 
$\Sigma^{(n)}$ containing $J$.  
The $n$-dimensional cone spanned by $C(J)$ and $\ell$ 
appears in $\Delta[I_k]$ for $k=1, 2\dots, r$ with 
the form $C[I_k](J\cup\{\star\})$.

\begin{lemm} \label{lemm:sum w[I_k]}
$\sum_{k=1}^r w[I_k](J\cup\{\star\})=0.$
\end{lemm}

\begin{proof}
Consider the projection of the cones $C(I_k)$'s 
on $N_\R^{C(J)}$. 
We define $\sign(I_k)=1$ or $-1$ according as the projection image of 
$C(I_k)$ disagrees or agrees with that of $\ell$. 
Applying (13.3) with $I=I_k$ and $I\backslash\{i\}=J$, one sees that 
\[
w[I_k](J\cup\{\star\})=\sign(I_k)w(I_k).
\]
On the other hand, it follows from the completeness of $\Delta$ that 
\[
\sum_{\sign(I_s)=1}w(I_s)=\sum_{\sign(I_t)=-1}w(I_t).
\]
These two equalities imply the lemma. 
\end{proof}

\begin{proof}[Proof of Theorem~\ref{theo:higher dim torus orbifold}] 
By lemma~\ref{lemm:connected} we may assume that $\Delta$ is connected. 
We choose a generic (rational) $1$-dimensional cone $\ell$ and decompose 
$\Delta$ using $\ell$ into minimal multi-fans $\Delta[I]$'s 
$(I\in\Sigma^{(n)})$.  
By Lemma~\ref{lemm:primitive} $\Delta[I]$ is realized by a torus orbifold, 
say $M[I]$, such that all its edge vectors are primitive.  
We consider the disjoint union of $M[I]$ over $I\in\Sigma^{(n)}$ 
and piece them together using equivariant connected sum in 
the following way. For each $i\in \Sigma^{(1)}$ we do equivariant
connected sum of $\{M[I]\mid i\in I\}$ successively along $M[I]_i$'s, 
and similarly do equivariant connected sum of all $M[I]$'s 
along $M[I]_\star$ as well. 
The resulting space is connected because $\Delta$ is connected, and 
becomes a torus orbifold.  
Its multi-fan is close to $\Delta$ but contains extra 
cones which are the cones spanned by $\ell$ and 
$C(J)$ for $J\in\Sigma^{(m)}$ with $m\le n-1$.  
For a fixed $J\in \Sigma^{(n-1)}$, 
it follows from Lemma~\ref{lemm:sum w[I_k]} that 
there are the same number of $T$-fixed points $p$ with $\epsilon_p=1$
and $q$ with $\epsilon_q=-1$ contained
in the union of $M[I_k]$ with $J\subset I_k$ and corresponding to 
the cone spanned by $\ell$ and $C(J)$. Hence one can do 
equivariant connected sum at pairs of $T$-fixed points $p$ and $q$
so that those 
$T$-fixed points will be eliminated.  
Doing this for each $J\in\Sigma^{(n-1)}$, we obtain 
a torus orbifold, say $M$, realizing $\Delta$. 
In fact, the characteristic suborbifolds $M[I]_\star$ turn into 
a codimension two suborbifold of $M$, which is 
fixed by the circle subgroup determined 
by $\ell$ but has no $T$-fixed point, so it is not a characteristic 
suborbifold of $M$ by definition.  This means that all the cones 
in $\Delta[I]$'s containing $\ell$ as an edge do not show up in the 
mulit-fan of $M$.  
\end{proof}

\begin{proof}[Proof of Theorem~\ref{theo:4-dim torus orbifold}] 
We already observed the \lq\lq only if" part, so we prove the 
\lq\lq if" part. 
By Lemma~\ref{lemm:connected} we may assume that our 
$\Delta$, which satisfies the conditions (1) and (2) 
in Theorem~\ref{theo:4-dim torus orbifold}, 
is connected.  Then (the realization of) $\Sigma$ is either 

\emph{Case 1.} a $1$-simplex, or 

\emph{Case 2.} the boundary of a $d$-gon where $d\ge 3$, 

\noindent
and that 
\begin{equation*}
\{w^+(I),w^-(I)\}=\begin{cases} \{1,1\} &\text{ in Case 1,}\\
\{1,0\} &\text{ in Case 2.}
\end{cases}
\end{equation*}
Using the latter statement in Lemma~\ref{lemm:primitive}, 
the same argument as in the proof of 
Theorem~\ref{theo:higher dim torus orbifold} shows that $\Delta$ in Case 2 
is geometrically realized.  As for Case 1, 
let $I\in\Sigma^{(2)}$ be the unique simplex. There exist a finite
covering $\tilT\to T$ whose kernel $H$ is isomorphic to  
$N/N_I$ where $N_I$ is the 
sublattice generated by the primitive vectors $v_i$'s for $i\in I$,
and a $2$-dimensional $\tilT$-module $V$ corresponding to 
the cone $C(I)$, as was explained in Section \ref{sec:orbifolds}. 
Then the one point compactification of $V/H$, i.e.,
the orbit space of $S^4$ by an action of $N/N_I$, 
realizes our $\Delta$ in Case 1.
\end{proof}
\vspace*{0.5cm}
\noindent{\it Appendix.\enskip Realization of minimal multi-fans 
by weighted projective spaces}.
\par
\vspace*{0.3cm}
We identify the $(n+1)$-dimensional torus 
$T^{n+1}=S^1\times \cdots \times S^1$ 
with the standard maximal torus of $GL(n+1,\C)$ 
consisting of diagonal matrices. We set $\tilT =T^{n+1}/D$ 
where $D$ denotes the subgroup of diagonal elements $(z,\ldots,z)$. 
It is a maximal torus in $PGL(n+1,\C)$ and acts effectively on 
the projective space $\Pn$. 
Let $\tilSi$ denote the $i$-th factor of $T^{n+1}$. It is mapped
injectively into $\tilT$. We shall denote by the same letter 
$\tilSi$ its
image in $\tilT$. We set $\tilde{M}_i=\{[z_0,\dots ,z_n]\mid z_i=0\}$,
for $i=0,\dots ,n$. They are the characteristic submanifolds of
$\Pn$ regarded as a torus manifold with the orientations induced
from the complex structure. If $H$ is a finite subgroup of
$\tilT$, then the quotient $M_H=\Pn /H$ is a torus orbifold 
acted on by $T=\tilT/H$ for which
$(M_H,\Pn,H,p)$ is an orbifold chart, where $p:\Pn \to M_H$ is the 
projection. It is called a weighted projective space.
Its characteristic suborbifolds are 
$M_i=p(\tilde{M}_i),\ i=0,\dots ,n$, and the corresponding circle
subgroups are $S_i=\pi(\tilSi)$, where $\pi :\tilT \to T$ is
the projection. The symmetric group 
$\mathcal{S}_{n+1}$ of degree $n+1$ acts on $T^{n+1}$ and also
induces an action on $\tilT$. It also acts on $\Pn$. If $H^{\sigma}$
denotes the transform of $H$ by an element $\sigma\in \mathcal{S}_{n+1}$,
then the transformation $\sigma :\Pn \to \Pn$ induces an isomorphism
of torus manifolds $M_H\to M_{H^{\sigma}}$. We set
\[ \mathcal{WP} =
 \{H\mid \text{finite subgroup of}\ \tilT\}/\mathcal{S}_{n+1}. \]
Every element in $\mathcal{WP}$ represents an isomorphism class
of weighted projective spaces. 

In order to describe the multi-fan $\Delta_H$ associated with
the torus orbifold $M_H$ we introduce the following notations:
\[ \tilN=\Z^{n+1}/\text{diagonal submodule},\quad
\tilv_i=\text{image of $e_i$ in}\ \tilN,\quad N=\Z^n, \]
where $e_i$ is the $i$-th fundamental unit vector in $\Z^{n+1}$.
$\tilN$ is canonically identified with $\Hom(S^1,\tilT)$.
If one chooses an identification of $\Hom(S^1,T)=H_2(BT;\Z)$ with $N$,
then the finite covering map $\pi:\tilT \to T$ induces an injective
homomorphism $\varphi:\tilN\to N$. The vectors $v_i=\varphi(\tilv_i)$
are the edge vectors of the $1$-dimensional cones of $\Delta_H$.
Note that they satisfy the equality
\begin{equation*}
  \sum_iv_i =0, \tag{13.4}
\end{equation*}
since the $\tilv_i$'s satisfy a similar equality. This implies
that $\Delta_H$ is a minimal multi-fan satisfying the conditon
(1) in Lemma \ref{lemm:Delta_bar}. It is also clear that
$(w^+(I),w^-(I))=(1,0)$. We shall denote by $\mathcal{MF}$ the set 
of minimal multi-fans satisfying the above two conditions.
If one chooses another identification of $\Hom(S^1,T)$ with
$N$, then $\varphi$ is transformed to $\psi\circ\varphi$ where
$\psi\in GL(n,\Z)$. $GL(n,\Z)$ acts on $\mathcal{MF}$ from left
by transforming the cones simultaneously by its elements.
Let $d_H\in \Z$ be the maximal common divisor of the edge vectors 
$v_i$ of $\Delta_H$. We get a correspondence
\[ \alpha:\mathcal{WP}/\mathcal{S}_{n+1} \to 
 GL(n,\Z)\backslash\mathcal{MF}\times \Z_{>0} \]
which sends $H$ to $(\Delta_H,d_H)$. 

\begin{lemm}\label{lemm:WP}
The correspondence $\alpha$ is a bijection. In particular, every 
minimal multi-fan $\Delta$ in $\mathcal{MF}$ is realizable.
\end{lemm}

\begin{proof}
We shall define a correspondence $\beta :
GL(n,\Z)\backslash\mathcal{MF}\times \Z_{>0}
\to \mathcal{WP}/\mathcal{S}_{n+1}$ which is to be the inverse of
$\alpha$. Take a multi-fan $\Delta$ in $\mathcal{MF}$ and 
$d\in \Z_{>0}$. It is easy to see there is a unique set $\{v_i\}$
of edge vectors of $\Delta$ such that $\sum_iv_i=0$ and the maximal
common divisor of $\{v_i\}$ is $d$. Define a homomorphism
$\varphi:\tilN\to N$ by requiring $\varphi(\tilv_i)=v_i$.
Then there is a unique finite covering map $\pi:\tilT\to T$ which 
induces $\varphi:\tilN=\Hom(S^1,\tilT)\to N=\Hom(S^1,T)$. Let
$H$ be the kernel of $\pi$. The homomorphism $\varphi$, hence 
$H$ either, does not
depend on the choice of identification $N=\Hom(S^1,T)$, but it 
depends on the numbering of $v_i$'s. So if we put $\beta(\Delta,d)
=\text{the class of $H$ in}\ \mathcal{WP}/\mathcal{S}_{n+1}$, 
it induces a correspondence $\beta$ as above. It is clear that
$\beta$ is in fact the inverse of $\alpha$.
\end{proof}

\begin{rema}
Let $a$ be a positive integer. The correspondence
$T^{n+1}\ni(z_0,z_1,\dots ,z_n)\mapsto(z_0^a,z_1^a,\dots ,z_n^a)
 \in T^{n+1}$ induces a homomorphism $\rho:\tilT\to \tilT$.
For a finite group $H$ of $\tilT$ define $H'=\rho^{-1}(H)$. The
edge vectors $\{v'_i\}$ corresponding to the torus manifold $M_{H'}$
are of the form $v'_i=av_i$, where $\{v_i\}$ correspond to $M_H$. 
Hence $\Delta_H=\Delta_{H'}$ and $d_{H'}=ad_H$. Let $g:\Pn \to \Pn$
be the map defined by 
$g[z_0,z_1,\dots ,z_n]=[z_0^a,z_1^a,\dots ,z_n^a]$.
Then it induces a homeomorphism $M_{H'}\to M_H$ which is 
equivariant with respect to the isomorphism of tori between 
$\tilT/H'$ and $\tilT/H$ induced by $\rho$.
If $M_H$ and $M_{H'}$ are considered as algebraic varieties 
then the homeomorphism becomes an equivalence. It is a fundamental
fact in the theory of toric varieties that to each fan corresponds 
a toric variety. The above equivalence gives an 
interpretation of this fact within this special case in our context. 
Related results are found in \cite{LT}. 
Related to the above remark, for a later use, we point out the following fact. 
Let $a_i,\dots ,a_n$ be positive integers, and let 
$\Z/a_i\subset S^1$ be the subgroup of $a_i$-th roots of unity. Set
$G=\prod_i\Z/a_i$. Then the map 
$\C^n\ni(z_1,\dots ,z_n)\mapsto (z_1^{a_1}\dots ,z_n^{a_n})\in\C^n$ 
induces an equivalence of affine algebraic varieties 
$\C^n/G\to \C^n$.
\end{rema} 

Let $M_H\in \mathcal{WP}$ and let $\{v_i\}$ be the edge vectors
corresponding to the orbifold structure as given above. Even if 
$d_H=1$, it may happen
that some of $v_i$'s are not primitive. We will show that there
always exists a torus orbifold structure on $M_H$ such that 
the corresponding edge vectors are all primitive. More generally
we have

\begin{lemm}
Let $M_H$ be a weighted projective space and $\{v_i\}$ the 
corresponding edge 
vectors satisfying $\sum_iv_i=0$ as given above. Suppose that 
$\{v'_i\}$ are vectors in $N$ such that $v_i=a_iv'_i$ with 
$a_i\in \Z_{>0}$. Then there 
is an orbifold structure on $M_H$ which admits $\{v'_i\}$ as the
corresponding edge vectors.
\end{lemm}

\begin{proof}
For each $x\in M_H$ let $\tilT_x\subset \tilT$ 
be the isotropy subgroup at $\tilde{x}$ of the $\tilT$-action
on $\Pn$ where $\tilde{x}\in p^{-1}(x)$. $\tilT_x$ does not depend
on the choice of $\tilde{x}$ in $p^{-1}(x)$. If $x$ lies in
$\Int M_I=M_I\setminus \bigcup_{J \supsetneqq I}M_J$ for 
$I\in \Sigma(M_H)^{(k)}$, then $\tilT_x =\tilde{S}_I=\prod_{i\in I}\tilSi$.
We put $H_x=H\cap \tilT_x$. We take a family 
$\{V_{x,\mu}|\mu \in \Z_{>0}\}$ of small $\tilT_x$-invariant
open neighborhoods of $\tilx$ such that $V_{x,\mu}$ converges to $\tilx$ 
when $\mu$ tends to infinity. We may assume that $V_{x,\mu}$ is
equivariantly diffeomorphic to an $\tilde{S}_I$-invariant open disk
in $\C^n$. It is possible to make $V_{x,\mu}$'s so
small that they satisfy the following condition:
\begin{equation*}
 H_x= \{h\in H\mid h\cdot V_{x,\mu}\cap V_{x,\mu}\not= \emptyset\}.\tag{13.5}
\end{equation*} 
Then $U_{x,\mu}=V_{x,\mu}/H_x$ is an open 
neighborhood of $x$ in $M_H$, and $(U_{x,\mu},V_{x,\mu},H_x,p|V_{x,\mu})$ 
is an orbifold chart of $M_H$ compatible with $(M_H,\Pn,H,p)$.

On the other hand the fact that $v_i=a_iv'_i$ implies that the
kernel of $p:\tilSi\to \Si$ contains $\Z/a_i$,which we denote
by $G_i$. Since $H$ is 
the kernel of $p:\tilT\to T$, $G_i$ is contained in $H$. We
put $G_I=\prod_{i\in I}G_i$ for $I\in \Sigma(M_H)^{(k)}$ and define
\[ V'_{x,\mu}=V_{x,\mu}/G_I,\enskip H'_x=H_x/G_I \enskip 
 \text{for}\enskip x\in \Int M_I .\]
$V'_{x,\mu}$ can be considered as an open disk in $\C^n$ as pointed
out in Remark above. The projection
$p|V_x:V_x\to U_x$ induces a map $p'_{x,\mu}:V'_{x,\mu}\to U_x$ 
which induces a homeomorphism $V'_{x,\mu}/H'_x \to U_x$.

We shall prove that the family 
$  \{(U'_{x,\mu},V'_{x,\mu},H'_x,p'_{x,\mu})\mid x\in M,\ \mu \in \Z_{>0}\} $
forms a set of orbifold charts of an orbifold structure on $M_H$.
For that purpose it suffices to show that, 
if $U'_{x.\mu}\subset U'_{y,\nu}$, then there are an injective homomorphism
$\rho:H'_x\to H'_y$ and a $\rho$-equivariant open embedding 
$\phi:V'_{x,\mu}\to V'_{y,\nu}$ such that
\begin{equation*}
 \rho(H'_x)=\{h\in H'_y\mid h\cdot \phi(V'_{x,\mu})\cap 
 \phi(V'_{x,\mu})\not= \emptyset\}. \tag{13.6}
\end{equation*}
The condition (13.5) implies that, if $x\in \Int M_I$ and
$y\in \Int M_J$ with $I \ \text{and}\ J \in \Sigma(M_H)$, and 
if $U'_{x.\mu}\subset U'_{y,\nu}$, then $I\supset J$. Therefore
\[ H_x\subset H_y \quad \text{and}\quad G_I\cap H_y=G_J. \]
It follows that the inclusion $H_x \to H_y$ induces an injective
homomorphism $\rho:H'_x=H_x/G_I \to H_y/G_J=H'_y$. If
$\tilx$ is taken in $V_{y,\nu}$, then $V_{x,\mu}$ is contained in
$V_{y,\nu}$. The inclusion induces an embedding 
$\phi:V'_{x,\mu}\to V'_{y,\nu}$. $\phi$ is clearly $\rho$-equivariant.
The condition (13.6) follows from (13.5).

If $x$ lies in $M_i$, then the action of $S_i$ lifts to the action
of $\tilde{S}'_i=\tilSi/G_i$ on $V'_{x,\mu}$ and the lifting is minimal. Hence
the edge vector of $C(i)$ corresponding to the orbifold structure
defined above must be $v'_i$.
\end{proof}

\bigskip
\bigskip

\providecommand{\bysame}{\leavevmode\hbox to3em{\hrulefill}\thinspace}

\end{document}